Department of Mathematics
TUM School of Computation, Information and Technology

Technical University of Munich

Master’s Thesis in Mathematics

# Fibrations in Directed Type Theory

Benno Lossin

I hereby declare that this thesis is entirely the result of my own work except where otherwise indicated. I have only used the resources given in the list of references.

Munich, 1. December 2025 Benno Lossin

**Abstract**

We study $\infty$-categories in the synthetic *simplicial type theory* developed by Riehl and Shulman. In particular, we define *cocartesian fibrations* and prove their closure properties using a novel equivalence between *LARI adjunctions* and *initial sections*. We formalize our work using the experimental proof assistant `rzk` and upstream our work to the formalization effort by Riehl et al. In addition to our new work, we also give an introduction to general type theory, homotopy type theory, and the simplicial type theory used by the rest of the thesis.

## Acknowledgments

I would like to thank my advisor Tashi Walde, who expended considerable effort and time discussing the material of this thesis. He helped me understand and dive into the theory without the fear of getting lost. I also would like to thank my supervisor Claudia Scheimbauer for introducing me to Tashi and for supervising this thesis.

## Contents

# 1 Introduction

Modern type theory collects several exceptional concepts and ideas into one foundational mathematical theory replacing logic and set theory. It is connected to many fundamental and very important disciplines of mathematics: higher category theory, homotopy theory, order theory, logic, foundational mathematics, formalized mathematics and constructive mathematics. It is rooted in foundational mathematics, logic and theoretical computer science and offers a high degree of correctness through its good compatibility with formalization. In addition, it offers one of the most intuitive interpretations of mathematics. Connections to homotopy theory and higher categories have only recently been discovered, which has led to many interesting explorations of the novel branch of *homotopy type theory* (HoTT) culminating in the work from The Univalent Foundations Program [1].

This thesis summarizes the *simplicial type theory* developed by E. Riehl and M. Shulman [2] as an extension of homotopy type theory. E. Riehl *et al.* [3] have formalized this theory using the proof assistant `rzk` created by N. Kudasov, A. Abounegm, and D. Danko [4] based on Riehl and Shulman's theory. Several extensions of the theory by U. Buchholtz and J. Weinberger [5] and C. Bardomiano Martínez [6] have also been partially formalized by E. Riehl *et al.* [3]. We continue their work by formalizing *cocartesian fibration* in this thesis and contributing it upstream at [3]. We in particular follow the work from U. Buchholtz and J. Weinberger [5] and J. Weinberger [7] who study cocartesian fibrations as *LARI-families*. We connect these to *initial sections* by showing that initial sections and *LARI adjunctions* are equivalent, which simplifies the formalization effort. We also prove that cocartesian fibrations are closed under products, pullbacks, and composition. Our formalizations can be found online at https://github.com/BennoLossin/sHoTT and when they have been accepted upstream, they will also be available at https://rzk-lang.github.io/sHoTT in a more readable form.

In the rest of this section, we give a brief overview of the most important historic developments in type theory. Then we cover the connections to the other branches of mathematics and discuss the merits of the synthetic approach. We also discuss the two different ways of adding directionality into homotopy type theory, in order to study higher categories. Finally we give an overview of the structure of this thesis.

## 1.1 History

Type theory originated from foundational mathematics and the solution to Russell's paradox in Principia Mathematica by A. N. Whitehead and B. Russel [8] in 1910. Later, A. Church [9] developed Lambda Calculus in 1940 as a "mathematical Turing machine" and introduced types to ensure termination. At this point the notion of a type already has a very similar meaning to today. Another gigantic leap in the theory was made by P. Martin-Löf *et al.* [10] in 1984 by introducing intuitionistic type theory. His impactful work is the basis of many type theories to this day and commonly called Martin-Löf-Type-Theory (MLTT). This theory leans into constructive mathematics and is a proper foundational system. MLTT was further developed by many logicians, mathematicians and theoretical computer scientists. The next milestone we focus on is the connection to homotopy theory and the formation of the new branch of *homotopy type theory* allowing the study of synthetic $\infty$-groupoids. V. Voevodsky [11] contributed the idea of *univalence* in 2010, which permits one to treat *weakly equivalent* objects as *equal*. This ultimately culminated in the publishing of "The HoTT Book" by The Univalent Foundations Program [1] in 2013. It gives a comprehensive overview of homotopy type theory along with deriving several mathematical results in it. Many further extensions of this theory have been explored. In particular, E. Riehl and M. Shulman [2] added directionality through a strict directed interval type in 2017. This allows the theory to define $(\infty, 1)$-categories in a synthetic manner.

We base the core of our work on two papers: U. Buchholtz and J. Weinberger [5] (2021) and J. Weinberger [7] (2024). They considered the notion of cocartesian fibration from several different directions in Riehl and Shulman's theory.

## 1.2 Connections to other Disciplines

Mathematics is a syntactical game of writing down the *correct* thing. The rules of the theory inform us about which syntactic constructions are valid. Type theory not only describes which syntactic constructions (also called *terms*) are valid, but it associates to each term a *type*. A type also is a syntactic construction, but on a different level than a term. This type describes what the capabilities of the value represented by the term are. For example, the type of a function $f : A \to B$ informs us about which values we can plug in and which values we can potentially get out.

Mathematicians often talk about the capabilities of values, for example "this function is a weak equivalence" or "this object is initial in the category of XYZ". In type theory these statements simply are that the value (or a derived value) is of a certain type. When mathematicians argue, they often recall and make use of these properties, for example "since this function is an equivalence" or "since this object is initial in XYZ". This corresponds to *using* the value as is dictated by the type.

Type theory naturally formalizes this notion of "capability" and makes working with them the primary way of doing mathematics. The central dogma of type theory is:

*"Types are propositions and terms are proofs."*

This theme is what connects most of the other disciplines to type theory. Giving a proof for a proposition boils down to writing down a term of the correct type. The ability of keeping track of *different* elements of types is what connects type theory to constructive and proof-relevant mathematics. A priori, elements of a type can be different, permitting multiple distinct proofs for the same statement. This gives us the connection to homotopy theory and higher categories, as we can track the additional structure using the terms and their non-trivial equalities are naturally found in the theory.

### 1.2.1 Foundational Mathematics and Logic

Type theory is special compared to many foundational systems before it, because it is not formulated in a logic. Type theory itself naturally embeds logical constructions and is directly formulated in the language logicians use to define logics. The Curry-Howard meta-theoretic isomorphism describes how logic and type theory relate. Table 1 shows a brief overview of how type theory incorporates logic. We discuss this in more detail in Section 2.

| **Concept in Logic** | **Concept in Type Theory** |
|---|---|
| Propositions | Types |
| Proof of a proposition $P$ | Term of type $P$ |
| True $\top$ | Unit type $\mathbb{1}$ |
| False $\bot$ | Empty type $\mathbb{0}$ |
| Conjunction $A \land B$ | Product type $A \times B$ |
| Disjunction $A \lor B$ | Disjoint union type $A + B$ |
| Implication $A \to B$ | Function type $A \to B$ |
| Existential quantifier $\exists x.P$ | $\Sigma$-type $\Sigma(x : X, P)$ |
| Universal quantifier $\forall x.P$ | $\Pi$-type $(x : X) \to P$ |

Table 1: Curry-Howard's meta-theoretic isomorphism between logic and type theory

Type theory is better suited as a foundational system for higher categories compared to set theory. This is because a cumulative hierarchy of bigger and bigger universes is built-in and equality is not a simple notion of "$x$ is equal to $y$", but also includes *how* the two elements are equal. Additionally, set theory has no built-in mechanism to track higher structures whereas type theory does.

We note that *type theory* is not a single theory by itself, but rather a collection of principles that one can use to build a theory. It is similar to how there is not a singular *logic* in mathematics, but different variants.

### 1.2.2 Formalized Mathematics

Type theory employs the language of logicians in a very specific way. The rules of type theory are written down in a way where each rule can be applied mechanically, without requiring any mathematical ideas or choices. Called *type checking*, the rules only dictate which terms are of which type and which terms are invalid. This makes them perfect for computers to reason about and develop programs that automatically check the correctness of statements. Formalized mathematics increase the trust we can put in our statements being correct. Additionally, the type theoretic approach puts formalized mathematics front and center as opposed to classical mathematics where formalization often is an afterthought.

Most proof assistants are built on top of Martin-Löf-Type-Theory developed by P. Martin-Löf *et al.* [10], which we already mentioned in Section 1.1. Examples of proof assistants are *Lean* (L. de Moura, S. Kong, F. Doorn, and J. Raumer [12], L. d. Moura and S. Ullrich [13] and Lean Contributors [14]), *Rocq* (T. Coquand and G. Huet [15] and The Coq Development Team [16]) and *Isabelle/HOL* (T. Nipkow and M. Wenzel [17] and Isabelle Contributors [18]).

There is an experimental proof assistant for the simplicial type theory by E. Riehl and M. Shulman [2]. It is called `rzk` and was created by N. Kudasov, A. Abounegm, and D. Danko [4]. We use it to formalize the majority of this thesis based on the work by E. Riehl *et al.* [3].

Proof assistants are a very useful tool for developing a good theory of higher categories. This is because higher categories are a very complex topic and one often needs to keep a lot of data in one's working memory. Additionally, they can in certain situations be unintuitive compared to ordinary categories. Proof assistants help to avoid these pitfalls and strengthen our trust in the theory.

### 1.2.3 Computer Science

Without a doubt one of the strongest connections from type theory is that to computer science. The study of programming languages – which employ various type theories – is a whole discipline of computer science and we cannot do it justice in this short paragraph. Functional programming languages are based on Church's Lambda Calculus and many modern languages borrow concepts from it. There are also many experimental programming languages using esoteric type theories. Computation and proving statements are inexplicably linked through the Curry-Howard isomorphism.

### 1.2.4 Constructive Mathematics and Proof-Relevance

Constructive mathematics concerns itself with using a restricted foundational system that does not permit conjuring values out of thin air. They for example do not assume the axiom of choice. More surprisingly, they also do not assume the axiom of excluded middle which states $A \vee \neg A$. In classical mathematics this is an indispensable tool used commonly for proofs by contradiction. The reason for dismissing this axiom from a constructive perspective is that from the axiom $A \vee \neg A$ one cannot derive *which* one of the two is true. This results in proofs by contradiction not being *constructive.* If one proves the existence of an object $x$ by contradiction, one has done so by proving that the non-existence of it leads to a contradiction and thus without ever constructing the $x$.

**Example 1.2.1**

A simple example exhibiting this behavior is the statement:

*Are there two irrational numbers $x$ and $y$ such that $x^y$ is rational?*

The answer is yes and we can show this easily using LEM (the law of excluded middle). We know that $\sqrt{2}$ is irrational. Now consider $x = y = \sqrt{2}$ and the following argument. Either $\sqrt{2}^{\sqrt{2}}$ is irrational or it is rational. If it is rational, we are done. Otherwise, consider $x = \sqrt{2}^{\sqrt{2}}$ and $y = \sqrt{2}$, which then satisfies the statement, since:

$$x^y = \left(\sqrt{2}^{\sqrt{2}}\right)^{\sqrt{2}} = \sqrt{2}^2 = 2$$

From this argument it is not clear which of the two cases hold, only that one of them holds. Thus such an argument is rejected by constructivists.

Aside from the philosophical reasons for adopting a constructive stance, one of the pragmatic reasons for it is that any proof can be turned into an algorithm producing the result of the theorem. This is particularly useful in the area of computer science and when applying mathematics to the real world.

Type theory invites constructivism through the principle of "types as propositions". We record *how* a proposition is proven through the term that we have to create in order to show inhabitance of a type. The $\Sigma$-type (the type theoretic analogue to the existential quantifier $\exists x.\ P(x)$) is proved to be inhabited by constructing a term of said $\Sigma$-type. This construction requires two elements: the first is the witness $x$ and the second the proof that the witness satisfies the condition $P(x)$. So if one is given a proof of an existential quantifier, one can always obtain a *concrete* element for which the condition holds. Similarly for conjunction, we require proof elements of both input statements to be provided to prove their conjunction.

This makes our theory not only constructive, but also *proof relevant.* We care about how a proposition can be proven. There might be multiple *different* ways of doing so. This is also the core of homotopy theory that studies *how* two things are equal and, if they are, how their equalities relate to each other.

We note that our theory still permits classical mathematics to take place, by considering only types which we call *propositions.* These types have at most one element, which restores proof irrelevance for *only* these types. We later note that we are allowed to assume LEM for propositions.

### 1.2.5 Homotopy Theory

Homotopy theory originates from topology and has connections to higher category theory. Its main object of study is the notion of *homotopy.* A homotopy essentially is a deformation between two topological maps $f, g : X \to Y$. Classically it is a continuous map $H : X \times [0,1] \to Y$ with $H(-,0) = f$ and $H(-,1) = g$. In this situation one writes $f \simeq g$ and says that they are *homotopic.* This is an "equivalence relation", however two homotopies $H_1, H_2 : f \simeq g$ of the same two functions are not necessarily the same. One can in fact study homotopies between homotopies $\alpha : X \times [0,1] \times [0,1] \to Y$ with $\alpha(-,-,0) = H_1$ and $\alpha(-,-,1) = H_2$. Iterating this process reveals the infinite structure of $\infty$-groupoids.

Homotopy type theory as introduced by The Univalent Foundations Program [1] describes this situation very well. It is a synthetic replica of the analytic homotopy theory from topology. The rules for this theory fit on an A5 sheet of paper and the higher structure is beautifully derived from the *elimination rule* of *equality types*. We later observe that a type with only a single proven element can be very interesting due to this higher structure. This example type resembles $S^1$ and contains the higher structure group $\mathbb{Z}$ as expected.

#### 1.2.6 Order Theory, Category Theory and Higher Category Theory

In the lecture by R. Harper *et al.* [19], an interesting connection is made between order theory, ordinary category theory and type theory. Many concepts are found in all of three theories, similar to the Curry-Howard meta-theoretic isomorphism we had before. For example, product types correspond to the category-theoretic product and in order theory this corresponds to the meet construction. The same applies to their duals and the Yoneda lemma. Additionally, Heyting algebras arise naturally from basic type theory and are one of the objects of study in order theory. They also naturally exist as special categories and appear in proof theory.

The natural continuation of this study is to take it from $\infty$-groupoids to $\infty$-categories. We now discuss two possible approaches to add directionality and make morphisms non-invertible.

### 1.3 Two Approaches to add Directionality to HoTT

We already mentioned that HoTT represents homotopy theory in a particularly simple and beautiful way. The same cannot be said about the approaches to adding directionality and exhibiting synthetic $\infty$-categories in a type theory. While the two theories are useful on their own, they are lacking when comparing to the situation that HoTT and $\infty$-groupoids are in. Producing an analogous situation would be the “holy grail” in this area of study and it is very likely to not exist. In this section we present two approaches that both are flawed in their own ways. The first violates certain category theoretic expectations and the second does not build upon HoTT, but requires modifications that make it inherently different from Martin-Löf-Type-Theory.

In this thesis, we study the first approach and explain in detail which category theoretic expectation is violated in Section 3.8. The benefit of this approach is that all statements from HoTT continue to hold. Additionally, the formalization effort using the proof assistant `rzk` [4] covers considerable parts of the theory in [3] including the work from U. Buchholtz and J. Weinberger [5] and C. Bardomiano Martínez [6].

#### 1.3.1 Extending HoTT with Shapes

We give a very quick summary of how E. Riehl and M. Shulman [2] add directionality: we add another layer to our type theory containing *shapes*, essentially simplices $\Delta^n$ and certain subshapes such as their boundaries. We then define $\hom_A(x, y)$ to be the type of maps from $\Delta^1 \to A$ that are $x$ on the initial vertex and $y$ on the terminal vertex. This construction results in types that are not categories, since a priori composition for these morphisms is not possible. However, if one restricts themselves to study *Rezk types*, the theory allows us to prove several statements about $(\infty, 1)$-categories.

In addition to the problem of types not being categories, we cannot construct a sensible category of categories. Even worse, the Grothendieck construction of straightening and unstraightening does not work in the categorical sense. We explore this in detail in Section 3.8.

The advantage of extending HoTT is that all existing statements continue to hold, which is great, since HoTT has been shown to be a useful theory.

### 1.3.2 Restricting Π-Types

The alternative approach is to essentially restrict the formation of Π-types. These are the type theoretic equivalent of the universal quantifier and generalize normal functions. Not all dependent types should be available to be turned into a Π-type. This is explored by D.-C. Cisinski, B. Cnossen, K. Nguyen, and T. Walde [20] in a non type-theoretic way, however it should be possible to extract a type theory from the axioms. The advantage of this approach is that it truly exhibits synthetic category theory. Their work does not even define a notion of "type". There is a category of categories and the Grothendieck construction works as expected.

The disadvantage of this approach is that only cocartesian fibrations support exponentiation (cf. [20, Axiom N]). This corresponds to the formation of Π-types in type theory, so the Π-type only exists for dependent types that are cocartesian families. This means in particular that not all statements from HoTT can be transferred and all need to be carefully examined (for example by proving them using a proof assistant) before relying on them. The work on this theory is not yet finished and there is no proof assistant that formalizes this theory and we will not study it any further here.

### 1.3.3 Holy Grail

Since both approaches give up something from the other world, we cannot reproduce the success of HoTT and $\infty$-groupoids. The category-friendly approach of restricting Π-types will most likely result in the more natural setting for category theory. However, it remains to be seen if formal methods also perform well there.

## 1.4 Motivation

A notion from higher category theory that appears in many places is that of cocartesian fibrations. They play an important role in D.-C. Cisinski, B. Cnossen, K. Nguyen, and T. Walde [20], where they appear in three axioms. In E. Riehl and M. Shulman [2], they are part of the central obstruction to a faithful Grothendieck construction discussed in Section 2.8. Additionally, cocartesian fibrations are in general of great interest in higher category theory, as they allow constructions of functors between infinity categories as discussed by M. Land [21].

For these reasons, we formalize cocartesian fibrations in `rzk` and upstream our work into [3]. In particular we follow work from U. Buchholtz and J. Weinberger [5] and J. Weinberger [7] who study cocartesian fibrations as a special case of *LARI-families*.

## 1.5 Synthetic vs Analytic

We already mentioned that we take a *synthetic* approach to mathematics. It is the dual of the analytic approach – instead of constructing a theory inside of set theory, a synthetic approach aims to axiomatically produce a definition behaving correctly. One famous synthetic theory for the natural numbers is the Peano axioms. These define "0" through just declaring that it exists; additionally there exists a successor for any natural number that is distinct from any number that comes before it. From these few properties (and more axioms such as induction) it is possible to derive many useful statements about natural numbers. In contrast, an analytic approach to defining the natural numbers is for example the von Neumann construction which is embedded in Zermelo Fraenkel set theory. The numbers are defined as sets, starting with $0 :\equiv \emptyset$ and the successor of $n$ as $n \cup \{n\}$. Addition, subtraction, and any other operation is defined using the basic set operations.

Analytic approaches allow uses of existing theories as well as show that the mathematical concept exists in that theory. However, the need to embed them in another theory makes them inherently more verbose to deal with. Often times one has to "look into" the definitions and use facts about the concrete representation to arrive at a proof for the desired statements. This is especially true when formalizing mathematics. When adding the real numbers 2 and 1 in a theory that constructs them as equivalence classes of Cauchy-sequences, then the proof assistant has two options. Either it aids in destructuring the sequences, adding them element-wise and then observing that the resulting sequence also is constant. Alternatively, the mathematician has to perform this extra work and the proof assistant only checks the correctness of that work. Synthetic approaches remedy this issue by building a theory specifically tailored for the concept which results in cleaner proofs and more generality. In any theory where the axioms can be proven analytically, one can rely on any proven statement from the synthetic theory.

## 1.6 Thesis Structure

Section 2 introduces homotopy type theory and covers many important concepts from general type theory. We follow The Univalent Foundations Program [1] and R. Harper *et al.* [19] and explain how `rzk` implements HoTT. Section 3 builds on top of that, adding the extension types from E. Riehl and M. Shulman [2] and introduces many important concepts from category theory. In Section 4, we cover orthogonal fibrations and families in order to define inner families. These are a fundamental notion that is needed for cocartesian fibrations. We also study the behavior of squares in our theory. Next, we cover the main work of our thesis in Section 5, Section 6 and Section 7, where we first introduce the notion of initial section as well as dependent initial sections. We show that they both are equivalent to LARI adjunctions, which are used by U. Buchholtz and J. Weinberger [5] to define LARI families. As our last step before covering cocartesian fibrations, we study LARI-cells, a notion coined by J. Weinberger [7] that generalizes cocartesian fibrations. We show that these are equivalent to our earlier definition of LARI-families. Finally, in Section 7.4, we collect the fruits of our abstract labor and discuss the concrete instantiation of $i_0$-LARI families as cocartesian fibrations. Lastly, in Section 8 we explain how we formalized the majority of our work using `rzk`. During our formalization effort, we discovered a bug in `rzk` that the maintainer promptly fixed. We include a glossary in Appendix A which contains a quick reference to look up symbolic notation.

# 2 Homotopy Type Theory

In this section, we give a summary of homotopy type theory (HoTT), following [1]. We cover everything that we make use of in later sections. We also explain how certain concepts are implemented in `rzk`. For concrete examples, more context and longer explanations, we direct the reader to [1] and [19].

This section has a strong focus on formal correctness and we often explicitly mention details that we gloss over later. We justify being more lax in the later sections by having proven our statements in the proof assistant `rzk` [4].

## 2.1 The Basics of Type Theory

HoTT is heavily based on the type theory developed by Martin-Löf [10]. While several decades have passed since his original paper and lots of other mathematicians have contributed, the original ideas are still visible in our theory. For this subsection, we mainly follow the lecture by R. Harper *et al.* [19].

There are several concepts that we need to cover before we can start writing down our theory. Type theory originates from logic and has connections to constructive mathematics and proof-relevant mathematics. For example, we will *not* assume that the law of excluded middle holds, as it stands in direct conflict with the *univalence axiom*, which we cover later. There is a fundamental connection between logic and type theory described by the Curry-Howard isomorphism. It puts notions from logic and type theory into a 1-1 correspondence:

| **Concept in Logic** | **Concept in Type Theory** |
|---|---|
| Propositions | Types |
| Proof of a proposition $P$ | Term of type $P$ |
| True $\top$ | Unit type $\mathbb{1}$ |
| False $\bot$ | Empty type $\mathbb{0}$ |
| Conjunction $A \wedge B$ | Product type $A \times B$ |
| Disjunction $A \vee B$ | Disjoint union type $A + B$ |
| Implication $A \to B$ | Function type $A \to B$ |
| Existential quantifier $\exists x.P$ | $\Sigma$-type $\Sigma(x : X, P)$ |
| Universal quantifier $\forall x.P$ | $\Pi$-type $(x : X) \to P$ |

Table 2: Curry-Howard's meta-theoretic isomorphism between logic and type theory

This table is an essential part of the motivation for the various types that we introduce in this section. They represent certain logical concepts that a useful mathematical theory ought to have.

Before we dive into the specifics of defining these types, we explain how one sets up the formal language to be able to write them down in the first place.

### 2.1.1 Inference Rules and Judgments

To write down a mathematical theory, several different kinds of notational systems exist. Set theory for example uses first order logic. Since type theory does not sit on top of any logic, we introduce its primitive notions using the same notational system that logicians usually use to describe first order logic (and many other kinds of logics as well).

This notational system is a formal language where we define *judgments* and *inference rules* (often just called *rules*). The judgments are the "*statements* that we can prove within our theory". The rules are the way in which we can derive judgments from other judgments. A rule has the following syntax:

$$\frac{J_1 \quad J_2 \quad J_3 \quad \cdots \quad J_n}{J} \; R_1$$

This rule is named $R_1$. The $J_i$ and $J$ are judgments. The $J_i$ that sit above the horizontal line are the *premises* of the rule and the $J$ on the bottom is the *conclusion*. We allow rules with zero premises, these rules essentially correspond to axioms of our theory, as they have no preconditions. A rule can be used to derive its conclusion given that all the premises hold. We have *variables* in our judgments, which are essentially "all-quantified" in our rules. This is not an internal quantifier of the theory, but the result of allowing syntactic substitution into these rules. Essentially they are *macros* where we can plug in any value for the variables and they hold for any value we can substitute.

Rules are used to construct proofs of judgments. To prove a judgment, we need to write down a series of steps where we are only allowed to apply inference rules. At the end we should arrive at the conclusion that the desired judgment holds. Logicians have invented special notation for this purpose: proof trees. A simple example can be found in Example 2.1.1. The inference rules can be combined like LEGO bricks, plugging a conclusion of one rule into the premise of another. If at the end all leafs of the tree are horizontal lines (signifying that they prove their conclusion without any premises), we can conclude that the judgment at the bottom can be derived from the axioms.

**Example 2.1.1:** Basic proof-tree

$$\frac{\dfrac{\dfrac{}{J_1} \quad \dfrac{}{J_2}}{J_3} \quad \dfrac{}{J_4} \quad \dfrac{\dfrac{}{J_1} \quad \dfrac{\dfrac{}{J_2} \quad \dfrac{}{J_4}}{J_6}}{J_5}}{J_7}$$

A playful interpretation of this is that of a game like chess. Our axioms are the starting positions of the pieces and the rules are the game rules that dictate how the pieces are allowed to move. A judgment is then just a certain position of the pieces. A proof is a sequence of rules that when applied in order to the starting position results in the position that one wanted to prove.

All rules of a theory need to be justified, as we state them to be true without proof. To be strict, we would also need to justify the system itself. The arguments that one can give about this are beyond the realm of mathematics. We leave them to philosophers and are satisfied with the simple fact that this system permits us to describe very interesting and useful mathematical concepts.

**Toy Logic**

Before we dive into the concrete rules and judgments of our type theory, we give a quick example of (incomplete) propositional logic. It is a simple system of judgments and inference rules and will give the reader the chance to obtain a more intuitive understanding of the logical language before our much more complicated rules of type theory later in this section. We also introduce the very important concept of entailment, contexts, and much more.

We first need to describe what the judgments of our toy logic are. The first kind of judgment that we have are just plain variables $A, B, C, \ldots$; we also allow indexed variables.

The next kind of judgment also appears later in our type theory (though in a more complicated form). It is called *entailment* and it requires us to define a new syntactic concept: *contexts*. A context is just a list of variables. We use capital Greek letters ($\Gamma, \Xi, \Theta$) to denote a context. Entailment then allows us to denote the fact that a judgment is derivable from a list of variables. We write $A_1, \ldots, A_n \vdash B$ where the $A_i$ are the premises of the entailment and $B$ is the conclusion. When we write $\Gamma \vdash B$, we speak it as "in context $\Gamma$, $B$ holds". Entailment is mainly governed by these two rules:

$$\frac{}{A \vdash A}\text{R} \qquad\qquad \frac{\Gamma \vdash A \quad \Gamma, A \vdash B}{\Gamma \vdash B}\text{T}$$

The first rule is that of *reflexivity:* all statements in our logic must follow from themselves. The second rule is *transitivity*, it corresponds to plugging in the proof of a lemma into the proof that assumes the lemma to be true. In addition to these two rules, there are three more that allow manipulating the context of an entailment:

$$\frac{\Xi \vdash B}{\Gamma, \Xi, \Theta \vdash B}\text{W} \qquad \frac{\Gamma, A, \Xi, A, \Theta \vdash B}{\Gamma, A, \Xi, \Theta \vdash B}\text{C} \qquad \frac{\Gamma \vdash A}{\pi(\Gamma) \vdash A}\text{P}$$

The first rule is called *weakening* and allows to add unnecessary variables to the context. The second is *contraction* which allows to discard duplicate variables. And the last rule is *permutation*, which allows reordering the context; in it $\pi(\Gamma)$ means "apply the permutation $\pi$ to the elements of $\Gamma$".

We phrased these rules in a way where it's nice to apply them. It would suffice to have the following different rules, but then we would have many more intermediate steps:

$$\frac{\Gamma \vdash B}{\Gamma, A \vdash B}\text{W}' \qquad \frac{\Gamma, A, A \vdash B}{\Gamma, A \vdash B}\text{C}' \qquad \frac{\Gamma, A, B, \Xi \vdash C}{\Gamma, B, A, \Xi \vdash C}\text{P}'$$

This concludes all the judgments in our toy logic and we can move onto the first primitive notion of our propositional logic: conjunction. We introduce the new syntax of expressions. Variables are expressions and given two expressions $E_1, E_2$ we declare that $E_1 \wedge E_2$ also is an expression. To give this syntax its semantic meaning (that matches the intuition of conjunction) we need three rules:

$$\frac{\Gamma \vdash A \quad \Gamma \vdash B}{\Gamma \vdash A \wedge B}\wedge\text{-I} \qquad \frac{\Gamma \vdash A \wedge B}{\Gamma \vdash A}\wedge\text{-E}_1 \qquad \frac{\Gamma \vdash A \wedge B}{\Gamma \vdash B}\wedge\text{-E}_2$$

The first rule is called the *introduction* rule, it allows us to form a conjunction given that both parts of the conjunction are true. The second and third rules are the *elimination* rules and they allow us to prove the two constituents from their conjunction.

Since in all of these rules the context $\Gamma$ appears as the sole context on the top and the bottom, it does not convey any new information. By this argument, we can omit it. We thus also permit writing rules in their *local rule form* if possible:

$$\frac{A \quad B}{A \wedge B}\wedge\text{-I} \qquad \frac{A \wedge B}{A}\wedge\text{-E}_1 \qquad \frac{A \wedge B}{B}\wedge\text{-E}_2$$

In addition to conjunction we also have implication in our toy logic. Given two expressions $E_1, E_2$, we declare $E_1 \to E_2$ to also be an expression. The introduction and elimination rules are:

$$\frac{\Gamma, A \vdash B}{\Gamma \vdash A \to B}\to\text{-I} \qquad \frac{A \to B \quad A}{B}\to\text{-E}$$

We note that the introduction rule has to be written with an explicit context, since we use two different ones. The elimination rule is commonly called "modus ponens".

Finally, we allow parenthesized expressions, so if $E$ is an expression, then so is $(E)$. We forbid any expressions that cause ambiguity in the order of operations (one can also define a precedence between the operations to avoid unnecessary parenthesis).

Using these rules we can now prove the statement $A \to (B \to (A \land B))$ from the empty context:

$$\dfrac{\dfrac{\dfrac{\dfrac{\dfrac{}{A \vdash A}\,\mathrm{R}}{A, B \vdash A}\,\mathrm{W} \qquad \dfrac{\dfrac{}{B \vdash B}\,\mathrm{R}}{A, B \vdash B}\,\mathrm{W}}{A, B \vdash A \land B}\,\land\text{-I}}{A \vdash B \to (A \land B)}\,\to\text{-I}}{\vdash A \to (B \to (A \land B))}\,\to\text{-I}$$

In addition to these concepts, we need to allow writing definitions. We use a colon and an equals ":=" for this operation. It introduces an abbreviation on the left side for the right. For example:

$$X :\equiv A \to B$$

This means that whenever one sees $X$ in the future, one should substitute it with $A \to B$. When substituting one needs to pay attention to operator precedence, which we do not define in our toy logic.

Two concepts from logic that come up later again are disjunction, false and true. The two constants have particularly simple rules:

$$\dfrac{}{\top}\,\top\text{-I} \qquad\qquad \dfrac{\bot}{A}\,\bot\text{-E}$$

Truth introduction says that in any context we can derive true. It does not have an elimination rule, we cannot prove anything from "true". Dually, false obviously doesn't have an introduction rule. Its elimination rule assumes false and allows one to show any primitive judgment $A$.

Disjunction should be familiar to the reader. The introduction rules are also very simple and reflect the duality to conjunction. One potential novel aspect is the elimination rule.

$$\dfrac{A}{A \lor B}\,\lor\text{-I}_1 \qquad \dfrac{B}{A \lor B}\,\lor\text{-I}_2 \qquad \dfrac{\Gamma, A \vdash C \quad \Gamma, B \vdash C}{\Gamma, A \lor B \vdash C}\,\lor\text{-E}$$

The elimination rule allows us to prove any statement $C$ given $A \lor B$, when we can prove $C$ from both $A$ and $B$. An interesting note here is that this adds some ambiguity to our theory, as we have to choose the $C$ when invoking this rule. When we later discuss product types, we lose this ambiguity *because* of the proof terms in our type theory. Before we move onto type theoretic entailment, we give one very important comment about the meaning of inference rules in type theory.

### 2.1.2 Rules and Type-Checking

Inference rules in type theory play a much less important role than in logic, because they are not used directly to construct proofs. The rules in type theory dictate which *terms* are valid; this is called *type checking*. However, the rules together with the syntax of terms are set up in a way where they can be applied automatically.

Because of this property, type theory is naturally suited to be formalized with the help of the computer. Type checking can be implemented via a program which leaves the user to provide valid terms.

In the rest of this thesis, we are not using these rules explicitly again. Instead, we claim that a certain term is valid and since it is just an algorithmic exercise to check that fact, we do not give proof trees beyond this section.

This is the neat part about type theory and a big step up from logic. A proof tree is a cumbersome object – a term is much lighter and easier to write down (it is still difficult to *come up* with the right term, just as it is difficult to come up with a *proof*).

This is further corroborated by the fact that in the proof assistant `rzk` [4] – which we use to formalize our results – there is no way to explicitly use an inference rule. They are all fully implemented by the type checker and not visible to the user. Type checking is only shown to the user when an error needs to be reported, because they formed a term that does not type-check. `rzk` prints a stacktrace of steps explaining the expected state of things. In Section 8 we explain this in more detail.

### 2.1.3 Type Theoretic Entailment

We begin our definition of type theory with a bit of syntax and then introduce contexts. They play a much more important role than in propositional logic, as they contain *more* information: they allow us to define *dependent types and terms* which play a crucial role in our theory.

We write $\Gamma$ **ctx** to *declare* a context. We denote the empty context by $\cdot$ when its absence would be confusing. In order to describe what a context *is*, we also need the syntactic construct of *type declaration*. We write $\Gamma \vdash A$ **type** to say that $A$ is a type in the context $\Gamma$. Additionally, we denote *term declaration* by $a : A$, which states that the term $a$ has type $A$. Creation of contexts then is governed by the following rules:

$$\frac{}{\cdot\ \textbf{ctx}}\ \textbf{ctx}\text{-I}_1 \qquad\qquad \frac{\Gamma\ \textbf{ctx} \qquad \Gamma \vdash A\ \textbf{type}}{\Gamma, a : A\ \textbf{ctx}}\ \textbf{ctx}\text{-I}_2$$

**Remark 2.1.2**

We note that specifying the rules for forming contexts as inference rules technically makes our type theory a multi-sorted theory. Since when we write down an inference rule later, the formation of the contexts there need to follow the rules we wrote above.

This is not uncommon, in fact [19] does it without even mentioning this.

When we write $a : A \vdash \boldsymbol{J}$ for some judgment $\boldsymbol{J}$, then this means that the variable $a$ is allowed to occur in that judgment. When we later add a notion of *universe*, we can get rid of this implicit dependence when using $\boldsymbol{J}$.

Next we give the rules of reflexivity and transitivity for entailment. The transitivity rule of entailment takes a slightly different form in type theory requiring a new syntactic concept called *substitution*. It is also more commonly called the substitution rule:

$$\frac{\Gamma \vdash A\ \textbf{type}}{\Gamma, x : A, \Xi \vdash x : A}\ \text{R} \qquad\qquad \frac{\Gamma \vdash A\ \textbf{type} \qquad \Gamma, x : A, \Xi \vdash \boldsymbol{J} \qquad \Gamma \vdash a : A}{\Gamma, \Xi[x \mapsto a] \vdash \boldsymbol{J}[x \mapsto a]}\ \text{S}$$

The syntax $[x \mapsto a]$ applies substitution to the variable directly written in front of it. Substituting changes all occurrences of $x$ to $a$. We allow to write $M[x \mapsto a, y \mapsto b]$ instead of $(M[x \mapsto a])[y \mapsto b]$.

**Remark 2.1.3**

Substitution complicates the theory due to having to distinguish free and bound variables. We fully glance over this technicality, as there is a canonical solution: do not substitute bound variables. This can be implemented by a computer. Alternatively, we could only ever use each variable once. We allow repeats of variables for convenience. To find the type of a variable one only considers the most recent declaration when walking the tree of subexpressions upwards.

A variable is *bound* if one does not need to add it to the context for the term to type-check. A *free* variable is one that needs to be present in the context for the term to type-check.

We do not have a single kind of weakening rule, instead we give a *rule schema*. We have one inference rule for each judgment $\boldsymbol{J}$ (that fits into the syntactic position) that takes the form given below:

$$\frac{\Gamma, \Xi \vdash \boldsymbol{J} \qquad \Gamma \vdash A \textbf{ type}}{\Gamma, x : A, \Xi \vdash \boldsymbol{J}}\,\mathrm{W}$$

Contraction and permutation also have to be rule schemas:

$$\frac{\Gamma \vdash A \textbf{ type} \qquad \Gamma, x : A, y : A, \Xi \vdash \boldsymbol{J}}{\Gamma, z : A, \Xi[x \mapsto z, y \mapsto z] \vdash \boldsymbol{J}[x \mapsto z, y \mapsto z]}\,\mathrm{C}$$

$$\frac{\Gamma \vdash A \textbf{ type} \qquad \Gamma \vdash B \textbf{ type} \qquad \Gamma, a : A, b : B, \Xi \vdash \boldsymbol{J}}{\Gamma, b : B, a : A, \Xi \vdash \boldsymbol{J}}\,\mathrm{P}$$

The permutation rule can only be applied when $a$ does not appear in $B$ and $b$ does not appear in $A$.

**Remark 2.1.4:** Dependent Types & Terms

The novelty compared to our toy logic is that we can have *dependent terms/types*. For example, if we have a type $\mathbb{N}$ of natural numbers, we can declare the dependent type of sequences of length $n$ where $n$ is a natural number. Written in our formal language we could express this via the judgment:

$$A \textbf{ type}, n : \mathbb{N} \vdash \mathrm{Seq}_n(A) \textbf{ type}$$

In `rzk`, there is no entailment. Instead of dependent types, it uses Π-types and universes (both of which we introduce later). Definitions are also allowed to specify a context which is translated into a Π-type.

**Remark 2.1.5**

The proof-theoretic perspective on type theory is that elements of types are proofs of the statements that the types represent. We can have multiple unequal proofs of the same statement, which makes our theory proof-relevant. Note that classical mathematics is proof-irrelevant, it is not possible from within the theory to distinguish two different proofs of the same statement. It is not even possible to compare them in the first place, as a judgment either holds or it does not.

A type that has an element is called *inhabited* and it corresponds with a true statement from classical mathematics. A proof without any elements is called *uninhabited* and it corresponds with false.

Another concept that our toy logic lacked was that of equality, so we introduce that next.

### 2.1.4 Judgmental Equality

The next concept that we go over is *judgmental equality*. It is also called *definitional equality*. It describes when two syntactic elements are "the same". This means that throughout our theory, we should be able to substitute one for the other anywhere without incurring a semantic change. We need to distinguish this kind of equality from the *propositional equality* which we introduce later as a *type* using the normal $=$ sign, whereas judgmental equality uses $\equiv$. The $=$-type is the most important type from homotopy type theory, giving it the interesting infinite structure.

We write $\Gamma \equiv \Gamma'$ to denote that the two contexts $\Gamma$ and $\Gamma'$ are judgmentally equal. For types we write $\Gamma \vdash A \equiv A'$ (if $\Gamma \vdash A$ **type** and $\Gamma \vdash A'$ **type**) to denote that $A$ and $A'$ are judgmentally equal in the context $\Gamma$. Finally, we also allow to equate elements of types, so given that $\Gamma \vdash a : A$ and $\Gamma \vdash a' : A$ we write $\Gamma \vdash a \equiv a'$ to denote that $a$ is judgmentally equal to $a'$. We thus also need three rules of reflexivity for $\equiv$.

$$\frac{\Gamma\ \mathbf{ctx}}{\Gamma \equiv \Gamma}\equiv\text{-R}_1 \qquad \frac{\Gamma \vdash A\ \mathbf{type}}{\Gamma \vdash A \equiv A}\equiv\text{-R}_2 \qquad \frac{\Gamma \vdash A\ \mathbf{type} \quad \Gamma \vdash a : A}{\Gamma \vdash a \equiv a}\equiv\text{-R}_3$$

Judgmental equality is a very fundamental notion and as such it requires compatibility rules with every other concept that we add to our formal language. For types (which we have not yet introduced) these rules are called $\eta$ and $\beta$ rules and we give them when defining each type. To make judgmental equality compatible with entailment, we need the following rules (we do not give names, as we do not use them explicitly):

$$\frac{\Gamma \vdash \boldsymbol{J} \quad \Gamma \equiv \Gamma'}{\Gamma' \vdash \boldsymbol{J}} \qquad \frac{\Gamma \equiv \Gamma' \quad \Gamma \vdash A \equiv A'}{\Gamma, x : A \equiv \Gamma', x : A'}$$

$$\frac{\Gamma \vdash a : A \quad \Gamma \vdash a' : A \quad \Gamma \vdash a \equiv a' \qquad \Gamma, x : A \vdash B\ \mathbf{type}}{\Gamma \vdash B[x \mapsto a] \equiv B[x \mapsto a']} \qquad \frac{\Gamma \vdash a : A \quad \Gamma \vdash a' : A \quad \Gamma \vdash a \equiv a' \qquad \Gamma, x : A \vdash b : B}{\Gamma \vdash b[x \mapsto a] \equiv b[x \mapsto a']}$$

The first rule is a schema that allows us to prove the same things from equal contexts. The second rule describes how we can prove equality of contexts. The third and fourth rules show equality of dependent types and terms given that the values that we substitute are equal.

In `rzk`, judgmental equality also uses the `≡` symbol. It has the same substitution rules in its type checker, thus permitting to exchange judgmentally equal values in all situations.

### 2.1.5 General Syntax

We use $:=$ to make definitions, writing $X := a$ to introduce $X$ as a shorthand notation for $a$ (here it makes little sense, but $a$ is allowed to be any expression, while $X$ needs to be an identifier).

In `rzk` a very similar notation is used:

```
#def X : A := a
```

Here we need to specify the type of $X$ in addition to its representation. In front of the `: A` that specifies the type, we are allowed to put context variables. The syntax requires us to put them in parenthesis:

```
#def X
  ( a : A)
  : A
  := a
```

In our theory, we also have parenthesized expressions, so if $E$ is an expression, $(E)$ also is a valid expression. In certain situations, we allow omitting parenthesis that would strictly speaking be required by the syntax. For example when applying a function to a tuple we allow writing $f(a, b)$ instead of $f((a, b))$. This is ambiguous with applying a function with two parameters to two values, but those function types are equivalent by currying anyways, so we are not running into problems.

### 2.1.6 Unit Type

The first type that we introduce is going to be a very simple one. It is called the *unit* type, written as $\mathbb{1}$ and it consists of exactly one element for which we use the syntax $*_{\mathbb{1}}$:

**Definition 2.1.6:** Unit type

$$\frac{}{\mathbb{1}\ \textbf{type}}\ \mathbb{1}\text{-F} \qquad\qquad \frac{}{* : \mathbb{1}}\ \mathbb{1}\text{-I}$$

Since we reuse the notation of $*$ for other types that only have a single member, we write $*_{\mathbb{1}}$ if there is ambiguity or we want to make it more explicit which $*$ we talk about.

This type only has a single $\eta$ rule and no $\beta$ rules:

$$\frac{x : \mathbb{1}}{x \equiv *}\ \mathbb{1}\text{-}\eta$$

This expresses that any value of the unit type must be judgmentally equal to $*$.

This type is not particularly interesting, as it is the most trivial type that we can construct. It represents the value of $\top$/true from logic. From category theory, this is the terminal object in $\mathrm{Cat}$. There is only one functor that maps into $\mathbb{1}$ (note that this uniqueness is up to judgmental equality!).

In `rzk` it takes the name `Unit` and the only value of the type is called `unit`. It also has the $\eta$ rule, so whenever we have a definition depending on `x : Unit`, we know that `x ≡ unit`.

### 2.1.7 The Empty Type

The empty type represents false, is written $\mathbb{0}$ and has no elements. In its presence we can prove the inhabitance of any type.

**Definition 2.1.7:** Empty type

$$\frac{}{\mathbb{0}\ \textbf{type}}\ \mathbb{0}\text{-F} \qquad\qquad \frac{A\ \textbf{type} \qquad x : \mathbb{0}}{x : A}\ \mathbb{0}\text{-E}$$

There are no $\beta$ or $\eta$ rules.

We note that this type is not available in `rzk`. Categorically, it is the initial object in $\mathrm{Cat}$ and due to the elimination rule we can always define a function into any arbitrary type as expected.

### 2.1.8 $\Sigma$ Types

The next type from Martin-Löf's type theory that we introduce is the $\Sigma$-type. It is a generalization of product types, which we cover first. The syntax to form product types is $A \times B$. We also have two projection functions which we call $\pi_1$ and $\pi_2$. The rules of the product type are (given in local rule form):

$$\frac{A\ \textbf{type}\qquad B\ \textbf{type}}{A\times B\ \textbf{type}}\times\text{-F}\qquad \frac{a:A\qquad b:B}{(a,b):A\times B}\times\text{-I}\qquad \frac{p:A\times B}{\pi_1(p):A}\times\text{-E}_1\qquad \frac{p:A\times B}{\pi_2(p):B}\times\text{-E}_2$$

The rules in order from left to right are called *formation*, *introduction*, and *elimination* rules for the product type. Most types have all three kinds of rules. The *formation* rule describes when we're allowed to even talk about the product type. We omit all the non-dependent "**type**" judgments from our rules, as they are very distracting and contain no new information (it can always be inferred from the context). The second rule is the *introduction* rule. It states how we can *introduce* a value of the product. The next two rules are its counterpart, the *elimination* rules. They dictate what one can do with a value of the product type. Note how we can specify all of these rules using their local rule form, since no contexts are used anywhere.

Note that while $\pi_1$ and $\pi_2$ look like functions, they strictly speaking are syntactic concepts that just *behave* like functions (which we have yet to define). This distinction does not matter in practice and we just treat them as functions in this thesis.

The $\beta$ and $\eta$ rules for compatibility with $\equiv$ are:

$$\frac{a:A\qquad b:B}{\pi_1((a,b))\equiv a}\times\text{-}\beta_1\qquad \frac{a:A\qquad b:B}{\pi_2((a,b))\equiv b}\times\text{-}\beta_2\qquad \frac{x:A\times B}{(\pi_1(x),\pi_2(x))\equiv x}\times\text{-}\eta$$

An alternative way to write down these rules is via equations, since the premises of the rules are obvious from the conclusion (essentially we just need to all-quantify all variables and then run a type inference algorithm):

$$\pi_1((a,b))\equiv a\qquad \pi_2((a,b))\equiv b\qquad (\pi_1(x),\pi_2(x))\equiv x$$

The product type is the type-theoretic analogue of conjunction. This is because when one proves the inhabitance of $A\times B$, one gives a proof for $A$ *and* a proof for $B$.

The generalization that $\Sigma$-types offer is that of *dependence*. The type of the second component is allowed to depend on the *value* of the first. For example, with $\Sigma$-types, we can form a single type of sequences of arbitrary length.

**Definition 2.1.8:** $\Sigma$-types

The syntax for the sigma type is $\Sigma(x:A,B)$ where $A$ is a type and $B$ is a dependent type over a variable $x:A$. This construction binds the variable $x$, which is allowed to occur in $B$. Later when we have universes, this dependence is more clearly labeled by writing $\Sigma(x:A,B(x))$ when we have $B$ being a function from $A$ into the universe, written $B:A\to\mathcal{U}$.

We have the following formation, introduction, and elimination rules:

$$\frac{\Gamma\vdash A\ \textbf{type}\qquad \Gamma,x:A\vdash B\ \textbf{type}}{\Gamma\vdash\Sigma(x:A,B)\ \textbf{type}}\Sigma\text{-F}\qquad \frac{\Gamma,x:A\vdash B\ \textbf{type}\qquad \Gamma\vdash a:A\qquad \Gamma\vdash b:B[x\mapsto a]}{\Gamma\vdash(a,b):\Sigma(x:A,B)}\Sigma\text{-I}$$

$$\frac{\Gamma,x:A\vdash B\ \textbf{type}\qquad \Gamma\vdash p:\Sigma(x:A,B)}{\Gamma\vdash\pi_1(p):A}\Sigma\text{-E}_1\qquad \frac{\Gamma,x:A\vdash B\ \textbf{type}\qquad \Gamma\vdash p:\Sigma(x:A,B)}{\Gamma\vdash\pi_2(p):B[x\mapsto a]}\Sigma\text{-E}_2$$

Additionally, the $\beta$ and $\eta$ rules are the same as for products (but if we now type-check this, we have $b:B[x\mapsto a]$ instead):

$$\pi_1((a,b)) \equiv a \qquad \pi_2((a,b)) \equiv b \qquad (\pi_1(x), \pi_2(x)) \equiv x$$

We allow writing $\Sigma(a : A, b : B, C)$ which should be interpreted as $\Sigma(a : A, \Sigma(b : B, C))$. We also do the same for terms, so $(a, b, c) \equiv (a, (b, c))$. In cases where the types inside take up lots of vertical space, we also use the bigger "$\sum$" symbol.

The literature on type theory often uses a slightly different notation for Σ-types:

$$\sum_{x:A} B \equiv \Sigma(x : A, B)$$

We opted for our notation, as it avoids *index-drift*, consider for example the following type in both notations:

$$\sum_{x : \sum_{y : \sum_{a : A} B} C} D \qquad \text{vs} \qquad \Sigma(x : \Sigma(y : \Sigma(a : A, B), C), D)$$

**Example 2.1.9**

We give a basic example of a Σ-type. Assume that we have the type of natural numbers $\mathbb{N}$ and a type of sequences of length $n : \mathbb{N}$ written $\mathrm{Seq}_n(A)$ where $A$ is a type.

We can now form the Σ-type:

$$\Sigma(n : \mathbb{N}, \mathrm{Seq}_n(A))$$

This type is the type of sequences of arbitrary length. An element of it consists of a length $n : \mathbb{N}$ and a sequence of said length.

The logical equivalent of this type is the statement "there exists a natural number such that there is a sequence of that length in $A$". This statement is very trivial and shows that all the interesting properties of this type must come from the *difference* between the elements.

Our notation also places the same priority on both types of the sum type whereas the summation notation makes the second component much more important. In addition, our notation is much closer to what is used in `rzk`:

```
Σ (a : A) , B
```

One disadvantage from using this notation is that at high nesting levels it becomes difficult to discern which parenthesis belong to which Σ. We use bigger parenthesis, vector notation or peculiar newlines/formatting to prevent this from impacting readability. When writing `rzk` code, the special formatting used by [3] ensures that any nesting level is easily recognized.

As we mentioned before, the Σ-type corresponds with the existential quantifier from logic: in order to give an element of $\Sigma(x : A, B)$ (and thus a proof that $\exists x : A.\ B$ holds), one first has to provide a witness $a : A$ and then a proof $b : B[x \mapsto a]$. This proof $b$ shows that for our witness $a$, the type $B$ (with $x$ replaced by $a$) is inhabited (logically speaking: the witness $a$ satisfies the condition $B$).

To recover the normal, independent product type, we just need to choose a non-dependent type (on $x$) for $B$. In this case, we allow writing $A \times B \equiv \Sigma(x : A, B)$ (with $x$ not occurring free in $B$).

The category theoretic version of Σ-types is called the "dependent sum" [22]. It is rather uncommon to learn, since it works directly with the fibration version of what we later call *type families*. This makes it significantly more verbose to work with. It generalizes the notion of coproducts. In our type theory, we can form any categorical limit using the $\Sigma$ type.

**Remark 2.1.10**

Sometimes one encounters exceedingly big $\Sigma$-types. In this case, writing elements as well as the type itself as vectors increases readability considerably, especially when the terms of each component contain parentheses and many other subexpressions. The vector always lists the items from top to bottom. For example:

$$\begin{pmatrix} a \\ b \\ c \\ d \\ e \end{pmatrix} : \Sigma \begin{pmatrix} a : A \\ b : B \\ c : C \\ d : D \\ E \end{pmatrix}$$

We in addition also allow omitting the name of elements in $\Sigma$-types when they don't occur in the later types. This could also be written using $\times$, but is in some circumstances more space efficient.

### 2.1.9 $\Pi$ Types

Next up we define $\Pi$-types. These are a generalization of function types and can be also interpreted as a product type from which they inherit their name. Their logical counterpart is the universal quantifier $\forall x.P$; if one has a value $x$ of type $A$, then the proof of "$\forall x.P$" allows one to obtain a proof of $P[x \mapsto a]$. The notation is very similar to $\Sigma$-types, but uses a capital pi instead $\Pi(a : A, B)$. Values of this type use the classical lambda notation $\lambda a\colon\ A.\ b$ where $x$ is a variable of type $A$ bound after the dot and $b$ is a term.

The product interpretation is that a value of $\Pi(i : I, E_i)$ consists of $e_i : E_i$ for each $i : I$. The type-theoretic interpretation of the $\Pi$-type as a dependent function type is the most useful one to us. That is because they exhibit a behavior similar to continuous maps – only allowing their values to vary in a *consistent* manner. We observe this somewhat nebulous behavior later in this section when we define the notion of contractibility and again much later in Section 5.

**Definition 2.1.11:** $\Pi$-types

$$\frac{A\ \textbf{type} \qquad \Gamma, x : A \vdash B\ \textbf{type}}{\Gamma \vdash (x : A) \to B\ \textbf{type}}\ \Pi\text{-F} \qquad \frac{\Gamma, x : A \vdash B\ \textbf{type} \qquad \Gamma, x : A \vdash b : B}{\Gamma \vdash \lambda x.\ b : (x : A) \to B}\ \Pi\text{-I}$$

$$\frac{\Gamma, x : A \vdash B\ \textbf{type} \qquad \Gamma \vdash f : (x : A) \to B \qquad \Gamma \vdash a : A}{\Gamma \vdash f(a) : B[x \mapsto a]}\ \Pi\text{-E}$$

The $\beta$ and $\eta$ rules are:

$$(\lambda x.\ b)(a) \equiv b[x \mapsto a] \qquad\qquad (\lambda x.\ b(x)) \equiv b$$

The $\eta$ rule (the second one) only holds when $x$ is not free in $b$.

We allow writing $(a : A, b : B, c : C) \to D$ instead of $(a : A) \to (b : B) \to (c : C) \to D$. In the independent case, we also write $A \to B \to C \to D$ which should not be confused with the usual mathematical meaning of a composed function from $A$ to $D$. When we want to write composi-tions, we always place names above the arrows. This notation for the Π-type is again similar to the notation in `rzk`. Classically, one would use the product symbol for this type:

$$\prod_{a:A} B$$

However, this also has the problem of index drift, hence we prefer our notation.

For the lambda expression, we also allow to put a type after the variable in order to aid the reader in understanding the type of a standalone lambda expression: "$\lambda x\colon A.\ b$". We also allow pattern matching of Σ-types the argument. Suppose that we want to write down a type with the domain $\Sigma(a : A, B)$, we then do it like this: $\lambda(a, b).\ x$ where $x$ is allowed to be any term.

As with Σ-types, we allow vector-like notation when writing down the type and when defining maps:

$$(a : A, b : B, c : C) \to D \equiv \begin{pmatrix} a : A \\ b : B \\ c : C \end{pmatrix} \to D$$

$$(a, b, c) \mapsto d \equiv \begin{pmatrix} a \\ b \\ c \end{pmatrix} \mapsto d$$

Similarly to Σ-types, we can recover the standard notion of non-dependent function type (which we just write $A \to B$).

Category theoretically, this is the dependent product [23], which similarly to the dependent sum is rarely taught and is formulated using fibrations rather than families.

**Definition 2.1.12:** Identity functions

For any type $A$, we write $\mathrm{id}_A$ for the identity function $A \to A$ defined as $\lambda a.\ a$.

In `rzk` it is defined like this:

```
#def identity
  ( A : U)
  : A → A
  := \ a → a
```

The syntax of Π types in `rzk` is `(a : A) → B` and for a lambda expression it uses a backslash followed by the arguments and then an arrow followed by the body of the lambda `\ x y → t`.

The `A : U` part declares that `A` is of the type `U`, which is the universe of types. We introduce it next as a replacement for the "$A$ **type**" judgment.

After the next section, we are able to define a *single* identity function and not require a schema definition for every type.

**Definition 2.1.13:** Function composition

Given two composable functions $f : A \to B, g : B \to C$, we define:

$$g \circ f := \lambda a.\ g(f(a))$$

### 2.1.10 Universes

A universe is a *type of types*— a type whose elements are types themselves. There are several justifications for why one would want to have universes in their theory. One of them is the Grothendieck construction, which we cover in Section 2.8; one cannot perform that construction without having type universes. We also want to have universes to take full advantage of dependent types: a universe turns types into elements of types, making them non-special values of the theory. Finally, universes allow us to reduce the number of inference rules needed by certain types.

**Example 2.1.14**

Consider the type `bool` with the following rules:

$$\frac{}{\texttt{bool}\ \textbf{type}} \qquad \frac{}{\texttt{t} : \texttt{bool}} \qquad \frac{}{\texttt{f} : \texttt{bool}}$$

$$\frac{\Gamma, b : \texttt{bool} \vdash C\ \textbf{type} \quad \Gamma \vdash x : C[b \mapsto \texttt{t}] \quad \Gamma \vdash y : C[b \mapsto \texttt{f}] \quad \Gamma \vdash v : \texttt{bool}}{\Gamma \vdash \texttt{if}(v, x, y) : C[b \mapsto v]}$$

$$\text{if}(\texttt{t}, x, y) \equiv x \qquad \text{if}(\texttt{f}, x, y) \equiv y$$

With our current rules it is impossible to define a type that is dependent on a `bool`. So the variable $b$ in $C$ is essentially useless and thus $x$ and $y$ have the same type. However, allowing them to be of entirely different types is useful and very much expressible if we were to introduce universes.

An alternative to universes in this case would be to add a type-level `if`-expression:

$$\frac{\Gamma \vdash X\ \textbf{type} \quad \Gamma \vdash Y\ \textbf{type} \quad \Gamma \vdash v : \texttt{bool}}{\Gamma \vdash \texttt{IF}(v, X, Y)\ \textbf{type}}$$

$$\texttt{IF}(\texttt{t}, X, Y) \equiv X \qquad \texttt{IF}(\texttt{f}, X, Y) \equiv Y$$

We can then use this to define a type dependent on a `bool`. Given $x : X$ and $y : Y$ we have

$$\texttt{if}(v, x, y) : \texttt{IF}(v, X, Y)$$

This approach is called “large elimination” and detailed in [19]. The approach is cumbersome, since now we have to add two constructors for all of these kinds of types whereas with universes, we can phrase this succinctly.

In our theory, we denote the universe using $\mathcal{U}$. To formally introduce it to our theory, we would like to add the following rules:

$$\frac{}{\mathcal{U} : \textbf{type}} \qquad \frac{\Gamma \vdash A\ \textbf{type}}{\Gamma \vdash A : \mathcal{U}}$$

However, this makes the judgment $\mathcal{U} : \mathcal{U}$ derivable, which leads us into the same conundrum that exists for set theory and a set of sets: Russell’s paradox, see [1] and [19]. To avoid the paradox, we use a cumulative hierarchy of universes. $\mathcal{U}_i$ is a universe of types of “size $i$”. The next level up, $\mathcal{U}_{i+1}$ contains all types from the previous level and $\mathcal{U}_i$ itself as well. So instead of the two inconsistent rules from above, we have the following rules:

$$\frac{}{\mathcal{U}_i : \mathcal{U}_{i+1}}\ \mathcal{U}\text{-F} \qquad \frac{A : \mathcal{U}_i}{A : \mathcal{U}_{i+1}}\ \mathcal{U}\text{-I}$$

We then need to replace our existing rules that use the $A$ **type** primitive with $A : \mathcal{U}_i$ (we can then also fully remove the **type** primitive from our theory). So for example, we have the following formation rule for $\Sigma$-types now:

$$\frac{\Gamma \vdash A : \mathcal{U}_i \quad \Gamma, x : A \vdash B : \mathcal{U}_i}{\Gamma \vdash \Sigma(x : A, B) : \mathcal{U}_i} \Sigma\text{-F}$$

Similar to how we can omit the context from our rules when it is self-evident, we can do the same for the level of the universe. In fact, one can implement an inference algorithm that automatically chooses the lowest possible numbers for the levels for all universes involved. Hence, we now write $\mathcal{U}$ without the index. If we write $\mathcal{U} : \mathcal{U}$, that means that the index of the left universe must be lower than the one on the right.

These levels are part of the meta-theory and unavailable to the inside: there is no way within the theory to make a level dependent on values from the theory. For example, the type $(n : \mathbb{N}) \to \mathcal{U}_n$ does not exist, it is a syntactic error to put anything except a level into the index of $\mathcal{U}$. This is important, because if it were possible, then we would have to have $\mathcal{U}_{\aleph_0}$ as well as higher ordinals.

Any definition or theorem that uses $\mathcal{U}$ is a *schema definition/theorem* because it is meta-generic over all level indices that occur. This will never come up in practice, since we can just treat $\mathcal{U}$ as a single type.

We can now use universes to exhibit dependent types as a type themselves:

**Definition 2.1.15:** Type Family

A type family over a given type $A : \mathcal{U}$ is just a function $B : A \to \mathcal{U}$.

Type families are more convenient compared to dependent types, because we can delegate all the complications that come from variable substitution to the $\lambda$ terms and just write $B(a)$ for the type corresponding with $a$.

**Remark 2.1.16**

From now on we get even more lazy with our notation, assuming the reader is able to derive certain conditions from the way we write our statements. For example, we might begin a statement with “Let $B : A \to \mathcal{U}$ be a type family.” in this case, the reader should assume that we also require $A : \mathcal{U}$ without any other restrictions.

We only use this when there is only “one sensible choice” for the unspecified variables. It should always be assumed that it has the minimal capabilities required from the context.

Since we use type families extensively, we introduce useful notation and notions for them:

**Definition 2.1.17**

For a type family $B : A \to \mathcal{U}$, we define
- The *total type* of $B$ is $\Sigma(a : A, B(a))$. We also write $\Sigma B$ for this type.
- The *projection* of the total type of $B$ is the function:

$$\pi_B : \Sigma B \to A \quad (a, b) \mapsto a$$

- The *type of dependent functions* of $B$ is $(a : A) \to B(a)$. We also write $\Pi B$ for this type. It is also called the *sections* of $B$, which we introduce and explain later in Section 5.

We can now define a single identity function:

**Definition 2.1.18:** The identity function

The identity function id is defined with the following type and term:

$$\mathrm{id} : (A : \mathcal{U}, a : A) \to A$$
$$\mathrm{id} :\equiv \lambda A \colon \mathcal{U}.\ \lambda a \colon A.\ a$$

We still continue to write and use $\mathrm{id}_A : A \to A$ meaning $\mathrm{id}(A)$ with the above definition.

**Remark 2.1.19**

In `rzk` we have in fact no way to specify a dependent type, since entailment does not exist. Instead, we always give type families. There also is no real distinction between putting a variable in the context of a definition and making it a parameter of a Π-type inside the definition. So our `identity` function that we defined like this before:

```
#def identity
  ( A : U)
  : A → A
  := \ a → a
```

can equivalently be defined like this:

```
#def identity
  : (A : U) → A → A
  := \ _ a → a
```

## 2.2 Homotopies

Having finished the most important concepts from Martin-Löf's type theory, we now introduce the most important notion of HoTT: *propositional equality*, also called *identity types*, *path types* or just $=$ *-types*. We still follow [1] and [19]. As we already discussed when introducing judgmental equality, these two notions are *different*. Judgmental equality is a *judgment*, so it either holds or it doesn't and propositional equality is a type, so we have an element of the type $x \underset{A}{=} y$ which proves the fact that $x$ is equal to $y$ (in the type $A$) in our theory. However, this proof is not something that "either exists or doesn't", but there is possibly some additional data attached to this proof. It stores *in which way* $x$ is equal to $y$ and there might be multiple *different* ways in which two values may be equal.

This brings us to an important observation: we can *nest* the propositional equality types, essentially asking the question if two proofs $p, q : x \underset{A}{=} y$ are equal themselves, so if $p \underset{x=y}{=} q$ is inhabited. This question can now be asked recursively, given two proofs that two proofs are equal. Writing it in equations, this gets difficult to parse relatively quickly:

$$p, q : x \underset{A}{=} y \qquad H_1, H_2 : p \underset{x \underset{A}{=} y}{=} q \qquad H_1 \underset{p \underset{x \underset{A}{=} y}{=} q}{=} H_2$$

This is a glimpse at the infinite structure of identity types. One might be worried about the index drift of this notation. But since we always need to specify the elements, the type within which we consider equality is almost always known from context. Additionally, we rarely consider multiple nestings of equality types and instead write our theorems in a generic form where one is allowed to plug in any type.

The main alternative interpretation of $=$-types is that of paths: $p : x =_A y$ is a path from $x$ to $y$ in $A$. This interpretation of course takes a homotopical viewpoint. All properties that are interesting for homotopy theory arise from this type's elimination rule. It is rather cryptic if one has never seen it before, so it will take some time to get used to.

**Definition 2.2.1:** $=$-types

We already gave the syntax of $=$ types. Given two elements $x, y : A$, we write $x =_A y$ for the type of *paths* or *equalities* between $x$ and $y$. When the type of $x$ and $y$ are clear from the context, as is usually the case, we also just write $x = y$. For the elimination rule, we need a special symbol for which we use $\mathsf{J}$.

The formation and introduction rules are straight-forward:

$$\frac{A : \mathcal{U} \quad x : A \quad y : A}{x = y : \mathcal{U}}\text{=-F} \qquad\qquad \frac{A : \mathcal{U} \quad x : A}{\mathsf{refl} : x = x}\text{=-I}$$

We also allow writing the type underneath the equals $x =_A y$ in case it is not obvious from the context; the same for the value, so $\mathsf{refl}_x : x = x$. The elimination is as already mentioned rather complicated and we give a more extensive explanation below:

$$\frac{A : \mathcal{U} \quad x : A \quad C : (y : A, p : x = y) \to \mathcal{U} \quad c : C(x, \mathsf{refl}_x)}{\mathsf{J}(C) : (y : A, p : x = y) \to C(y, p)}\text{=-E}$$

Our compatibility with $\equiv$ is given via the following rule:

$$\mathsf{J}(C)(x, \mathsf{refl}_x) \equiv c$$

In the elimination rule, we start out with any type $A$ and then fix a starting point $x : A$. Next, we have a type family $C$ that is allowed to depend on a path in $A$ that starts at $x$, but with a free end-point which is called $y$. As our last ingredient to the rule we have a proof $c$ that shows that the type family $C$ is inhabited for the end-point $x$ and path $\mathsf{refl}_x$. The conclusion of the elimination rule is that $C$ is inhabited for any path starting at $x$. This rule is also called (based) path induction, as showing the (dependent) statement $C$ for $\mathsf{refl}_x$ shows it for all other paths $p : x = y$.

Another perspective for the elimination rule requires us to jump a bit ahead and use the notion of *equivalence* that we have not yet introduced. The elimination rule allows us to make our theory invariant under equivalences. Consider the following commutative diagram:

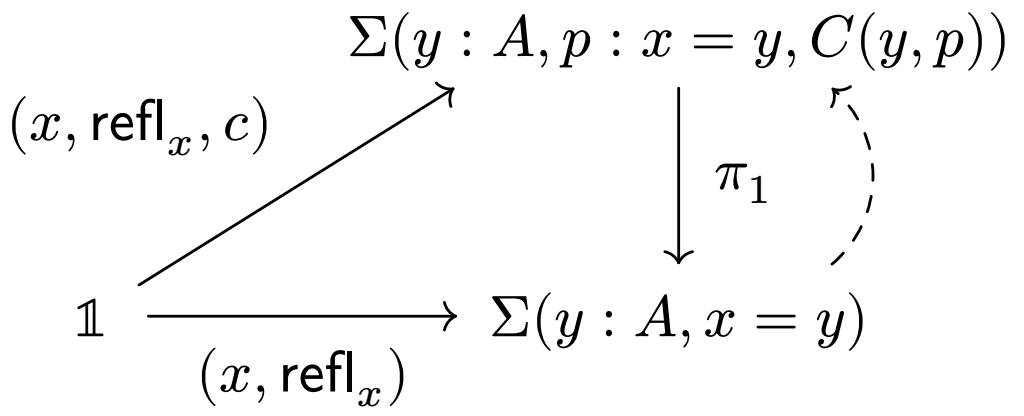

We later show that the bottom map is an equivalence. Through this equivalence, we should be able to obtain a section of the right vertical map that fits into the diagram as the dashed map. From this perspective, the dashed map is the elimination rule for $=$-types, as we have to only show existence of the commutative triangle (which exists as soon as we have a $c$). When we later also discuss the univalence axiom, which equates equivalences and propositional equalities, this interpretation will gain even more credibility.

Note that our formulation of the elimination rule is a bit unusual, the classical variant uncurries the function and allows the family to vary over the starting point $x$:

$$\frac{A : \mathcal{U} \quad C : (x : A, y : A, p : x = y) \to \mathcal{U} \quad c : (a : A) \to C(a, a, \mathsf{refl}_a) \quad x : A \quad y : A \quad p : x = y}{\mathsf{J}(C, p) : C(x, y, p)}$$

While there are subtle differences between our based path induction rule and this non-based path induction rule (depending on the exact way we specify it, as a rule or as a theorem in our theory), they are equivalent. This is shown in [1, Section 1.12.2] and we will not go into the details here. We note that one can only fix *one* of the two endpoints, since if we could fix both endpoints, we would be able to conclude that $p = \mathsf{refl}_x$ for any $p : x = y$, which would collapse all the higher structure.

**Remark 2.2.2**

Judgmental equality of two elements implies propositional equality. That is because judgmental equality is compatible with dependent types and terms. So given $x \equiv y$ we have $(x = y) \equiv (x = x) \equiv (y = y)$ and that $\mathsf{refl}_x : x = y$ as well as $\mathsf{refl}_y : x = y$. We also observe the two refl to be the same: $\mathsf{refl}_x \equiv \mathsf{refl}_y$.

**Remark 2.2.3**

In `rzk`, the syntax is the same without the specified type: `x = y` and with a type: `x =_{A} y`. This is the *only* type that supports type inference. refl is also called `refl` and supports inferring the term. In case one needs to specify it, the syntax is `refl_{x : A}` where `A` is the type of `x`.

### 2.2.1 Groupoidal Structure

The identity types exhibit a groupoidal structure that we can fully derive from the elimination rule. We first have to define operations to invert and concatenate paths. This in turn proves that $x = y$ is an equivalence relation together with the fact $\mathsf{refl} : x = x$.

**Theorem 2.2.4:** $=$ is an equivalence relation

Given a type $A : \mathcal{U}$ with three elements $x, y, z : A$ there exist these two functions:

$$(-)^{-1} : (x = y) \to (y = x)$$
$$- \cdot - : (x = y) \to (y = z) \to (x = z)$$

The first one is called *path inversion* and the second one *path concatenation*.

**Proof.** We construct both functions using the elimination rule for identity types. In the first case we apply it to the input argument and thus we need to show there exists an element of $x = x$ (as $y$ has now been replaced by $x$ through invoking path induction). But that is of course possible, since $\mathsf{refl}_x$ always exists.

Constructing the second function works similarly, we only need to apply path induction twice at which point we have $z \equiv y \equiv x$ and we can again choose $\mathsf{refl}_x$. □

These two functions satisfy the groupoid laws which is proven in [1, Section 2.1]. We cite their result here:

**Theorem 2.2.5:** $x = y$ forms a groupoid

In any type $A : \mathcal{U}$ with $x, y, z, w : A$ and $p : x = y, q : y = z, r : z = w$ we have:

$$p \cdot p^{-1} = \mathsf{refl}_x \qquad p^{-1} \cdot p = \mathsf{refl}_y \qquad \left(p^{-1}\right)^{-1} = p \qquad \mathsf{refl}_x^{-1} \equiv \mathsf{refl}_x$$
$$p \cdot \mathsf{refl}_y = p \qquad \mathsf{refl}_x \cdot p = p \qquad (p \cdot q) \cdot r = p \cdot (q \cdot r) \qquad \mathsf{refl}_x \cdot \mathsf{refl}_x \equiv \mathsf{refl}_x$$

Since all of these only hold propositionally (except the rules about refl), all of this structure only holds "up to higher coherence". This is exactly what is expected of $\infty$-groupoids. Since we can choose $A$ to be an identity type itself, we also get all of this structure for the recursive identity types.

Having defined all of this theoretically, we give the example of $S^1$.

**Example 2.2.6:** The Circle

To define the circle, we use *generators* and *relations* as described in [1, Chapter 6]. The type $S^1$ consists of a single point called $*$ (or $*_{S^1}$ if there is ambiguity) and a path called $\mathsf{loop} : * = *$.

The generators and relations mechanism a priori does not produce any further higher coherences, so we have that $\mathsf{loop} \not\equiv \mathsf{loop} \cdot \mathsf{loop}$ and $\mathsf{loop} \not\equiv \mathsf{loop}^{-1}$. In fact, [1] explicitly proves that $\mathsf{loop} \neq \mathsf{refl}_*$. Additionally, it does not add the elimination rule of $\mathbb{1}$ (the rule proving that $*_{\mathbb{1}} \equiv x$ when $x : \mathbb{1}$), as that would be detrimental to the structure of $* \underset{S^1}{=} *$. This means that while $S^1$ has a canonical element $*_{S^1}$, we cannot prove that it is the *sole* element.

This might seem confusing, but it fits perfectly with the homotopical view of $S^1$. It consists of a single connected component so all points are equivalent, but we cannot assume that there is only one point to begin with.

[1] proves that $\left(* \underset{S^1}{=} *\right) \simeq \mathbb{Z}$ as one would expect.

### 2.2.2 Functoriality

Functions satisfy functoriality with respect to the groupoid structure automatically. This fact is recorded by the following definition:

**Definition 2.2.7:** Function action on paths

Let $f : A \to B$ be a function and $x, y : A$. We define the following function using path induction:

$$\mathsf{ap}_f : (x = y) \to (f(x) = f(y))$$

That sends $\mathsf{refl}_x$ to $\mathsf{refl}_{f(x)}$.

We also need to show that this operation is indeed functorial:

**Theorem 2.2.8:** (cf. [1, Lemma 2.2.2])

For functions $f : A \to B, g : B \to C$ and paths $p : x \underset{A}{=} y, q : y \underset{A}{=} z$, the following statements hold:

$$\mathsf{ap}_f(p \cdot q) = \mathsf{ap}_{f(p)} \cdot \mathsf{ap}_{f(q)} \qquad \mathsf{ap}_f\left(p^{-1}\right) = \mathsf{ap}_{f(p)}^{-1}$$
$$\mathsf{ap}_g\left(\mathsf{ap}_f(p)\right) = \mathsf{ap}_{g \circ f}(p) \qquad \mathsf{ap}_{\mathrm{id}_A}(p) = p$$

## 2.3 Contractibility

In the rest of this section, we introduce important terminology for homotopy type theory that we use later in the thesis. We start with the notion of *contractibility*.

**Definition 2.3.1:** Contractible

We call a type $A : \mathcal{U}$ *contractible* if the type of *contractions* is inhabited:

$$\text{is-contr}(A) :\equiv \Sigma(x : A, (a : A) \to (a = x))$$

The first component $x$ is called the *center* of the contraction and the second one is called the *homotopy* of the contraction.

A contractible type essentially consists of only one element. At first it might seem that this definition only ensures the first level of identity types collapses and higher structure is still permitted. However, that is not the case. Our example of $S^1$ is **not** contractible. This might seem strange at first, but it is due to the *continuity* of Π-types that we described earlier in the section. To prove that $S^1$ cannot be shown to be contractible, we give the following lemma that also justifies the name *contractible* for our definition:

**Lemma 2.3.2**

If $A : \mathcal{U}$ is contractible, then any identity type in $A$ also is contractible.

**Proof.** Let $a$ be the center of the contraction of $A$. We write $p_x : a = x$ for the homotopy of this contraction to any other point $x : A$.

We have to show that $(x : A, y : A) \to \text{is-contr}(x = y)$ holds. We fix $x, y : A$ and define a contraction: as the center, we choose $p_x^{-1} \cdot p_y : x = y$. For the homotopy, we fix a $q : x = y$ and apply path induction to $q$. We observe that the center now is refl: $p_x^{-1} \cdot p_x \equiv \mathsf{refl}_x$ and thus we can use $\mathsf{refl}_{\mathsf{refl}_x}$ as the equality between the center and $q \equiv \mathsf{refl}_x$. □

A very important statement that we use extensively later is that the type of paths with a fixed starting point is contractible. This follows directly from the elimination rule.

**Lemma 2.3.3**

For any type $A$ and any value $a : A$ we have that $\Sigma(x : A, a = x)$ is contractible with center $(a, \mathsf{refl}_a)$.

## 2.4 Propositions

**Definition 2.4.1**

We call a type $A : \mathcal{U}$ a *proposition* if any two inhabitants are canonically equal:

$$\text{is-prop}(A) :\equiv (x : A, y : A) \to \text{is-contr}(x = y)$$

We also define the "universe of propositions" as $\mathbb{P} :\equiv \Sigma(A : \mathcal{U}, \text{is-prop}(A))$ and allow implicit coercions of elements of this type to their first component.

If a type is a proposition, that means that there is only one way that the statement represented by the type can be true. In proof-theoretic terms, any two proofs of a proposition are equivalent. There are multiple alternative ways to define the notion of proposition. As a type where all elements are equal $(x : A, y : A) \to x = y$ or as a type that is contractible if inhabited $A \to \text{is-contr}(A)$.

Since $\to$ corresponds with logical implication, the latter definition in particular is very meaningful in garnering an intuitive understanding of the notion of proposition. A type is a proposition if it is either uninhabited (so there is no proof for it) or it is contractible. The former definition is very useful when one wants to show a certain type to be a proposition.

We note that we can recover classical logic from propositions if we were to introduce the law of excluded middle for them. Later when discussing the univalence axiom, we observe that we cannot do this for arbitrary types as the axioms conflict with each other.

For a given proposition $A$ we say "$A$ holds" to mean that $A$ is inhabited. The convention in `rzk` and type theory in general is to only give a definition the name "is-xyz" if it is a proposition. This is the case for our "is-contr" and "is-prop" definitions. We generally do not prove this fact explicitly. [3] often provides proofs in cases where it is not immediately evident that a type is a proposition.

Notably, being a proposition is closed under several of our type constructors:

**Theorem 2.4.2:** Closure properties of propositions

Given a type family $B : A \to \mathbb{P}$ that is fiberwise a proposition (this is required by it having the codomain $\mathbb{P}$), then

- The associated dependent function type $\Pi B \equiv (a : A) \to B(a)$ is also a proposition.
- If additionally $A$ is a proposition, the total type $\Sigma B \equiv \Sigma(a : A, B(a))$ is also a proposition.

Additionally, any contractible type is a proposition.

## 2.5 Equivalences

The notion of equivalence allows us to equate types in yet another different notion to $=$ and $\equiv$. It is the weakest of the three. Category theorists are very familiar with this notion (where it is often called weak equivalence in the case of higher categories), as it permits equating types that contain vastly different numbers of elements. The deciding factor are the identity types of the compared types. The connected components of the two types are in bijection and the higher structure is itself identified up to equivalence.

**Definition 2.5.1:** Equivalence

Given two types $A, B : \mathcal{U}$ and a map $f : A \to B$, we say that $f$ is an equivalence, if the following proposition holds:

$$\begin{aligned} \text{is-equiv}(f) \coloneqq{} & \Sigma(g_1 : B \to A, (b : B) \to f(g(b)) = b) \\ & \times \Sigma(g_2 : B \to A, (a : A) \to g(f(a)) = a) \end{aligned}$$

$g_1$ is called a *section* of $f$ and $g_2$ a *retraction*.

We also write $f : A \simeq B$ for the type of equivalences from $A$ to $B$.

That this definition is indeed a proposition is not immediately visible from the definition. And the naive definition that just uses a single function that is the inverse is *not* a proposition. However, there is a map from this naive definition to the proper one, that just uses the inverse twice.

The identity function is trivially an equivalence, by using refl as the homotopy value. Additionally, any path between types $A =_{\mathcal{U}} B$ yields via path induction an equivalence $A \simeq B$. We also can compose equivalences:

**Theorem 2.5.2:** (cf. [1, Lemma 2.4.12])

Let $f : A \to B$ and $g : B \to C$ be equivalences. Their composition is also an equivalence (by composing their respective sections/retractions and homotopies).

Another alternative way of obtaining an equivalence is showing that a map has contractible fibers:

**Theorem 2.5.3:** (cf. [1, Theorem 4.4.5])

A function $f : A \to B$ is an equivalence iff $\mathrm{fib}_f(b)$ is contractible for every $b : B$.

Equivalences satisfy the usual 2-out-of-3 property:

**Theorem 2.5.4:** (cf. [1, Theorem 4.7.1])

Given $f : A \to B, g : B \to C$. If any two of $f$, $g$ and $g \circ f$ are equivalences, so is the third.

Equivalences between propositions simply require giving maps in both directions:

**Theorem 2.5.5:** (cf. [1, Lemma 3.3.3])

Given propositions $A, B : \mathbb{P}$ and $f : A \to B, g : B \to A$, then $A \simeq B$.

### 2.5.1 Important Equivalences

There are several important equivalences that we make use of later. Additionally many axioms in type theory are formulated using equivalences. The axiom of function extensionality states that pointwise equal functions are propositionally equal. A similar axiom is added for "extension types" in the next section. Later on in this section we discuss univalence, which also is an axiom declaring an equivalence.

**Axiom 2.5.6:** Function Extensionality (cf. [1, Axiom 2.9.3])

For any type $A : \mathcal{U}$ and any two dependent maps $f, g : \Pi B$ with $B : A \to \mathcal{U}$, the natural map defined by path induction is an equivalence:

$$(f = g) \to \big((a : A) \to (f(a) = f(g))\big)$$

This axiom is provable from the univalence axiom which we do not assume, but discuss later. We make ample use of this axiom, sometimes without explicitly invoking it. This is because in `rzk`, one always has to introduce the axiom per file and explicitly use it in every statement that depends on it. So there is no confusion which statements require it and which do not.

The next equivalences that we are going to cover are trivial to prove and often glossed over in practice, since they occur so often. They allow us to move arguments in functions around, swap elements in Σ-types and describe how to exchange Σ and Π types. The first is called the "axiom of choice" in classical mathematics. In our theory it is a statement that we can prove. This is because of the constructive nature of type theory. The sigma carries with it a proof along with the witness of the property. For this reason we are able to construct a choice function from the elementwise proof of existence. [1] also gives a more sensible type-theoretic axiom of choice that actually is an axiom. For this reason we refer to the following as the "theorem of choice".

**Proposition 2.5.7:** Theorem of choice (cf. [1, Theorem 2.15.7])

Given $B : A \to \mathcal{U}$ and $C : (a : A, b : B(a)) \to \mathcal{U}$, the following map is an equivalence:

$$\big((a : A) \to \Sigma(b : B(a), C(a, b))\big) \to \sum\big(b : (a : A) \to B, (a : A) \to C(a, b(a))\big)$$
$$F \mapsto (\pi_1 \circ F, \pi_2 \circ F)$$

**Proposition 2.5.8:** Associativity of $\Sigma$-types

Given $B : A \to \mathcal{U}, C : \Sigma B \to \mathcal{U}$, the canonical map between the two types below is an equivalence:

$$\Sigma C \equiv \Sigma(x : \Sigma B, C(x)) \equiv \Sigma(x : \Sigma(a : A, B(a)), C(x)) \simeq \Sigma(a : A, b : B(a), C(a, b))$$

**Proposition 2.5.9:** Commutativity of $\Sigma$-types

Given $B, C : A \to \mathcal{U}$, the canonical map between the two types below is an equivalence:

$$\Sigma(a : A, b : B(a), c : C(a)) \simeq \Sigma(a : A, c : C(a), b : B(a))$$

**Proposition 2.5.10:** Currying of $\Pi$-types

Given $B : A \to \mathcal{U}, C : \Sigma B \to \mathcal{U}$, the canonical map between the two types below is an equivalence:

$$\big((a : A, b : B(a)) \to C(a, b)\big) \simeq (x : \Sigma B) \to C(x)$$

Before we study the transport equivalence, we give a useful lemma to proving that two dependent propositions are equivalent without requiring fiberwise equivalence:

**Lemma 2.5.11**

Given two types $A, A' : \mathcal{U}$ and propositional families over each $B : A \to \mathbb{P}, B' : A' \to \mathbb{P}$. If we additionally have:

- an $f : A \to A'$ with $F : (a : A, b' : B'(f(a))) \to B(a)$ and
- a $g : A' \to A$ with $G : (a' : A', b : B(g(a'))) \to B'(a')$.

Then $\Pi B \simeq \Pi B'$.

**Proof.** Since both $\Pi B$ and $\Pi B'$ are propositions, it suffices to give maps back and forth:

$$\Pi B \to \Pi B' \qquad\qquad \Pi B' \to \Pi B$$
$$b \mapsto \lambda a'.\ G(a', b(g(a'))) \qquad\qquad b' \mapsto \lambda a.\ F(a, b'(f(a)))$$

□

**Transport**

As our last important equivalence, we turn to *transport*. It is an equivalence that relates the fibers of type families with each other given a path between the value the fiber is over.

**Definition 2.5.12**

Given a type family $B : A \to \mathcal{U}$ and a path in the base $p : x =_A y$, we define the transport of $p$ as the map $\mathsf{tr}^p_B : B(x) \to B(y)$ defined by path induction on $p$, sending $\mathsf{refl}_x$ to $\mathrm{id}_{B(x)}$.

Transport satisfies many useful properties, it is an equivalence with the inverse being transport along the reversed path. The composition of two transports is the transport of the concatenation of the paths and pullbacks turn into composing with Definition 2.2.7. For the details, refer to [1, Lemma 2.3.9-11]

## 2.6 Embeddings

An embedding is a weakened version of an equivalence. Instead of requiring a bijection on connected components along with compatible groupoidal structure, an embedding is only injective on connected components along with compatibility of groupoidal structures.

**Definition 2.6.1:** Embedding

A function $f : A \to B$ is called an embedding if the following holds:

$$\text{is-emb}(f) :\equiv \text{is-equiv}\left(\mathsf{ap}_f\right)$$

We later introduce the notion of a full embedding which corresponds with fully faithful functors from category theory.

## 2.7 Univalence

In this subsection, we discuss the *univalence axiom*. Because E. Riehl *et al.* [3] do not assume it and all of our proofs do not require it, we do not add it as an axiom to our theory. Nevertheless, it is an integral part of homotopy type theory and there are attempts at adding a directed version of univalence [24] to type theory.

The statement of the axiom requires the notion of equivalence and is rather simple to write down:

**Axiom 2.7.1**

For any two types $A, B : \mathcal{U}$ we have that the canonical map of equalities between $A$ and $B$ to equivalences between the two defined by path induction is an equivalence:

$$(A = B) \xrightarrow{\simeq} (A \simeq B)$$

This axiom is what allows us to treat *isomorphic* objects as *the same* not only from an intuitive perspective, but also directly from the formalism of our theory! It states that all equivalences arise from paths in $\mathcal{U}$.

Univalence is incompatible with the classical law of excluded middle and the law of double negation. The latter is shown in [1, Theorem 3.2.2] and the former in [1, Corollary 3.2.7]. We are however allowed to assume both axioms for propositions. But we again do not need to resort to using them, so we omit them from our theory.

## 2.8 Grothendieck Construction

Now we discuss the very interesting construction of straightening and unstraightening in our theory. For this subsection, we need to assume the univalence axiom. With it, we are able to prove the straightening theorem from higher category theory (albeit only for groupoids, since we do not yet have *directionality* in our theory).

Before we get to that, we define several notions and auxiliary statements that do not yet require univalence.

**Definition 2.8.1:** Fiber

Given a map $f : A \to B$, we define its fiber:

$$\mathrm{fib}_f(b) :\equiv \Sigma(a : A, f(a) = b)$$

**Lemma 2.8.2**

Let $B : A \to \mathcal{U}$ be a family, the fiber of the projection $\pi_B$ at $a$ is equivalent to the value of the family at $a$ via this map:

$$B(a) \to \mathrm{fib}_{\pi_B}(a) \quad b \mapsto ((a, b), \mathsf{refl}_a)$$

**Proof.** The inverse of the equivalence is given by $\lambda(x : \Sigma B, p : \pi_b(x) = a).\ \mathsf{tr}_B^p(\pi_2(x))$. The homotopy at $B(a)$ is $\mathsf{refl}_b$ and the homotopy at $\mathrm{fib}_{\pi_B}(a)$ is given by path induction over $p$. □

For this reason, we also call $B(a)$ the fiber of the family $B$.

**Lemma 2.8.3**

Let $f : A \to B$ be any map, then the following map is an equivalence:

$$\Sigma\mathrm{fib}_f \equiv \Sigma\big(b : B, \mathrm{fib}_f(b)\big) \to A$$
$$(b, (a, p)) \mapsto a$$

The Grothendieck construction shows that type families over a type $B$ are equivalent to fibrations over that type $B$ (where fibrations are just type maps). For this we introduce the shorthand notations:

$$\mathrm{Fam}(B) :\equiv B \to \mathcal{U} \qquad \mathrm{Fib}(B) :\equiv \Sigma(A : \mathcal{U}, A \to B)$$

**Theorem 2.8.4:** Straightening (cf. [1, Theorem 4.8.3], [5, Theorem 2.5.1])

For any type $B : \mathcal{U}$, the (un)straightening maps form an equivalence:

$$\mathrm{Un}_B : \mathrm{Fam}(B) \to \mathrm{Fib}(B) \qquad \mathrm{St}_B : \mathrm{Fib}(B) \to \mathrm{Fam}(B)$$
$$P \mapsto (\Sigma P, \pi_P) \qquad (A, f) \mapsto \mathrm{fib}_f$$

This construction still makes sense when one does not assume the univalence axiom; it just will not be an equivalence. Despite us not assuming it, we make heavy use of this construction, since it allows us to formulate every notion using type families. These are much more well-behaved from a notational perspective, as they often lead to *strict* versions of their fibration counterparts. Strict meaning that the homotopy that occurs in the fibration counterpart can be chosen to be refl.

**Remark 2.8.5**

This Grothendieck construction is the faithful representation of the version from $\infty$-groupoids. We later discuss the same construction for $\infty$-categories and note that there the story is not successful.

As our last definition of this section, we have composable type families:

**Definition 2.8.6**

Given a type family $B : A \to \mathcal{U}$ and a second family over the first, so $C : \Sigma B \to \mathcal{U}$, we can form the composite of $B$ and $C$:

$$\begin{aligned} C \odot B &: A \to \mathcal{U} \\ a &\mapsto \Sigma(b : B(a), C((a, b))) \end{aligned}$$

This corresponds with normal function composition when applying the Grothendieck construction.

In the next section we move onto *simplicial type theory*, which allows us to discuss $\infty$-categories and not only $\infty$-groupoids.

# 3 Simplicial HoTT and Category Theory

In the previous section, we introduced homotopy type theory and its core concepts. Now we want to add *directionality* to our theory. For this, we study the *simplicial homotopy type theory* by E. Riehl and M. Shulman [2]. It is a multi-layered type theory, each layer can only be dependent on itself or previous layers. There are three layers, the first two are *cubes* and *topes*. These are used to define *shapes*, which can then be used to map into classical types.

For example, the type of arrows $\mathrm{arr}(A)$ in a type $A$ is defined as mapping the 1-dimensional cube into $A$. The directionality in this theory comes from the 1-dimensional cube, it consists of two elements $0, 1$ with an ordering $0 < 1$.

As we already covered in the introduction, this approach allows us to keep all $\Pi$-types with the drawback that we have types that aren't categories. We thus have to introduce an internal notion of category in our theory, which is the notion of a *Rezk type*.

We introduce several category theoretic concepts in this section. Starting with exhibiting the usual notion of a category with objects, morphisms and their compositions. We then define isomorphisms, initial and final objects, fully faithful functors and adjunctions. Additionally, we analyze how squares decompose into triangles and how they relate to natural transformations.

Lastly, we observe that the universe $\mathcal{U}$ is not the category of categories. This fails already on the level of objects which are types and not necessarily *Rezk*. On top of that, any definition of $\mathrm{Cat}$ will not be *Segal*, making composition of morphisms in $\mathrm{Cat}$ impossible. Additionally, the Grothendieck construction is not faithfully represented.

While we give many of the inference rules of the theory, we will not be exhaustive and choose more concise, but potentially ambiguous notations. In particular, we omit $A$ **type** judgments and other distracting notation from our rules. We do this to more quickly introduce the theory in a simpler way without confusing the reader with unintelligible formalisms. If any questions arise, we direct the reader to E. Riehl and M. Shulman [2] or the formalization by E. Riehl *et al.* [3]. Since our work is formalized using `rzk` [4], we do not have any formal issues with being more lax.

## 3.1 Cubes, Topes and Shapes

The theory consists of three layers: *cubes*, *topes* and *types*. The type layer is the same as we have discussed before with a single addition: *extension types*. They make use of *shapes* which we define using the previous two layers.

Each of the layers specifies a context, we use the same symbol for entailment "$\vdash$" in all layers and collect the contexts on the left side as usual. To separate the contexts, we use a vertical line. We always use $\Xi$ for the context in the cube layer, $\Phi$ for the context in the tope layer and $\Gamma$ for the type layer. So a judgment with entailment using all contexts looks like this: $\Xi \mid \Phi \mid \Gamma \vdash \mathbb{1}\ \mathbf{type}$. We often omit this context where possible and apply local rule form even more aggressively. The cube context is most of the time the same in all judgments of a rule, so we generally omit it.

Cubes can be thought of as the synthetic equivalent of cubes in $\mathbb{R}^n$. Topes are essentially just booleans, dependent topes then give us a way to select vertices, edges, faces, and higher dimensional structures from within the cubes. We call these dependent topes *shapes*. Four basic examples of shapes that we can form are given in Figure 1. We call the first one $\Delta^1$, the second one is $\Lambda^2_0$, the third is $\Delta^2$ and the fourth one is $\Delta^3$.

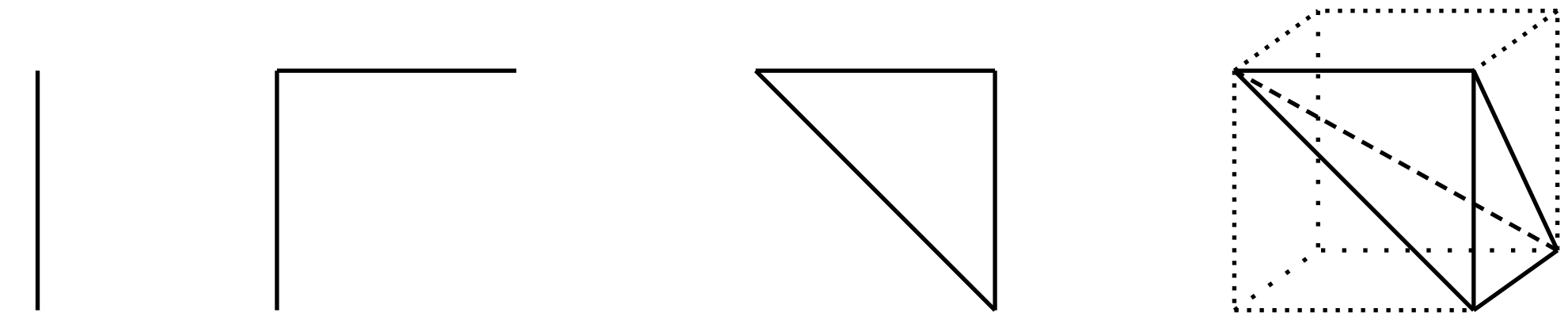

$t : \mathbb{2} \vdash \top$ $\quad (t, s) : \mathbb{2} \times \mathbb{2} \vdash t \equiv 0 \vee s \equiv 0$ $\quad (t, s) : \mathbb{2} \times \mathbb{2} \vdash s \leq t$ $\quad ((t, s), r) : (\mathbb{2} \times \mathbb{2}) \times \mathbb{2} \vdash s \leq t \wedge r \leq s$

Figure 1: Basic shapes and their representing judgments

In this section we do not go into all the small details of the theory described in [2]. We do not study how we have to phrase weakening in the various layers or contraction, we also do not list other technical rules like substitution. However, we still give an overview of the elements of this type theory that are not as "obvious" as these technical rules. In addition, we directly introduce the full theory in one go. Riehl and Shulman opt to first study "type theory with shapes" before they introduce the interval type and comparison tope.

### 3.1.1 Cubes

The first layer consists of *cubes*. We have the primitive 0-dimensional cube and a 1-dimensional one. To obtain higher dimensional ones, we allow finite products. The context in the cube layer is called $\Xi$. We also allow local rule forms, but only if the contexts are exactly the same in all cases (including all contexts of our theory). In the cube layer we only have access to the cube context, so most of our rules here can be stated in the local form.

**Definition 3.1.1:** Cube layer

The *cube layer* consists of the following syntax and rules: $I$ **cube** is the notation that we use to declare that $I$ is a cube. $I \times J$ is the notation used for products. $\mathbb{2}$ is the 1-dimensional cube and $\mathbf{1}$ the 0-dimensional one. Note that $\mathbf{1}$ is *different* from $\mathbb{1}$ which is the unit type from the third layer (also don't confuse both of them with 1, which is an element of $\mathbb{2}$).

$$\frac{(t : I) \in \Xi}{\Xi \vdash t : I} \qquad \frac{}{\mathbf{1}\ \textbf{cube}} \qquad \frac{}{\star : \mathbf{1}} \qquad \frac{}{\mathbb{2}\ \textbf{cube}} \qquad \frac{}{0 : \mathbb{2}} \qquad \frac{}{1 : \mathbb{2}}$$

$$\frac{I\ \textbf{cube} \quad J\ \textbf{cube}}{I \times J\ \textbf{cube}} \qquad \frac{s : I \quad t : J}{(s, t) : I \times J} \qquad \frac{t : I \times J}{\pi_1(t) : I} \qquad \frac{t : I \times J}{\pi_2(t) : J}$$

Note that $\pi$ is overloaded, as it appears in both the cube and the type layer.

In `rzk`, the unit cube and directed interval cube both exist along with forming products of cubes. Additionally, `CUBE` is used for the "universe of cubes" and then `I : CUBE` instead of $I$ **cube**. Pattern matching is also supported for elements of product cubes.

The fact that $\mathbb{2}$ is an interval is represented in the tope layer. Judgmental equality of cube elements also resides there. We are now finished with the first layer of our theory and can move onto the second.

### 3.1.2 Topes

The tope layer consists of *topes*, these essentially are “cube-dependent booleans” and behave like “strict propositions”. There are constants for the false and true topes, with the appropriate elimination and introduction rules (true can always be shown and from false one can derive anything). Additionally, we have conjunction and disjunction. Importantly, we do not have negation or implication. In this layer, we also add judgmental equality for cubes, which is going to be a tope so the next layer can use it as a primitive notion.

The tope context uses the letter $\Phi$. From now on when we write a rule in local rule form, one needs to take care to infer the correct context. For example, the tope declaration $X$ $\textbf{tope}$ only is allowed to have a cube context, whereas tope disjunction $\varphi \vee \psi$ has both a cube and a tope context.

**Definition 3.1.2:** Tope layer

The tope layer includes many more constructions than the cube layer. It consists of (cube dependent) disjunction, conjunction and equality. We begin with the simple constant terms:

$$\frac{}{\top\ \textbf{tope}} \qquad \frac{}{\top} \qquad \frac{}{\bot\ \textbf{tope}} \qquad \frac{\bot}{\psi}$$

Next we define tope-conjunction:

$$\frac{\varphi\ \textbf{tope} \quad \psi\ \textbf{tope}}{(\varphi \wedge \psi)\ \textbf{tope}} \qquad \frac{\varphi \quad \psi}{\varphi \wedge \psi} \qquad \frac{\varphi \wedge \psi}{\varphi} \qquad \frac{\varphi \wedge \psi}{\psi}$$

And disjunction:

$$\frac{\varphi\ \textbf{tope} \quad \psi\ \textbf{tope}}{(\varphi \vee \psi)\ \textbf{tope}} \qquad \frac{\varphi}{\varphi \vee \psi} \qquad \frac{\psi}{\varphi \vee \psi} \qquad \frac{\Phi, \varphi \vdash \chi \quad \Phi, \psi \vdash \chi \quad \Phi \vdash \varphi \vee \psi}{\Phi \vdash \chi}$$

We also have a tope for cube-equality. Which essentially is judgmental equality for cubes. We give it the same symbol that is used by judgmental equality in the type layer, since we can only compare elements of cubes or elements of types, there is no room for ambiguity. It also includes the $\eta$ and $\beta$ rules for the product cube, $\mathbf{1}$ and $\mathbb{2}$.

$$\frac{s : I \quad t : I}{(s \equiv t)\ \textbf{tope}} \qquad \frac{\Phi \vdash s : I}{s \equiv s} \qquad \frac{s \equiv t}{t \equiv s} \qquad \frac{s \equiv t \quad t \equiv v}{s \equiv v} \qquad \frac{t : \mathbf{1}}{t \equiv \star} \qquad \frac{0 \equiv 1}{\bot}$$

$$\frac{s : I \quad t : J}{\Phi \vdash \pi_1((s,t)) \equiv s} \qquad \frac{s : I \quad t : J}{\Phi \vdash \pi_2((s,t)) \equiv t} \qquad \frac{t : I \times J}{\Phi \vdash (\pi_1(t), \pi_2(t)) \equiv t}$$

$$\frac{\Xi, x : I \vdash \psi\ \textbf{tope} \quad \Xi \mid \Phi \vdash t \equiv s \quad \Xi \mid \Phi \vdash \psi[x \mapsto s]}{\Xi \mid \Phi \vdash \psi[x \mapsto t]}$$

In addition to equality, we also have to give $\mathbb{2}$ an ordering through a comparison tope. We give the axioms of this comparison tope in equation notation rather than inference rules, as they are much easier to understand. In the cube context $x : \mathbb{2}, y : \mathbb{2}$, we declare $(x \leq y)$ $\textbf{tope}$ called the *inequality tope*. We give the axioms of the strict interval: for any $x, y, z : \mathbb{2}$ we have:

$$(x \leq x) \qquad (0 \leq x) \qquad (x \leq 1)$$

$$(x \leq y), (y \leq z) \vdash (x \leq z) \qquad (x \leq y), (y \leq x) \vdash (x \equiv y) \qquad (x \leq y) \vee (y \leq x)$$

These rules are called from left to right and top to bottom: reflexivity, minimum, maximum, transitivity, anti-symmetry and totality.

The strict interval cube $\mathbb{2}$ only declares that 0 and 1 are elements of it. However, the axioms do not forbid further elements in between them. In fact this is what makes the theory interesting, as one of the valid models of this theory uses $[0, 1] \subset \mathbb{R}$ (or $[0, 1] \cap \mathbb{Q}$) as the interval. The intuition for $\mathbb{2}$ should be that of a continuous interval.

Similar to cubes, `rzk` has a tope universe called `TOPE`. A shape is now a "function" `I → TOPE` where `I` is a cube. This is also interpreted as a function, but is in no relation to Π-types from the next layer. Parameters of these "tope functions" can also be tope functions themselves. This is how one can define subshapes, given `ψ : I → TOPE`, we can declare a subshape `φ : ψ → TOPE`. The constants $\top$, $\bot$ as well as conjunction, disjunctions, equality and inequality topes are also available.

### 3.1.3 Shapes

As we already explained earlier, shapes are specific cube-dependent topes, which we use to define terms and types depending on them. Concretely, a shape is a tope in a *singleton* context. We use a set-like notation for shapes. [2] mentions that we could formalize the notion by adding the following rule:

$$\frac{I \textbf{ cube} \qquad t : I \vdash \varphi \textbf{ tope}}{\{t : I \mid \varphi\} \textbf{ shape}}$$

One very important class of shapes are simplices:

$$\Delta^n := \{(t_1, ..., t_n) : \mathbb{2}^n \mid t_n \leq ... \leq t_1\}$$

We only need this shape for small $n$.

In addition to simplices, we define the notion of *subshape*: we write $\varphi \subset \psi$ to denote that $\varphi$ is a subshape of $\psi$. Formally, a subshape relation is declared as follows:

$$\frac{\{t : I \mid \psi\} \textbf{ shape} \qquad \{t : I \mid \varphi\} \textbf{ shape} \qquad t : I \mid \varphi \vdash \psi}{\varphi \subset \psi}$$

So both shapes must be declared in the same cube context and whenever the subshape picks out a point, the supershape also has to include it.

As one would expect, we can define the inner and outer horns of the simplices. Since we only need low-dimensional ones, we only give the most important one:

$$\Lambda^2_1 := \{(t, s) : \mathbb{2}^2 \mid t \equiv 1 \lor s \equiv 0\}$$

Another subshape that is rather important is the boundary of $\Delta^1$:

$$\partial\Delta^1 := \{t : \mathbb{2} \mid t \equiv 0 \lor t \equiv 1\}$$

**Remark 3.1.3**

U. Buchholtz and J. Weinberger [5] reinterpret shapes as types, which proves to be advantageous in several ways. For example, one is able to define maps between shapes which one can then use to define maps on types. $\Lambda^2_1$ is a subshape of $\Delta^2$ and postcomposition of the inclusion map

$$(\Delta^2 \to A) \to (\Lambda^2_1 \to A)$$

plays an important role when we later talk about Segal types. In our interpretation where we do not allow shapes to be coerced to types, we are not able to talk about this map as a postcomposition. Instead, we have to define it like this:

$$\lambda f \colon \Delta^2 \to A.\ \lambda t \colon \Lambda^2_1.\ f(t)$$

While this is straightforward to define, there are other cases that are more complicated.

Our reason for keeping shapes and types separate notions is that `rzk` also has this behavior and we could not formalize any notion that would rely on such a coercion existing.

**Shape Operations**

From the basic simplices, we can form more complex shapes by using the following operations on shapes: union, product and subshape product.

**Definition 3.1.4**

Given two shapes $\{t : I \mid \psi\}$ and $\{t : I \mid \zeta\}$ over the same cube, we can form their union:

$$\psi \cup \zeta := \{t : I \mid \psi \vee \zeta\}$$

From two shapes, we can form their product by giving a new shape over the product of their cubes.

**Definition 3.1.5**

Given two shapes $\{t : I \mid \psi\}$ and $\{s : J \mid \zeta\}$ over the cubes $I$ and $J$, we define their product as:

$$\psi \times \zeta := \{(t, s) : I \times J \mid \psi \wedge \zeta\}$$

Since we often consider subshape inclusions, we also need to be able to form a relative product of two subshapes in the product of their supershapes. This relative product is different from the actual product of the two subshapes, because it does not have the right property when one wants to use it in the currying theorem for extension types.

**Definition 3.1.6**

Given two shape inclusions $\varphi \subset \{t : I \mid \psi\}$ and $\chi \subset \{s : J \mid \zeta\}$ we define the *pushout product* of $\varphi$ and $\chi$ (we also call this the *subshape product*):

$$\begin{aligned} \varphi \otimes \chi &:= (\varphi \times \zeta) \cup (\psi \times \chi) \\ &\equiv \{(t, s) : I \times J \mid (\varphi \wedge \zeta) \vee (\psi \wedge \chi)\} \end{aligned}$$

**Example 3.1.7**

A very commonly used product shape is $\Delta^1 \times \Delta^1$. Its tope condition is not interesting though since it is just $\top$. However, if we consider $\partial\Delta^1 \subset \Delta^1$ in both cases instead, we get an interesting subshape product:

$$\partial\Delta^1 \otimes \partial\Delta^1 \equiv \{(t, s) : \mathbb{2} \times \mathbb{2} \mid t \equiv 0 \vee t \equiv 1 \vee s \equiv 0 \vee s \equiv 1\}$$

Which is the full boundary of the square.

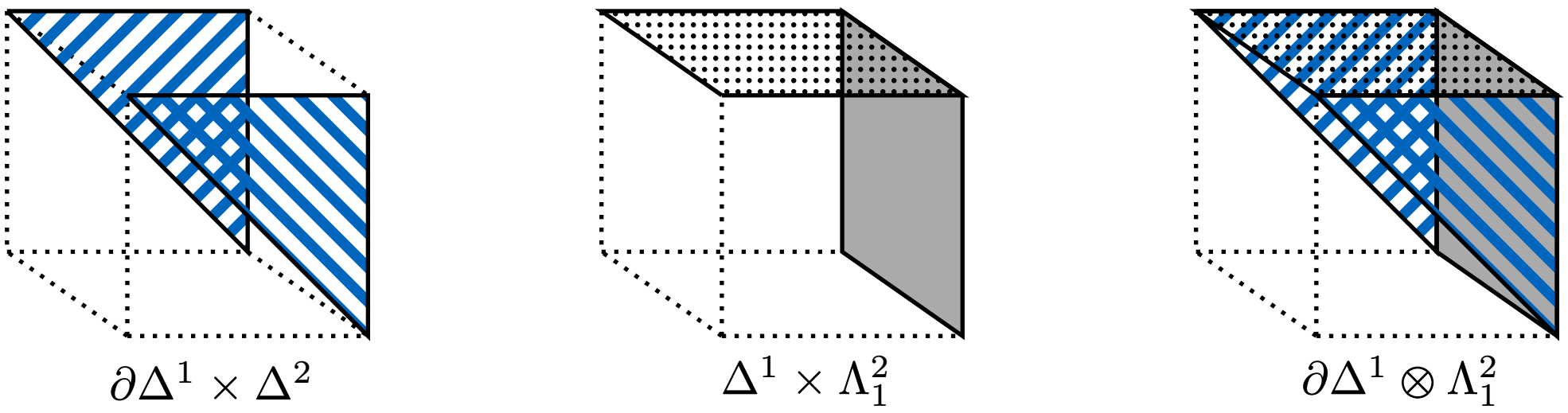

Figure 2: Decomposition of $\partial\Delta^1 \otimes \Lambda^2_1$

A more complex example is $\partial\Delta^1 \subset \Delta^1$ and $\Lambda^2_1 \subset \Delta^2$ shown in Figure 2. The shape $\Delta^1 \times \Delta^2$ is the entire prism (the dotted lines are the rest of the cube in which this shape lives). All faces except the diagonal "underside" of the prism are included in the subshape product.

## 3.2 Shape-Dependent Terms and Extension Types

Our theory now has the types that we introduced in Section 2 as the third layer, which means that any type can be dependent on shapes. Since our study in the previous section did not have the shape notion introduced, there is not yet any type constructor that allows us to use shapes. E. Riehl and M. Shulman [2] invent the notion of *extension type*, which makes use of shapes. It essentially is a Π-type where the domain is a shape and the codomain a shape-dependent type. In addition, we can specify a constraint on a subshape of the domain where we want the functions to take certain values in the type. Additionally, we add type-layer recursors for tope-disjunction and the ⊥-tope. These allow us to define shape-dependent types that we further use to define concrete extension types.

**Definition 3.2.1:** Type-recursors for tope constructions

Since $\bot$ corresponds to false from logic, we add the same recursor as for the empty type $\mathbb{0}$ we discussed before. The only difference is that instead of requiring an element $x : \mathbb{0}$, we have $\bot$ in the premises:

$$\frac{\bot}{\mathrm{rec}_\bot : A} \qquad \frac{\bot \quad a : A}{\mathrm{rec}_\bot \equiv a}$$

The type-recursor for tope-disjunction has these four rules. $\mathrm{rec}_\vee$ is essentially an `if` with the addition that the two values for the branches must agree on the conjunction of the topes.

$$\frac{\varphi \vee \psi \quad \varphi \vdash x : A \quad \psi \vdash y : A \quad \varphi \wedge \psi \vdash x \equiv y}{\mathrm{rec}_\vee^{\varphi,\psi}(x, y) : A} \qquad \frac{\varphi \vee \psi \quad a : A}{\mathrm{rec}_\vee^{\varphi,\psi}(a, a) \equiv a}$$

$$\frac{\varphi \vee \psi \quad \varphi \vdash x : A \quad \psi \vdash y : A \quad \varphi \wedge \psi \vdash x \equiv y}{\varphi \vdash \mathrm{rec}_\vee^{\varphi,\psi}(x, y) \equiv x} \qquad \frac{\varphi \vee \psi \quad \varphi \vdash x : A \quad \psi \vdash y : A \quad \varphi \wedge \psi \vdash x \equiv y}{\psi \vdash \mathrm{rec}_\vee^{\varphi,\psi}(x, y) \equiv y}$$

To make our notation shorter, we also allow writing $x \vee y$ instead of $\mathrm{rec}_\vee^{\varphi,\psi}(x, y)$ when the shapes are obvious from context.

Since we use the recursor of $\partial\Delta^1$ (which decomposes into $t \equiv 0 \vee t \equiv 1$) very often, also with values that themselves are already $\vee$-recursors, we give it a special notation:

$$\mathrm{rec}_\vee^{t\equiv 0,\, t\equiv 1}(x, y) \equiv [x, y]$$

We will also use this notation for $\Lambda_1^2$, which we define later.

Using these recursors, we can now define shape-dependent types and terms. To fully lift these to the type-layer, we also need to allow turning shape-dependent types into proper types. For this reason we introduce extension types.

**Definition 3.2.2:** Extension type

Because we have to declare two shapes and their subshape relationship in our rules for this type, they become quite long and might seem cryptic at first glance. But the contents of an extension type are exactly the same as a Π-type. Given two shapes $\varphi$ and $\psi$ over the cube $I$ with $\varphi \subset \psi$, a shape-dependent type $t : I \mid \psi \vdash A$ **type** and a shape-dependent term $t : I \mid \varphi \vdash a : A$ that is declared over the subshape, we write

$$\langle (t : \psi) \to A \mid {}^{\varphi}_{a} \rangle$$

for the extension type.

An extension type is like a Π-type, with two major differences. The first is that instead of the domain being a type, it is a shape. The second difference is that an extension type specifies a *constraint* on a subshape. All functions of an extension type evaluate to the constraint when plugging in a value that is part of the subshape.

If we had this constructor also for normal Π-types, then we could form functions like "a function from ℝ to ℝ that is 0 on negative inputs".

Formally, we have the following formation rule:

$$\frac{\{t : I \mid \psi\}\ \textbf{shape} \qquad \{t : I \mid \varphi\}\ \textbf{shape} \qquad t : I \mid \varphi \vdash \psi \qquad t : I \mid \psi \vdash A\ \textbf{type} \qquad t : I \mid \varphi \vdash a : A}{\langle (t : \psi) \to A \mid {}^{\varphi}_{a} \rangle\ \textbf{type}}$$

We note that in this case, writing the rule with no explicit contexts leads to confusion. The shapes $\psi$ and $\varphi$ are **not** allowed to depend on any external context. $A$ and $a$ are allowed to do so. Because this behavior is not found in `rzk` and because we do not make use of it, omitting the context will not get us into trouble.

The introduction rule uses almost the same premises, the only addition is a dependent term over the supershape $\psi$ that matches the constraint on the subshape (imposed strictly via judgmental equality). From now on we also omit the shape declaration judgments as well as the explicit cube variable $t : I$ to save space.

$$\frac{\varphi \subset \psi \qquad \psi \vdash A\ \textbf{type} \qquad \varphi \vdash a : A \qquad \psi \vdash b : A \qquad \varphi \vdash a \equiv b}{\lambda t.\ b : \langle (t : \psi) \to A \mid {}^{\varphi}_{a} \rangle}$$

We also have the usual $\beta$ and $\eta$ rules that Π-types have:

$$\frac{\varphi \subset \psi \qquad \psi \vdash A\ \textbf{type} \qquad \varphi \vdash a : A \qquad f : \langle (t : \psi) \to A \mid {}^{\varphi}_{a} \rangle}{f \equiv \lambda t.\ f(t)}$$

The $\eta$ rule requires us to talk about inputs, so we need to reintroduce the explicit cube variable for this rule:

$$\frac{t : I \mid \varphi \vdash \psi \quad t : I \mid \psi \vdash A \textbf{ type} \quad t : I \mid \varphi \vdash a : A \quad t : I \mid \psi \vdash b : A \quad t : I \mid \varphi \vdash a \equiv b \quad s : I \quad \psi[t \mapsto s]}{(\lambda t.\ b)(s) \equiv b[t \mapsto s]}$$

Additionally, we also have a rule that exhibits the strict nature of the extension type when we plug in a cube element from the subshape:

$$\frac{t : I \mid \varphi \vdash \psi \quad t : I \mid \psi \vdash A \textbf{ type} \quad t : I \mid \varphi \vdash a : A \quad f : \left\langle (t : \psi) \to A \,\middle|\, {}^{\varphi}_{a} \right\rangle \quad s : I \quad \varphi[t \mapsto s]}{f(s) \equiv a[t \mapsto s]}$$

We also use several shorthand notations. Firstly, instead of writing the tope in the extension type, we also allow writing the entire shape. For example, we write $\left\langle (t : \Delta^1) \to A \,\middle|\, {}^{\partial \Delta^1}_{a} \right\rangle$. When $a$ and $A$ are not depending on $t$, then we allow omitting the $t$ variable: $\left\langle \psi \to A \,\middle|\, {}^{\varphi}_{a} \right\rangle$. When $a$ depends on $t$, but $A$ doesn't, we also allow $\left\langle \psi \to A \,\middle|\, {}^{\varphi}_{\lambda t.\ a} \right\rangle$.

When $\varphi$ is $\bot$, we effectively don't have a subshape and thus allow writing $(t : \psi) \to A$ instead of $\left\langle (t : \varphi) \to A \,\middle|\, {}^{\bot}_{\mathrm{rec}_\bot} \right\rangle$. In the case that $A$ isn't depending on $t$, we also write $\psi \to A$.

Additionally to our notation, we allow certain coercions via judgmental equalities:

$$\frac{f : \left\langle (t : \psi) \to A(t) \,\middle|\, {}^{\varphi}_{a(t)} \right\rangle}{f : (t : \psi) \to A(t)} \qquad \frac{f : \left\langle (t : \chi) \to A(t) \,\middle|\, {}^{\varphi}_{a(t)} \right\rangle}{f : \left\langle (t : \psi) \to A(t) \,\middle|\, {}^{\varphi}_{a(t)} \right\rangle}$$

The first rule allows us to forget any constraint and the second rule allows us to restrict any function to a subshape of its domain.

#### Remark 3.2.3

`rzk` also has the two recursors, but uses a slightly different syntax for $\mathrm{rec}_\vee$:

```
recOR(«tope_1» ↦ «term_1», ..., «tope_n» ↦ «term_n»)
```

This term is well-typed when all intersections judgmentally agree. Extension types exist as a special case of Π-types by using the *restriction type constructor*. Using the same syntax as `recOR`, one can write:

```
«type» [«tope_1» ↦ «term_1», ..., «tope_n» ↦ «term_n»]
```

This term corresponds with the term `«type»` when all topes are false and otherwise the singleton type containing the element for which the tope is true. While this is a departure from the original paper by [2], there is some literature on including these types directly by J. Sterling and C. Angiuli [25, Notation 3] and D. Gratzer, J. Sterling, C. Angiuli, T. Coquand, and L. Birkedal [26, Section 3.1]. These types clearly allow more definitions and their exact behavior in this theory is not yet fully understood. We do not make use of them directly.

Using extension types, we can define our notion of morphism in a type, as well as all other category related definitions.

#### Definition 3.2.4: Morphisms

The type of morphisms between $x$ and $y$ in $A$ is defined as:

$$\hom_A(x, y) := \left\langle \Delta^1 \to A \,\middle|\, {}^{\partial \Delta^1}_{[x,y]} \right\rangle$$

We define the type of arrows $\mathrm{arr}(A) :\equiv \Delta^1 \to A$. So $\hom_A(x, y)$ is just the type of arrows, restricted to those elements that start at $x$ and end at $y$.

Before we can define how composition works for these morphisms, we study extension types themselves from a type theoretic perspective. For this we give several important equivalences.

### 3.2.1 Equivalences of Extension Types

The similarity of Π-types and extension types is also reflected in their respective equivalences. We now give several propositions that also exist for Π-types. Riehl and Shulman already explore most of these and we only have a few new additions. We only give proofs for the new ones, as the already known equivalences are proven in [2].

#### Currying, Commuting, and the Theorem of Choice

Extension types allow for currying with both themselves and ordinary Π-types. The only complication is the constraint where we use the previously defined subshape product when currying two extension types.

**Proposition 3.2.5:** (cf. [2, Theorem 4.1])

Given a shape inclusion $\varphi \subset \psi$, a dependent type family $B : \psi \to A \to \mathcal{U}$ and a constraint $b : (t : \varphi, a : A) \to B(t, a)$, then there is an equivalence:

$$\left\langle (t : \psi) \to (a : A) \to B(t, a) \;\middle|\; {}^{\varphi}_{b(t)} \right\rangle \simeq (a : A) \to \left\langle (t : \psi) \to B(t, a) \;\middle|\; {}^{\varphi}_{b(t,a)} \right\rangle$$

**Proposition 3.2.6:** (cf. [2, Theorem 4.2])

Given two shape inclusions $\varphi \subset \psi$, $\chi \subset \zeta$, a type family $A : \psi \to \zeta \to \mathcal{U}$ and a shared constraint $a : ((t, s) : \varphi \otimes \chi) \to A(t, s)$, then there is an equivalence:

$$\begin{aligned} &\left\langle (t : \psi) \to \left\langle (s : \zeta) \to A(t, s) \;\middle|\; {}^{\chi}_{a(t,s)} \right\rangle \;\middle|\; {}^{\varphi}_{\lambda s.\ a(t,s)} \right\rangle \\ \simeq &\left\langle ((t, s) : \psi \times \zeta) \to A(t, s) \;\middle|\; {}^{\varphi \otimes \chi}_{a(t,s)} \right\rangle \\ \simeq &\left\langle (s : \zeta) \to \left\langle (t : \psi) \to A(t, s) \;\middle|\; {}^{\varphi}_{a(t,s)} \right\rangle \;\middle|\; {}^{\chi}_{\lambda t.\ a(t,s)} \right\rangle \end{aligned}$$

Extension types are also compatible with Σ-types, so similar to Proposition 2.5.7 for Π-types, we can observe the same equivalence for extension types:

**Proposition 3.2.7:** (cf. [2, Theorem 4.3])

Given a shape inclusion $\varphi \subset \psi$, a dependent type over the bigger shape $A : \psi \to \mathcal{U}$, a family over it $B : (t : \psi, a : A(t)) \to \mathcal{U}$ and a constraint function $f : (t : \varphi) \to \Sigma(a : A(t), B(t, a))$, the following map is an equivalence:

$$\begin{aligned} &\left\langle (t : \psi) \to \Sigma(a : A(t), B(t, a)) \;\middle|\; {}^{\varphi}_{f} \right\rangle \\ &\quad \to \sum \left( a : \left\langle (t : \psi) \to A(t) \;\middle|\; {}^{\varphi}_{\pi_1 \circ f} \right\rangle, \left\langle (t : \psi) \to B(t, a(t)) \;\middle|\; {}^{\varphi}_{\pi_2 \circ f} \right\rangle \right) \\ &F \mapsto (\pi_1 \circ F, \pi_2 \circ F) \end{aligned}$$

### Propositional Constraints

Now we observe that a propositional constraint is equivalent to a strict constraint. Before we can state that proposition, we give an auxiliary lemma that shows an unconstrained extension type is equivalent to a constraint on the subshape in addition to a constrained extension type.

**Lemma 3.2.8**

Given a shape inclusion $\varphi \subset \psi$ and a dependent type $A : \psi \to \mathcal{U}$, we have that the following map is an equivalence:

$$\big((t : \psi) \to A(t)\big) \to \sum\Big(f : (t : \varphi) \to A(t), \Big\langle (t : \psi) \to A(t) \;\Big|\; {}^{\varphi}_{f} \Big\rangle\Big)$$
$$F \mapsto (F, F)$$

**Proof.** We give a strict inverse $(f, F) \mapsto F$. □

**Proposition 3.2.9**

Given a shape inclusion $\varphi \subset \psi$, a dependent type $A : \psi \to \mathcal{U}$ and a constraint $f : (t : \varphi) \to A(t)$, we have that the following map is an equivalence:

$$\Big\langle (t : \psi) \to A(t) \;\Big|\; {}^{\varphi}_{f} \Big\rangle \to \sum\Big(g : \Pi A, g \underset{(t:\varphi)\to A(t)}{=} f\Big)$$
$$F \mapsto (F, \mathsf{refl})$$

We note the use of implicit coercions of $F$, $f$, and $g$ ($F$ loses its constraint and $f$ and $g$ lose their constraints and are restricted to $\varphi$).

**Proof.** We apply Lemma 3.2.8 to the codomain of the function and obtain the following equivalent type:

$$\sum\left(\begin{array}{c} g : \Sigma\Big(f' : (t : \varphi) \to A(t), \Big\langle (t : \psi) \to A(t) \;\Big|\; {}^{\varphi}_{f'} \Big\rangle\Big) \\ \pi_1(g) \underset{(t:\varphi)\to A(t)}{=} f \end{array}\right)$$

Using Proposition 2.5.8 and Proposition 2.5.9, we have that this type is equivalent to the following type:

$$\sum\left(\begin{array}{c} f' : (t : \varphi) \to A(t) \\ f' \underset{(t:\varphi)\to A(t)}{=} f \\ \Big\langle (t : \psi) \to A(t) \;\Big|\; {}^{\varphi}_{f'} \Big\rangle \end{array}\right)$$

Since based paths are contractible by Lemma 2.3.3, we obtain an equivalence to our domain. Inspecting this equivalence now reveals that it is the map we claimed to be an equivalence. □

### Extension Extensionality

Similarly to function extensionality, we need the analogous axiom for extension types.

**Axiom 3.2.10:** Extension Extensionality (cf. [2, Axiom 4.6])

Given a shape inclusion $\varphi \subset \psi$, a type family $A : \psi \to \mathcal{U}$ and a constraint $a : (t : \varphi) \to A(t)$, such that every $A(t)$ is contractible, we have that the following type is contractible:

$$\langle (t : \psi) \to A(t) \mid {}^{\varphi}_{a} \rangle$$

Most of the time we use this axiom for proving that two elements of an extension type are the same by using the following derived proposition:

**Proposition 3.2.11:** (cf. [2, Proposition 4.8])

Given a shape inclusion $\varphi \subset \psi$, a type family $A : \psi \to \mathcal{U}$ and a constraint $a : (t : \psi) \to A(t)$, moreover let $f, g : \langle (t : \psi) \to A(t) \mid {}^{\varphi}_{a} \rangle$. Then the following map is an equivalence if we also assume Axiom 3.2.10:

$$(f = g) \to \langle (t : \psi) \to f(t) = g(t) \mid {}^{\varphi}_{\mathsf{refl}} \rangle$$

In particular, if $f$ and $g$ are equal pointwise (with the proof being $\mathsf{refl}$ on the subshape) then we have $f = g$.

In `rzk` this axiom is also introduced like the other axioms using the `#assume` directive. All the statements that require it list it either explicitly in the `uses` list or directly in their body.

## 3.3 Categorical Structure

We have already defined arrows/morphisms as $\hom_A(x, y) \equiv \left\langle (t : \Delta^1) \to A \;\middle|\; {}^{\partial\Delta^1}_{[x,y]} \right\rangle$, but we have not yet exhibited the categorical structure of the type $A$. In this subsection we explore how composition works for arrows as well as introduce adjacent concepts.

**Definition 3.3.1:** Identity Arrow

We define the identity arrow $\mathsf{id}_x$ for any $x : A$ as $\mathsf{id}_x := \lambda t.x$

Classical higher category theory formulated in ZFC via simplicial sets defines composition of morphisms by requiring the existence of unique fillers. In our theory we also follow this definition. Types in which we can compose arrows are called Segal types.

### 3.3.1 Segal Types

**Definition 3.3.2**

Given a type $A$, three elements $x, y, z : A$ and three arrows $f : \hom_A(x, y)$, $g : \hom_A(y, z)$, $h : \hom_A(x, z)$ we define the type of triangles that have the boundary fixed at $f, g, h$:

$$\hom^2_A \left( \begin{array}{ccc} & y & \\ {}^{f}\nearrow & & \searrow^{g} \\ x & \xrightarrow[h]{} & z \end{array} \right) := \left\langle ((t, s) : \Delta^2) \to A \;\middle|\; \begin{array}{l} \overbrace{s \equiv 0 \vee t \equiv 1 \vee t \equiv s}^{\partial\Delta^2} \\ f \vee g \vee h \end{array} \right\rangle$$

We allow omitting any of the ingredients in our notation, which is then supposed to mean that one takes the $\Sigma$-type qualifying over all of the missing parts. For example:

$$\hom_A^2\left(\begin{array}{ccc} & y & \\ {}^{f}\nearrow & & \searrow^{g} \\ x & \longrightarrow & z \end{array}\right) \equiv \Sigma\left(h : \hom_A(x,z), \hom_A^2\left(\begin{array}{ccc} & y & \\ {}^{f}\nearrow & & \searrow^{g} \\ x & \xrightarrow[h]{} & z \end{array}\right)\right)$$

We note that a vertex can also be unspecified if the morphism already ascribes which value should be taken. In this case, we use $\cdot$ in the place of the value.

Using this definition we can write down what it means for a type to be Segal.

**Definition 3.3.3:** Segal type

We call a type $A : \mathcal{U}$ a *Segal type* if the following proposition holds:

$$(x : A, y : A, z : A, f : \hom_A(x,y), g : \hom_A(y,z))$$
$$\to \text{is-contr}\left(\hom_A^2\left(\begin{array}{ccc} & y & \\ {}^{f}\nearrow & & \searrow^{g} \\ x & \longrightarrow & z \end{array}\right)\right)$$

So for any two composable morphisms $f, g$ there exists a unique filler constructing a composite.

In the context of a Segal type, we overload the symbol $\circ$ to also apply to arrows. The value of $g \circ f$ is the first component of the $\Sigma$-type above. The second component guarantees that this value is unique over all possible compositions.

There are lots of interesting properties that can be proven about Segal types. We however only consider the following one, as we need it later.

**Lemma 3.3.4:** (cf. [2, Proposition 5.10])

For any $f, g : \hom_A(x,y)$ in a Segal type $A$, the natural maps defined by path induction:

$$(f = g) \to \hom_A^2\left(\begin{array}{ccc} & x & \\ {}^{\mathsf{id}_x}\nearrow & & \searrow^{f} \\ x & \xrightarrow[g]{} & y \end{array}\right) \qquad (f = g) \to \hom_A^2\left(\begin{array}{ccc} & y & \\ {}^{f}\nearrow & & \searrow^{\mathsf{id}_y} \\ x & \xrightarrow[g]{} & y \end{array}\right)$$

are equivalences.

### Functoriality of Functions

Similarly to Definition 2.2.7 from HoTT, we can define an action of functions on arrows and thus exhibit functions between Segal types as functors preserving the categorical structure.

**Definition 3.3.5**

Let $F : A \to B$ be a function. We define its effect on hom-types:

$$F_\# : \hom_A(x,y) \to \hom_B(F(x), F(y))$$
$$f \mapsto \lambda t.\ F(f(t))$$

**Theorem 3.3.6:** Functoriality (cf. [2, Proposition 6.1])

Let $F : A \to B$ and $G : B \to C$ be functions between Segal types. Moreover, let $x, y, z : A$ and $f : \hom_A(x,y)$, $g : \hom_A(y,z)$. Then the following hold:

$$F_\#(\mathsf{id}_x) = \mathsf{id}_{F(x)} \qquad F_\#(g \circ f) = F_\#(g) \circ F_\#(f) \qquad (G \circ F)_\# = G_\# \circ F_\#$$

Note that the composition in the second formula is composition of arrows and in the third formula composition of functions.

**Proof.** The first two are shown by [2]. The third property is easily derived by unfolding definitions and applying function and extension extensionality. □

### 3.3.2 Dependent Arrows

Since we are using type families, we also need the notion of a dependent arrow. This arrow lives over an arrow in the base.

**Definition 3.3.7:** Dependent arrow

Given a type family $B : A \to \mathcal{U}$, elements $x, y : A$, an arrow $f : \hom_A(x, y)$ and elements in the fibers $\overline{x} : B(x), \overline{y} : B(y)$ we define the type of dependent arrows in $B$ over $f$:

$$\mathrm{dhom}_B^f(\overline{x}, \overline{y}) :\equiv \left\langle (t : \Delta^1) \to B(f(t)) \;\middle|\; {}^{\partial\Delta^1}_{[\overline{x}, \overline{y}]} \right\rangle$$

Dependent arrows over the identity are just normal arrows in the fiber.

**Lemma 3.3.8**

Given $B : A \to \mathcal{U}$ and $a : A$, $x, y : B(a)$ we have:

$$\mathrm{dhom}_B^{\mathsf{id}_a}(x, y) \equiv \hom_{B(a)}(x, y)$$

**Proof.** Unfolding the definitions, we get:

$$\mathrm{dhom}_B^{\mathsf{id}_a}(x, y) \equiv \left\langle (t : \Delta^1) \to B(\mathsf{id}_a(t)) \;\middle|\; {}^{\partial\Delta^1}_{[x, y]} \right\rangle$$
$$\hom_{B(a)}(x, y) \equiv \left\langle \Delta^1 \to B(a) \;\middle|\; {}^{\partial\Delta^1}_{[x, y]} \right\rangle$$

Since $\mathsf{id}_a(t) \equiv a$, we obtain that they are judgmentally equivalent. □

One place where dependent arrows naturally occur is arrows in total types of families:

**Lemma 3.3.9**

Given a type $A : \mathcal{U}$ and a family $B : A \to \mathcal{U}$, then there is an equivalence for every $x, y : \Sigma B$:

$$\hom_{\Sigma B}(x, y) \to \sum\left(f : \hom_A(x_1, y_1), \mathrm{dhom}_B^f(x_2, y_2)\right)$$
$$f \mapsto (t \mapsto (\mathrm{pr}_1(f(t)), \mathrm{pr}_2(f(t))))$$

**Proof.** This is just an instance of Proposition 3.2.7. □

### 3.3.3 Isomorphisms and Rezk Types

The notion of isomorphism requires a Segal type, since an arrow is an isomorphism if it has a left and a right inverse with respect to composition.

**Definition 3.3.10:** Isomorphism

An arrow $f : \hom_A(x, y)$ in a Segal type is an isomorphism if the following holds:

$$\mathsf{isiso}_A(f) :\equiv \left(\sum(g : \hom_A(y, x), g \circ f = \mathsf{id}_x)\right) \times \left(\sum(h : \hom_A(y, x), f \circ h = \mathsf{id}_y)\right)$$

This definition is a proposition similar to how our definition of equivalence is a proposition and specifically must be formulated to have separate left and right inverses (see [2, Proposition 10.2] for a proof).

We write $\mathrm{Iso}_A(x, y) :\equiv \Sigma(f : \hom_A(x, y), \mathsf{isiso}_A(f))$ for the type of isomorphism.

**Remark 3.3.11**

As one expects, isomorphisms compose and $\mathsf{id}_x$ is an isomorphism. However, isomorphisms are a strictly category theoretic notion that exists in *addition* to the type theoretic notion of propositional equality. Being embedded in type theory, we of course get that propositional equality implies that the two equal objects are isomorphic. We can define this function via path induction:

$$\mathsf{idtoiso} : (x = y) \to \mathrm{Iso}_A(x, y) \qquad \mathsf{refl}_x \mapsto \mathsf{id}_x$$

However, isomorphic objects are not necessarily equal. For this reason we introduce Rezk types in which the above map is an equivalence.

**Definition 3.3.12:** Rezk Type

A Segal type $A$ is called *Rezk* if $\mathsf{idtoiso} : (x = y) \to \mathrm{Iso}_A(x, y)$ is an equivalence for any $x, y : A$.

Since all isomorphisms in a Rezk type come from an equality, we can transfer the induction principle to isomorphisms.

**Theorem 3.3.13:** Isomorphism Induction

Given a Rezk type $A$ with $x : A$ and a dependent type family:

$$C : (y : A, \mathrm{Iso}_A(x, y)) \to U$$

If we also have a value $c : C(x, \mathsf{id}_x)$, there is a value of the type

$$(y : A, f : \mathrm{Iso}_A(x, y)) \to C(y, f)$$

Functions preserve isomorphisms.

**Lemma 3.3.14**

Let $F : A \to B$ be a map between Rezk types. Then for any $f : \hom_A(x, y)$ there is a map:

$$\mathsf{isiso}_A(f) \to \mathsf{isiso}_B(F_\#(f))$$

**Proof.** Unraveling the definition of isiso, we have:

$$\mathsf{isiso}_A(f) \equiv \left(\sum(g : \hom_A(y, x), g \circ f = \mathsf{id}_x)\right) \times \left(\sum(h : \hom_A(y, x), f \circ h = \mathsf{id}_y)\right)$$

We can apply $F_\#$ to $g$ and $h$ and $\mathsf{ap}_F$ to the equality. Concatenating with the equality from Theorem 3.3.6 we obtain an element of $\mathsf{isiso}_B(F_\#(f))$. □

## 3.4 Squares

Squares naturally occur when one examines the morphisms of hom-types. Consider for example $\hom_{\hom_A(x,y)}(f, g)$. When unfolding the definitions, we obtain $\left\langle \Delta^1 \to \left\langle \Delta^1 \to A \,\middle|\, {}^{\partial\Delta^1}_{[x,y]} \right\rangle \,\middle|\, {}^{\partial\Delta^1}_{[f,g]} \right\rangle$, which is what we call a *curried* square.

The interesting result of this section is that squares have the same induction principle as identity types. Given a square where two opposing sides of the square are arrows originating from paths, we obtain an induction principle where one can show any property about such squares by only considering the degenerate square $\mathsf{id}_f$ with $f : \hom_A(x, y)$. This result is used later for proving that $\mathrm{Iso}_A(x, y) \to \hom_A(x, y)$ is a full embedding.

There are two possible ways of defining a square, the first being the curried form we naturally obtain through nesting hom-types. The other being the natural definition of a square as a map from $\Delta^1 \times \Delta^1$.

**Definition 3.4.1:** Square

The type of *squares* is $\Delta^1 \times \Delta^1 \to A$. The type of *curried squares* is $\Delta^1 \to \Delta^1 \to A$.

Given $x, y, z, w : A$ and $f_1 : \hom_A(x, y), g_1 : \hom_A(x, z), f_2 : \hom_A(z, w), g_2 : \hom_A(y, w)$, we write:

$$\mathsf{square}_A\begin{pmatrix} x & \xrightarrow{f_1} & y \\ g_1\big\downarrow & & \big\downarrow g_2 \\ z & \xrightarrow[f_2]{} & w \end{pmatrix} :\equiv \left\langle \Delta^1 \times \Delta^1 \to A \,\middle|\, {}^{\Delta^1 \times \partial\Delta^1 \,\vee\, \partial\Delta^1 \times \Delta^1}_{[f_1, f_2] \,\cup\, [g_1, g_2]} \right\rangle.$$

We do the same for curried squares with:

$$\mathsf{csquare}_A\begin{pmatrix} x & \xrightarrow{f_1} & y \\ g_1\big\downarrow & & \big\downarrow g_2 \\ z & \xrightarrow[f_2]{} & w \end{pmatrix} :\equiv \left\langle (t : \Delta^1) \to \left\langle \Delta^1 \to A \,\middle|\, {}^{\partial\Delta^1}_{[f_1(t), f_2(t)]} \right\rangle \,\middle|\, {}^{\partial\Delta^1}_{[g_1, g_2]} \right\rangle$$

These two notions are canonically equivalent via (un)currying.

**Lemma 3.4.2**

The following map is an equivalence:

$$\mathsf{square}_A\begin{pmatrix} x & \xrightarrow{f_1} & y \\ g_1\big\downarrow & & \big\downarrow g_2 \\ z & \xrightarrow[f_2]{} & w \end{pmatrix} \to \mathsf{csquare}_A\begin{pmatrix} x & \xrightarrow{f_1} & y \\ g_1\big\downarrow & & \big\downarrow g_2 \\ z & \xrightarrow[f_2]{} & w \end{pmatrix}$$

$$\sigma \mapsto \lambda t.\ \lambda s.\ \sigma(t, s)$$

**Proof.** We give an inverse: $\lambda\sigma.\ \lambda(t, s).\ \sigma(t)(s)$ (note that we explicitly chose the notation $\sigma(t)(s)$ instead of our usual $\sigma(t, s)$ to avoid confusion). The homotopies of the equivalence are:

$$(\lambda(t,s).\ (\lambda t.\ \lambda s.\ \sigma(t,s))(t)(s)) = \sigma$$
$$(\lambda t.\ \lambda s.\ (\lambda(t,s).\ \sigma(t)(s))(t,s)) = \sigma$$

However, these already hold judgmentally, by the computation rules of $\lambda$. □

Additionally, we have another expected result that a square is just made up of two triangles that agree on the diagonal.

**Lemma 3.4.3**

The following map is an equivalences:

$$\mathsf{square}_A\begin{pmatrix} x & \xrightarrow{f_1} & y \\ g_1\big\downarrow & & \big\downarrow g_2 \\ z & \xrightarrow[f_2]{} & w \end{pmatrix} \to \sum\begin{pmatrix} d : \hom_A(x,w) \\ \hom_A^2\begin{pmatrix} & f_1\nearrow \ y \ \searrow g_2 & \\ x & \xrightarrow[d]{\qquad} & w \end{pmatrix} \\ \hom_A^2\begin{pmatrix} & g_1\nearrow \ z \ \searrow f_2 & \\ x & \xrightarrow[d]{\qquad} & w \end{pmatrix} \end{pmatrix}$$

$$\sigma \mapsto \begin{pmatrix} \lambda t.\ \sigma(t,t) \\ \lambda(t,s).\ \sigma(t,s) \\ \lambda(t,s).\ \sigma(s,t) \end{pmatrix} =: \begin{pmatrix} d \\ \tau_1 \\ \tau_2 \end{pmatrix}$$

**Proof.** We give an inverse:

$$\lambda(d,\tau_1,\tau_2).\ \mathsf{rec}_{\vee}^{\Delta^2,\nabla^2}(\tau_1, \lambda(t,s).\ \tau_2(s,t))$$

Where we write $\nabla^2$ for the lower left triangle in $\Delta^1 \times \Delta^1$ formally defined as $\nabla^2 :\equiv \lambda(t,s){:}\ \mathbb{2}\times\mathbb{2}.\ t \leq s$. The $\mathsf{rec}_{\vee}$ is well-defined, since $\tau_1$ and $\tau_2$ agree on the diagonal. The homotopies are again trivial by applying the computation rules of $\mathsf{rec}_{\vee}$ and $\lambda$. □

Since triangles with one side being the identity are equivalent to paths (in Segal types), we can combine the above equivalence with this fact and obtain that squares with one side being the identity are equivalent to triangles.

**Lemma 3.4.4**

Let $A$ be a Segal type, then there is an equivalence:

$$\mathsf{square}_A\begin{pmatrix} x & \xrightarrow{f_1} & y \\ \mathsf{id}_x\big\downarrow & & \big\downarrow g_2 \\ x & \xrightarrow[f_2]{} & z \end{pmatrix} \simeq \hom_A^2\begin{pmatrix} & f_1\nearrow \ y \ \searrow g_2 & \\ x & \xrightarrow[f_2]{\qquad} & z \end{pmatrix}$$

**Proof.** First we use Lemma 3.4.3, then we apply commutativity of independent $\Sigma$ types followed by associativity of $\Sigma$ types. Then we use Lemma 3.3.4 followed by Lemma 2.3.3:

$$
\begin{aligned}
&\mathsf{square}_A\left(\begin{array}{ccc} x & \xrightarrow{f_1} & y \\ {\scriptstyle \mathsf{id}_x}\downarrow & & \downarrow{\scriptstyle g_2} \\ x & \xrightarrow[f_2]{} & z \end{array}\right) \\
&\simeq \sum\left(d : \hom_A(x,z), \hom^2_A\left(\begin{array}{ccc} & y & \\ {\scriptstyle f_1}\nearrow & & \searrow{\scriptstyle g_2} \\ x & \xrightarrow[d]{} & z \end{array}\right), \hom^2_A\left(\begin{array}{ccc} & x & \\ {\scriptstyle \mathsf{id}_x}\nearrow & & \searrow{\scriptstyle f_2} \\ x & \xrightarrow[d]{} & z \end{array}\right)\right) \\
&\simeq \sum\left(d : \hom_A(x,z), \hom^2_A\left(\begin{array}{ccc} & x & \\ {\scriptstyle \mathsf{id}_x}\nearrow & & \searrow{\scriptstyle f_2} \\ x & \xrightarrow[d]{} & z \end{array}\right), \hom^2_A\left(\begin{array}{ccc} & y & \\ {\scriptstyle f_1}\nearrow & & \searrow{\scriptstyle g_2} \\ x & \xrightarrow[d]{} & z \end{array}\right)\right) \\
&\simeq \sum\left((d,_) : \sum\left(d : \hom_A(x,z), \hom^2_A\left(\begin{array}{ccc} & x & \\ {\scriptstyle \mathsf{id}_x}\nearrow & & \searrow{\scriptstyle f_2} \\ x & \xrightarrow[d]{} & z \end{array}\right)\right), \hom^2_A\left(\begin{array}{ccc} & y & \\ {\scriptstyle f_1}\nearrow & & \searrow{\scriptstyle g_2} \\ x & \xrightarrow[d]{} & z \end{array}\right)\right) \\
&\simeq \sum\left((d,_) : \sum(d : \hom_A(x,z), f_2 = d), \hom^2_A\left(\begin{array}{ccc} & y & \\ {\scriptstyle f_1}\nearrow & & \searrow{\scriptstyle g_2} \\ x & \xrightarrow[d]{} & z \end{array}\right)\right) \\
&\simeq \hom^2_A\left(\begin{array}{ccc} & y & \\ {\scriptstyle f_1}\nearrow & & \searrow{\scriptstyle g_2} \\ x & \xrightarrow[f_2]{} & z \end{array}\right)
\end{aligned}
$$

□

If we have a square with two opposite sides being the identity, we can thus conclude that it is equivalent to a path of the two other sides:

**Corollary 3.4.5**

Let $A$ be a Segal type, then there is an equivalence:

$$
\mathsf{square}_A\left(\begin{array}{ccc} x & \xrightarrow{f} & y \\ {\scriptstyle \mathsf{id}_x}\downarrow & & \downarrow{\scriptstyle \mathsf{id}_y} \\ x & \xrightarrow[g]{} & y \end{array}\right) \simeq (f = g)
$$

**Remark 3.4.6**

We note that Lemma 3.4.2 allows us to use `csquare` instead of `square` in the previous corollary. But such a square exactly is $\hom_{\hom_A(x,y)}(f,g)$. For this reason we use that type in the next theorem.

Now we obtain the first step towards our main result of this subsection, by showing that squares with two constant opposing sides have the same induction scheme as identity types.

**Theorem 3.4.7:** Square Induction

Given a Segal type $A$ and $x, y : A$, $f : \hom_A(x,y)$ as well as a dependent type family:

$$C : \big(g : \hom_A(x,y), \hom_{\hom_A(x,y)}(f,g)\big) \to U$$

If there is a value $c : C(f, \mathsf{id}_f)$, then there is a value of type:

$$\left(g : \hom_A(x, y), \sigma : \hom_{\hom_A(x,y)}(f, g)\right) \to C(g, \sigma)$$

**Proof.** From Corollary 3.4.5 and Lemma 3.4.2 we obtain an equivalence:

$$\gamma : (f = g) \simeq \hom_{\hom_A(x,y)}(f, g)$$

We apply path induction in $\hom_A(x, y)$ with start $f$ and the type family being $\lambda g.\ \lambda p\colon f = g.\ C(g, \gamma(p))$. The proof for $(f, \mathsf{refl})$ is given by $c$ and thus we obtain a function of the following type:

$$(g : \hom_A(x, y), p : f = g) \to C(g, \gamma(p)).$$

Precomposing this function with $\gamma^{-1}$ in the second component and postcomposing with transport of $\gamma(\gamma^{-1}(\sigma)) = \sigma$ (since $\gamma$ is an equivalence) we obtain a function of type:

$$\left(g : \hom_A(x, y), \sigma : \hom_{\hom_A(x,y)}(f, g)\right) \to C(g, \sigma).$$

□

**Corollary 3.4.8**

Theorem 3.4.7 also holds for plain squares by Lemma 3.4.2.

## 3.5 Initial and Final Objects

**Definition 3.5.1:** Initial/Final Object

We call a value $x : A$ *initial* / *final* in $A$ iff the following holds:

$$\begin{aligned} \text{is-initial}_A(x) &:= (a : A) \to \text{is-contr}(\hom_A(x, a)) \\ \text{is-final}_A(x) &:= (a : A) \to \text{is-contr}(\hom_A(a, x)) \end{aligned}$$

**Lemma 3.5.2:** (cf. [2, Lemma 9.8])

Let $A$ be a Segal type and $a : A$. Then we have that $(a, \mathsf{id}_a)$ is initial in $\Sigma(x : A, \hom_A(a, x))$.

## 3.6 Full Embeddings

In Section 2.6, we introduced the notion of embedding – a function that is an embedding equates the groupoidal structures of values taken on by the function. In HoTT, this is sufficient, as there is no further structure that a type can exhibit. In simplicial HoTT however, we have additional structure in the form of hom types. The notion of a full embedding remedies this missing terminology. Instead of asking that the groupoidal structures match, we ensure that the hom types are equivalent. This then provides all higher coherences as well. Full embeddings thus preserve concepts from category theory such as initial and final objects. We use full embeddings exactly for this purpose later, as they *detect* initial objects, which is very useful to determining that a certain map is an *initial section.*

A more common name for full embeddings is fully faithful functors. Since we inherited the notion of embeddings from HoTT, we use the name full embeddings, as they properly generalize that notion as well.

**Definition 3.6.1:** Full Embedding

A map $f : A \to B$ is a full embedding, iff for every $x, y : A$ we have that

$$f_\# : \hom_A(x, y) \to \hom_B(f(x), f(y))$$

is an equivalence.

To justify the name "full embedding", we now show that any full embedding is an embedding, so it also preserves the groupoidal structure of types.

**Lemma 3.6.2**

Let $f : A \to B$ be a full embedding between Rezk types. Then $f$ is an embedding.

**Proof.** We need to show for any $x, y : A$ that the following map is an equivalence:

$$\mathsf{ap}_f : (x = y) \to (f(x) = f(y))$$

Consider the following diagram:

$$\begin{array}{ccc}
(x = y) & \xrightarrow{\mathsf{ap}_f} & (f(x) = f(y)) \\
\simeq \downarrow & & \downarrow \simeq \\
\mathrm{Iso}_A(x, y) & \xrightarrow{(f_\#, \beta)} & \mathrm{Iso}_B(f(x), f(y)) \\
\downarrow & & \downarrow \\
\hom_A(x, y) & \xrightarrow[f_\#]{\simeq} & \hom_B(f(x), f(y))
\end{array}$$

Here $\beta : \mathsf{isiso}_A(g) \to \mathsf{isiso}_B(f_\#(g))$ is the map defined by Lemma 3.3.14. The outer square commutes by [2, Lemma 10.8], the lower square commutes by inspection (the vertical maps are just forgetting the isomorphism property) and the upper square commutes because of the other two squares commuting and the fact that $\mathsf{isiso}(_)$ a proposition.

Now all that's left to show is that $(f_\#, \beta)$ is an equivalence. Since $f_\#$ is an equivalence, one way to proceed is to show that $\beta$ also is an equivalence. Since $\mathsf{isiso}(_)$ is a proposition, we only need to produce a map in the backwards direction. We recall the definition of $\mathsf{isiso}_A(g)$ and $\mathsf{isiso}_B(f_\#(g))$:

$$\sum \begin{pmatrix} h_1 : \hom_A(y, x) \\ g \circ h_1 = \mathsf{id}_x \end{pmatrix} \times \sum \begin{pmatrix} h_2 : \hom_A(y, x) \\ h_2 \circ g = \mathsf{id}_y \end{pmatrix}$$

$$\sum \begin{pmatrix} h_1 : \hom_B(f(y), f(x)) \\ f_\#(g) \circ h_1 = \mathsf{id}_{f(x)} \end{pmatrix} \times \sum \begin{pmatrix} h_2 : \hom_B(f(y), f(x)) \\ h_2 \circ f_\#(g) = \mathsf{id}_{f(y)} \end{pmatrix}$$

The function $\beta$ from Lemma 3.3.14 applies $f_\#$ in the first components of each sigma type. Since it is an equivalence, to map back, we can use its inverse. The second components are mapped by $\beta$ using $\mathsf{ap}_{f_\#}$ and concatenation with functoriality (Theorem 3.3.6). Since concatenation is an equivalence and $\mathsf{ap}$ of an equivalence also is an equivalence, we can use their composed inverses to map the second components of the sigma types.

Finally we have a map $\mathsf{isiso}_B\big(f_\#(g)\big) \to \mathsf{isiso}_A(g)$ and thus conclude that $\beta$ is an equivalence. □

Any equivalence naturally is a full embedding.

**Lemma 3.6.3:** (cf. [3, `is-equiv-ap-hom-is-equiv`])

Let $f : A \to B$ be an equivalence, then $f$ is a full embedding.

As one would expect, full embeddings compose.

**Lemma 3.6.4**

Let $f : A \to B$ and $g : B \to C$ be full embeddings, then their composition also is a full embedding.

**Proof.** Theorem 3.3.6 proves that $(-)_\#$ is functorial, so we have $(g \circ f)_\# = g_\# \circ f_\#$ and since equivalences compose, $g \circ f$ is a full embedding. □

Now to our main use-case of a full embedding. They detect initial objects:

**Lemma 3.6.5**

Let $f : A \to B$ be a full embedding. If $f(a)$ is initial in $B$, then $a$ is initial in $A$.

**Proof.** For $a$ to be initial in $A$, we have to show that $\hom_A(a, x)$ is contractible for any $x : A$. Since $f$ is a full embedding, we have $\hom_A(a, x) \simeq \hom_B(f(a), f(x))$, which is contractible since $f(a)$ is initial in $B$. □

As the last statement of this section we prove that a specific map is a full embedding. We use this theorem later when we prove that LARI adjunctions are equivalent to initial sections. We prove the underlying mechanism of the theorem using the next lemma: a square in a Rezk type with two opposing sides being isomorphisms already consists of isomorphisms (when slicing vertically).

**Lemma 3.6.6**

Let $A$ be a Rezk type and we are given a square in $A$:

$$F : \mathsf{csquare}_A \begin{pmatrix} x_1 & \xrightarrow{f} & y_1 \\ h\big\downarrow & & \big\downarrow k \\ x_2 & \xrightarrow[g]{} & y_2 \end{pmatrix}$$

where $h$ and $k$ are isomorphisms. Then the following type is contractible:

$$(t : \Delta^1) \to \mathsf{isiso}_A(F(t))$$

**Proof.** We apply Theorem 3.3.13 to $h$ and $k$ to obtain the following situation:

$$F : \mathsf{csquare}_A \begin{pmatrix} x & \xrightarrow{f} & y \\ \mathsf{id}_x\big\downarrow & & \big\downarrow \mathsf{id}_y \\ x & \xrightarrow[g]{} & y \end{pmatrix}$$

We still have to show that $(t : \Delta^1) \to \mathsf{isiso}_A(F(t))$ is contractible. Now we can apply Theorem 3.4.7 by showing that $(t : \Delta^1) \to \mathsf{isiso}_A(\mathsf{id}_f(t))$ is contractible, but that is trivial, since $\mathsf{id}_f(t) \equiv \mathsf{id}_{f(t)}$ which always is an isomorphisms. □

Finally, we can derive our theorem to show that the inclusion of $\mathrm{Iso}_A \to \mathrm{hom}_A$ is a full embedding. We give a much more general result than only considering this basic inclusion. This is useful later to show that this full embedding also translates to the situation when we sum over $x$ and $y$.

**Theorem 3.6.7**

Let $A$ be a Rezk type, $x, y : \Delta^1 \to A$ and $f : \mathrm{Iso}_A(x(0), y(0)), g : \mathrm{Iso}_A(x(1), y(1))$. Then the following map is an equivalence:

$$\left\langle (t : \Delta^1) \to \mathrm{Iso}_A(x(t), y(t)) \;\middle|\; {}^{\partial\Delta^1}_{[f,g]} \right\rangle \to \left\langle (t : \Delta^1) \to \mathrm{hom}_A(x(t), y(t)) \;\middle|\; {}^{\partial\Delta^1}_{[\pi_1 f, \pi_1 g]} \right\rangle$$
$$\sigma \mapsto \lambda t.\ \lambda s.\ \sigma(t, s)$$

**Proof.** We note that the map of the theorem arises as the composition of the equivalence in Proposition 3.2.7 and projecting away the second component of the sigma type, which also is an equivalence, since the second component is contractible by Lemma 3.6.6. □

**Corollary 3.6.8**

Let $A$ be a Rezk type and $x, y : A$, then the following map is a full embedding:

$$\mathrm{Iso}_A(x, y) \to \mathrm{hom}_A(x, y)$$
$$f \mapsto \pi_1(f)$$

## 3.7 Adjunctions

There are several ways of defining adjunctions in our theory. Riehl and Shulman show through [2, Theorem 11.23] that they are equivalent when considering Rezk types. Since we are mostly interested in synthetic categories, we only concern ourselves with transposing adjunctions, as these generalize nicely to LARI adjunctions later.

**Definition 3.7.1:** Transposing Adjunctions

Given two types $A, B : \mathcal{U}$, we call the following data a *transposing adjunction between $A$ and $B$*:

$$f : A \to B \qquad\qquad g : B \to A$$
$$\varphi : (a : A, b : B) \to (\mathrm{hom}_B(f(a), b) \simeq \mathrm{hom}_A(a, g(b)))$$

We still need to be able to talk about the unit and counit of an adjunction which we define in terms of the equivalence $\varphi$.

**Definition 3.7.2:** Unit and Counit of an Adjunction

Given a transposing adjunction $f : A \to B, g : B \to A$ and

$$\varphi : (a : A, b : B) \to \mathrm{hom}_B(f(a), b) \simeq \mathrm{hom}_A(a, g(b))$$

We define the unit and counit as:

$$\eta_\varphi : (a : A) \to \hom_A(a, g(f(a))) \qquad \varepsilon_\varphi : (b : B) \to \hom_B(f(g(b)), b)$$
$$a \mapsto \varphi(a, f(a))\left(\mathsf{id}_{f(a)}\right) \qquad b \mapsto \varphi(g(b), b)^{-1}\left(\mathsf{id}_{g(b)}\right)$$

## 3.8 Cat and Category-theoretic Grothendieck Construction

Up until now, we studied how category theory arises from introducing directionality in our type theory. This approach leaves us with some types whose structure is non-categorical. In addition, the universe is that of all *types* and not of categories. We need to explicitly construct the type of categories (again with levels) and then hope that this type is itself a category (in a higher level). We ultimately fail to produce a category of categories and give an argument that in this theory even if we could, the categorical Grothendieck construction is impossible using any definition of type of categories. When comparing this situation with HoTT and $\infty$-groupoids, we can really appreciate the beauty of that theory: the universe in HoTT only contains groupoids (all types are groupoids after all) and for the same reason, the universe itself also is a groupoid.

### 3.8.1 Cat is not Rezk

We start with the sensible candidate of taking all Rezk types from our universe:

$$\mathrm{Cat} :\equiv \Sigma(C : \mathcal{U}, \text{is-rezk}(C))$$

To prove that Cat is Rezk, we need to first show that it is Segal. Consider the commutative diagram:

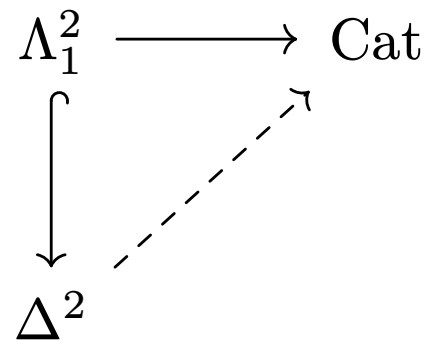

Where we need to produce a unique dashed map. In HoTT, this difficulty is circumvented by equipping each type with the groupoidal structure from the get-go. One does not need to know the structure of $x =_A y$, as it already satisfies all the required axioms of $\infty$-groupoids.

**Argument 3.8.1:** Cat is not Segal

For this argument, we assume the existence of $\Lambda^2_1$ and $\Delta^1$ as types. Consider the following maps:

$$\mathrm{id}_{\Lambda^2_1} : \Lambda^2_1 \to \Lambda^2_1 \quad t \mapsto t$$
$$\mathrm{id}_{\Delta^2} : \Delta^2 \to \Delta^2 \quad t \mapsto t$$
$$i_{\Lambda^2_1} : \Lambda^2_1 \to \Delta^2 \quad t \mapsto t$$

Their family equivalents under the Grothendieck construction are:

$$A : \Lambda^2_1 \to \mathcal{U} \quad t \mapsto \mathbb{1}$$
$$B : \Delta^2 \to \mathcal{U} \quad t \mapsto \mathbb{1}$$
$$C : \Delta^2 \to \mathcal{U} \quad t \mapsto \texttt{if } t \in \Lambda^2_1 \texttt{ then } \mathbb{1} \texttt{ else } \mathbb{0}$$

The codomains of $A, B, C$ all are in Cat, since propositions are Rezk. Additionally, $A, B$ and $C$ agree on $\Lambda^2_1$ and thus $B$ and $C$ both extend $A$. However, $B$ and $C$ are clearly not equivalent. We thus have that Cat is not even Segal, since the uniqueness condition of the filler is violated.

This defect remains even if we replace the condition of is-rezk with something else that still allows propositions.

We note that it is possible to define a Rezk type of groupoids that satisfies directed univalence. D. Gratzer, J. Weinberger, and U. Buchholtz [24] show under significant modifications of the theory that there is a Rezk type $\mathcal{S}$ whose elements are Rezk types that are groupoids i.e. $\hom_A(x, y) \simeq (x =_A y)$. In addition it satisfies directed univalence:

$$\hom_{\mathcal{S}}(A, B) \simeq (A \to B)$$

### 3.8.2 Grothendieck

In addition to being an $\infty$-groupoid itself, the universe in HoTT exhibits the Grothendieck construction. Any map into the universe can be equivalently represented by a fibration over the domain of the map as we explored in Section 2.8. This fibration does not require any type-theoretical restrictions, as any map automatically is a functor with respect to the groupoidal structure.

Even if our Cat would be Rezk (or at least Segal), we run into the problem of a faithful Grothendieck construction. We would like to have the same as Theorem 2.8.4 but with $\mathcal{U}$ replaced by Cat, so:

**Conjecture 3.8.2**

For any Rezk type $C$ : Cat, the (un)straightening maps form an equivalence:

$$\begin{aligned} \overline{\mathrm{Un}}_C : \overline{\mathrm{Fam}}(C) &\to \overline{\mathrm{Fib}}(C) & \overline{\mathrm{St}}_C : \overline{\mathrm{Fib}}(C) &\to \overline{\mathrm{Fam}}(C) \\ P &\mapsto (\Sigma P, \pi_P) & (A, f) &\mapsto \mathrm{fib}_f \end{aligned}$$

With $\overline{\mathrm{Fam}}(C) :\equiv C \to \mathrm{Cat}$ and $\overline{\mathrm{Fib}}(C) :\equiv \Sigma(A : \mathrm{Cat}, A \to C)$. We use a line over the letters to distinguish them from the family and fibration notion used in the type-theoretic Grothendieck construction.

This statement however is not the equivalent of the categorical straightening theorem which reads:

**Theorem 3.8.3:** Categorical Straightening (cf. [21, Theorem 3.3.10])

For every $\infty$-category $\mathcal{C}$, there is an equivalence of $\infty$-categories:

$$\mathrm{CoCart}(\mathcal{C}) \simeq \mathrm{Fun}(\mathcal{C}, \mathrm{Cat}_\infty)$$

Where $\mathrm{CoCart}(\mathcal{C})$ is the subcategory of $(\mathrm{Cat}_\infty)_{/\mathcal{C}}$ with objects restricted to cocartesian fibrations and morphisms to the morphisms of cocartesian fibrations (we have and will not define them here).

This theorem is considerably different from our conjecture. While we can associate $\mathrm{Fun}(\mathcal{C}, \mathrm{Cat}_\infty)$ with $\overline{\mathrm{Fam}}(C)$, the category of cocartesian fibrations over $\mathcal{C}$ is not accurately represented by $\overline{\mathrm{Fib}}(C)$. As a more accurate representation, we could opt for defining:

$$\overline{\mathrm{Fib}}(C) :\equiv \Sigma(p : \mathrm{Fib}(C), \text{is-cocartesian-fib}(p))$$

Using this notion of fibration in our Conjecture 3.8.2 now gives the right version of categorical straightening. However, this version cannot be proven in our theory due to the following fundamental limitation.

**Argument 3.8.4**

We begin our argument by analyzing $\overline{\mathrm{Fam}}(C)$. Using Proposition 2.5.7 we obtain:

$$\overline{\mathrm{Fam}}(C) \simeq \Sigma(B : \mathrm{Fam}(C), (c : C) \to \text{is-rezk}(B(c)))$$

Applying the type theoretic Straightening Theorem 2.8.4, we obtain the equivalent[1] type:

$$\Sigma\big(p : \mathrm{Fib}(C), (c : C) \to \text{is-rezk}\big(\mathrm{fib}_p(c)\big)\big)$$

This almost is $\overline{\mathrm{Fib}}(B)$, with a stark difference in the second component: the condition we want is "is-cocartesian-fib", which **cannot** be checked fiberwise. So even if we were to swap our definition of Cat for something different, in general $\mathrm{Cat} \equiv \Sigma(C : \mathcal{U}, P(C))$ for some $P : \mathcal{U} \to \mathbb{P}$. Then we would arrive at the same issue here that:

$$\Big((c : C) \to P\big(\mathrm{fib}_p(c)\big)\Big) \not\simeq \text{is-cocartesian-fib}(p)$$

Despite this flaw, the theory is quite powerful in describing statements about $(\infty, 1)$-categories that do not involve $(\infty, 2)$-categories. In particular we proceed with defining cocartesian fibrations. However, if we take our usual family-viewpoint, we want to have *cocartesian families* instead.

Cocartesian families are those families that allow us to define a *directed transport*. We again note that in HoTT, the transport principle is available for all type families. The later chapters of this thesis formalize this notion of cocartesian family along with its closure properties.

[1] We did not assume univalence in our theory, so this formally isn't an equivalence. But we're only interested in showing the meta-theoretic argument that the Cat we defined doesn't have the desired properties, even if we were to assume univalence.

# 4 Orthogonal Families and Fibrations

In this section, we study orthogonal families and observe that they are equivalent to orthogonal fibrations. They have been studied extensively in the past and are a very natural notion from topology. There are lots of different ways to define what an orthogonal fibration is and we study three of them in addition to the notion of orthogonal family. `sHoTT` already defines orthogonal fibrations, but in a different way compared to Buchholtz and Weinberger [5]. Thus, the first part of this section is proving their equivalence. Later in this section, we study *inner families* – a special case of orthogonal families – these are the reason for our study of orthogonal fibrations, as we need them later for initial sections.

Inner families essentially are "dependent Segal types", so given a type family $B : A \to \mathcal{U}$ that is inner, we can compose dependent morphisms. For this, one needs a triangle in the base $A$ and then a horn that lives dependently above that triangle, see Figure 3. The composition is then a *filler* of that horn (the blue striped area and the dotted line in the upper triangle), which in inner families is going to be contractible— hence we have a unique composition. Note that we don't require the base $A$ to be Segal.

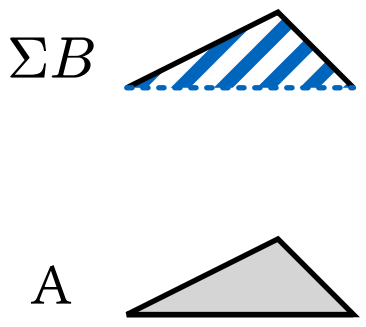

Figure 3: Inner family

**Definition 4.0.1:** Right orthogonal family

We call a type family $B : A \to \mathcal{U}$ a *right orthogonal family* with respect to a shape inclusion $Y \subset X$ if:

$$\text{is-orthogonal-family}_Y(B) :\equiv \begin{pmatrix} a : X \to A \\ f : (y : Y) \to B(a(y)) \end{pmatrix} \to \text{is-contr}\left(\left\langle (x : X) \to B(a(x)) \;\middle|\; {}^{Y}_{f} \right\rangle\right)$$

holds.

This definition essentially says that for every diagram of shape $X$ in the base and a diagram of shape $Y$ living above it, there is a unique extension to a diagram of shape $X$ living above the base. The fibration version is expectedly more complicated:

**Definition 4.0.2:** Right orthogonal fibration (cf. [5, Section 3])

A map $f : E \to B$ is a *right orthogonal fibration* with respect to a shape inclusion $Y \subset X$ if for any two maps $p : Y \to E$ and $q : X \to B$ that satisfy the commutative diagram:

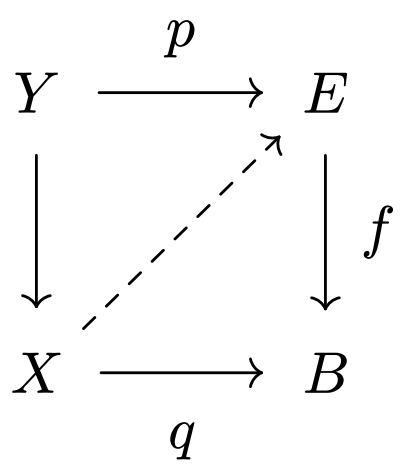

there exists a unique lift $h : X \to E$. Type-theoretically, we can express this as the proposition:

$$\begin{aligned} \text{is-right-orthogonal}_Y(f) :\equiv \left(p : Y \to E, q : \left\langle X \to B \;\middle|\; {}^{Y}_{f \circ p} \right\rangle\right) \\ \to \text{is-contr}\left(\Sigma\left(p' : \left\langle X \to E \;\middle|\; {}^{Y}_{p} \right\rangle, \left\langle (x : X) \to f(p'(x)) = q(x) \;\middle|\; {}^{Y}_{\lambda_.\ \mathsf{refl}} \right\rangle\right)\right) \end{aligned}$$

Given two maps $p, q$ that commute with $f$ and the shape inclusion, we have that the map along with the proof that it commutes is unique. Here we encoded the commutativity directly in the extension types (we could also have used propositional equality, but Proposition 3.2.9 tells us that that's equivalent. However, that would make this definition even longer).

Naturally, the family notion and fibration notion agree when we consider a type family and its canonical fibration.

**Lemma 4.0.3**

A type family $B : A \to \mathcal{U}$ is right orthogonal iff $\pi_B$ is right orthogonal.

**Proof.** We show the iff in one step by giving an equivalence of the two types that ought to be contractible for the two notions. Plugging in a family for the fibration definition yields:

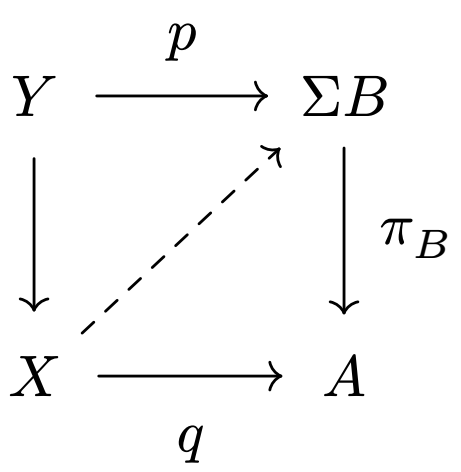

$$\Sigma\left(p' : \left\langle X \to \Sigma B \;\middle|\; {}^{Y}_{p} \right\rangle, \left\langle (x : X) \to \pi_B(p'(x)) = q(x) \;\middle|\; {}^{Y}_{\lambda_.\ \mathsf{refl}} \right\rangle\right)$$

Now consider the following equivalence:

$$\begin{aligned} & \Sigma\left(p' : \left\langle X \to \Sigma B \;\middle|\; {}^{Y}_{p} \right\rangle, \left\langle (x : X) \to \pi_B(p'(x)) = q(x) \;\middle|\; {}^{Y}_{\lambda_.\ \mathsf{refl}} \right\rangle\right) \\ \simeq & \left\langle (x : X) \to \Sigma(e : \Sigma B, \pi_B(e) = q(x)) \;\middle|\; {}^{Y}_{\lambda y.\ (p(x),\ \mathsf{refl})} \right\rangle \\ \equiv & \left\langle (x : X) \to \mathrm{fib}_{\pi_B}(q(x)) \;\middle|\; {}^{Y}_{\lambda y.\ (p(x),\ \mathsf{refl})} \right\rangle \\ \simeq & \left\langle (x : X) \to B(q(x)) \;\middle|\; {}^{Y}_{\pi_2 \circ p} \right\rangle \end{aligned}$$

The first equivalence comes from Proposition 2.5.7 and the second one from Lemma 2.8.2.

The final type of the equivalence is exactly the type that is asked to be contractible in the definition of right orthogonal family. ☐

sHoTT uses a different definition for right orthogonal fibrations:

**Definition 4.0.4:** Right orthogonal fibration (2)

A map $f : E \to B$ is a *right orthogonal fibrations* with respect to a shape inclusion $Y \subset X$ if the following square is homotopy-cartesian:

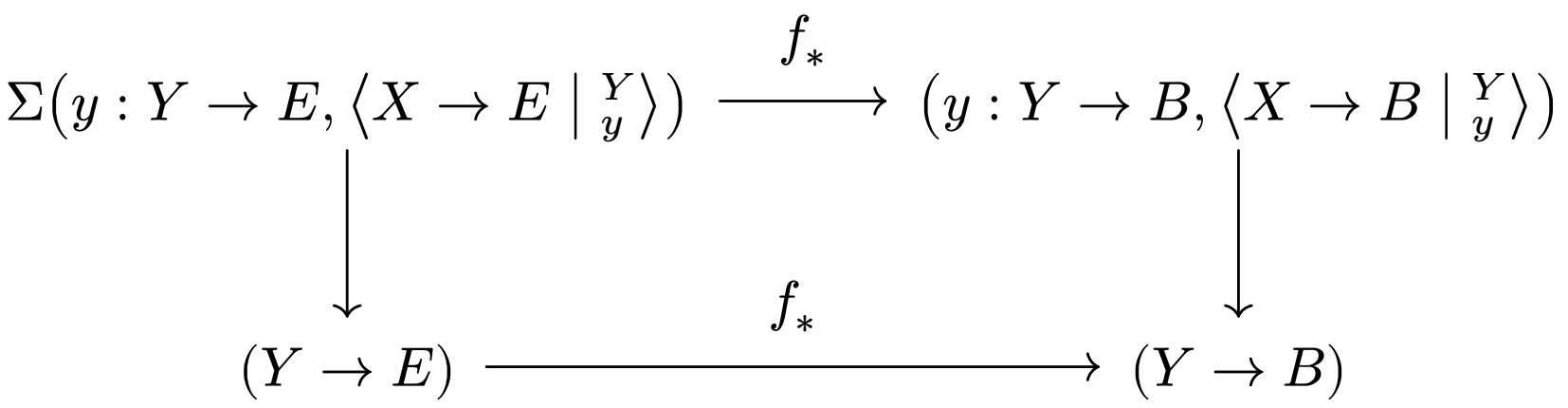

Being homotopy-cartesian means the following: for any $y : Y \to E$ we have that the induced map on the fiber is an equivalence:

$$f_* : \left\langle X \to E \;\middle|\; {}^{Y}_{y} \right\rangle \to \left\langle X \to B \;\middle|\; {}^{Y}_{f^*(y)} \right\rangle$$

This notion is equivalent to our other two definitions via the following lemma.

**Lemma 4.0.5:** [3, `has-contr-relative-extension-types-iff-is-right-orthogonal`]

A map $f : E \to B$ a is right orthogonal fibration (the `sHoTT` definition 4.0.4) iff the following holds:

$$\left(\begin{array}{c}\alpha : Y \to E\\ \beta : \left\langle X \to B \mid {}^{Y}_{f\circ\alpha}\right\rangle\end{array}\right) \to \text{is-contr}\left(\sum\left(\begin{array}{c}\alpha' : \langle X \to E \mid {}^{Y}_{\alpha}\rangle \\ \langle (x:X) \to f(\alpha'(x)) = \beta(x) \mid {}^{Y}_{\lambda_.\ \mathsf{refl}}\rangle\end{array}\right)\right)$$

This is exactly the definition of right orthogonal fibration by Buchholtz and Weinberger [5] also given above.

**Remark 4.0.6**

The whole theory about right orthogonal fibrations/families is dualizable by defining the notion of left fibration through the same commutative square – only by asking it to be the left map. The only issue for this is that we asked the left map to be a shape inclusion and for the right to not be one. So we would have to flip those two as well. This is yet another place where the distinction of shapes and types is biting us.

Since we do not make use of the dual – *left* orthogonal families – we call right orthogonal families/fibrations just *orthogonal families/fibrations*.

## 4.1 Leibniz Cotensor Map

Another alternative formulation of orthogonal fibrations is via the *leibniz cotensor map*. It is used extensively by Buchholtz and Weinberger [5], as it is the definition that most naturally carries over to LARI families.

**Definition 4.1.1:** Leibniz cotensor map (cf. [5, Definition 2.4.3])

Given a shape inclusion $Y \subset X$ and a map $f : E \to B$, we call the universal map from $E^X$ in to the pullback given below the *leibniz cotensor map:*

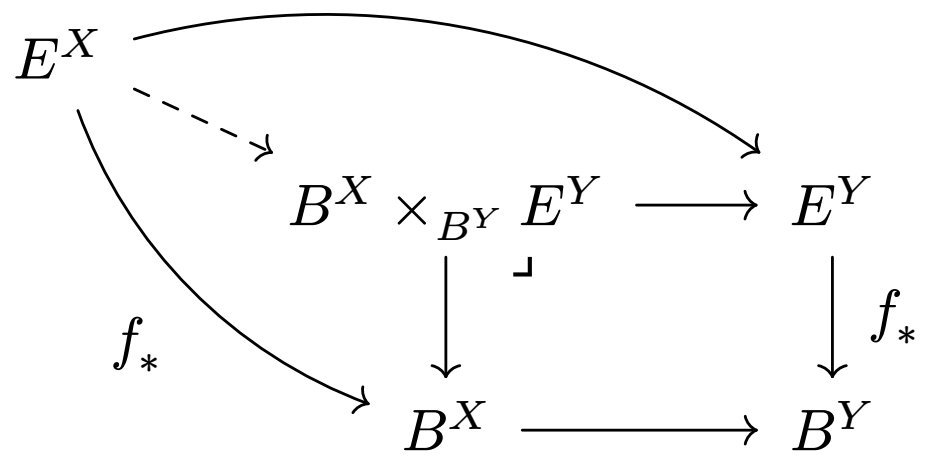

We denote the map $E^X \to B^X \times_{B^Y} E^Y$ with $\widehat{\pitchfork}_Y f$.

Since Buchholtz and Weinberger allow interpreting shapes as types, they have to record a map $j : Y \to X$ in the leibniz cotensor map and write $j \widehat{\pitchfork} f$ instead.

If the map $f$ is an orthogonal fibration (orthogonal to $Y \subset X$), then the leibniz cotensor map is an equivalence.

**Lemma 4.1.2:** cf. [5, Proposition 3.1.1]

A map $f : E \to B$ is an orthogonal fibration wrt. $Y \subset X$ iff the leibniz cotensor map $\widehat{\pitchfork}_Y f$ is an equivalence.

## 4.2 Closure Properties

Buchholtz and Weinberger [5] show several closure properties for orthogonal fibrations. These are also proven in [3] although with the different definition. We only cite them here to compare later with the closure properties of LARI adjunctions.

**Theorem 4.2.1:** Closure properties orthogonal fibrations

Let $Y \subset X$ be a shape inclusion. Then:

*Equivalences* any equivalence $f : B' \simeq B$ is an orthogonal fibration ([5, Proposition 3.1.3]),

*Products* given a family of orthogonal families: $B : I \to \mathcal{U}$ and $P : (i : I) \to B(i) \to \mathcal{U}$ with $P(i)$ being an orthogonal family, then $(i : I) \to P(i)$ is an orthogonal family ([5, Proposition 3.1.4]),

*Composition & left cancellation* Given two families $P : B \to \mathcal{U}$ and $Q : \Sigma P \to \mathcal{U}$. If $P$ is orthogonal, then the composite $Q \odot P : B \to \mathcal{U}$ is orthogonal iff $Q$ is.

*Pullbacks* The pullback of an orthogonal fibration $f : E \to B$ along any map $k : A \to B$ is an orthogonal fibration. ([5, Theorem 3.1.8])

## 4.3 Inner Families

**Definition 4.3.1:** Inner Family

A type family $B : A \to \mathcal{U}$ is called *inner* if it is orthogonal to the shape inclusion of the inner horn into the 2-simplex: $\Lambda^2_1 \subset \Delta^2$.

**Lemma 4.3.2**

Let $B : A \to \mathcal{U}$ be an inner family, then $B(a)$ is Segal for any $a : A$.

**Proof.** Fix $a : A$ and $x, y, z : B(a)$ and $f : \hom_{B(a)}(x, y), g : \hom_{B(a)}(y, z)$. We now have to show that the following type is contractible:

$$\hom^2_{B(a)}\left(\begin{array}{ccc} & y & \\ {\scriptstyle f}\nearrow & & \searrow{\scriptstyle g} \\ x & \longrightarrow & z \end{array}\right)$$

Since $B$ is an inner family, we can apply the inner condition to the following data: the lower triangle living in $A$ is the constant map $\lambda t \colon \Delta^2.\ a$ and the horn above is $[f, g]$. Thus, we obtain that the following type is contractible:

$$\left\langle \Delta^2 \to B(a) \,\middle|\, {}^{\Lambda^2_1}_{[f,g]} \right\rangle$$

Through inspecting this type it is clear that it is equivalent to the $\hom^2$ type we gave above. □

Exactly analogously to Segal types, we can give an equivalence between paths of dependent morphisms and homotopies of dependent morphisms. For this purpose, we first need to define "homotopies between dependent morphisms":

**Definition 4.3.3**

Given $B : A \to \mathcal{U}$ and $\tau : \Delta^2 \to A$ with lifts $x : B(\tau(0,0)), y : B(\tau(1,0)), z : B(\tau(1,1))$ and $f : \mathrm{dhom}_B^{\lambda t.\ \tau(t,0)}(x, y)$, $g : \mathrm{dhom}_B^{\lambda t.\ \tau(1,t)}(y, z)$, $h : \mathrm{dhom}_B^{\lambda t.\ \tau(t,t)}(x, z)$. We define:

$$\mathrm{dhom}_B^{2,\tau}\left(\begin{array}{ccc} & y & \\ {}^{f}\nearrow & & \searrow^{g} \\ x & \xrightarrow[h]{} & z \end{array}\right) := \left\langle (t : \Delta^2) \to B(\tau(t)) \;\middle|\; {}^{\partial\Delta^2}_{[f,g,h]} \right\rangle$$

**Lemma 4.3.4**

Let $B : A \to \mathcal{U}$ be an inner family and $f : \hom_A(x, y)$ as well as $g, h : \mathrm{dhom}_B^f(x', y')$. Then the following natural map is an equivalence:

$$(g = h) \to \mathrm{dhom}_B^{2,\lambda t.\ f(\pi_1 t)}\left(\begin{array}{ccc} & x' & \\ {}^{\mathsf{id}_{x'}}\nearrow & & \searrow^{g} \\ x' & \xrightarrow[h]{} & y' \end{array}\right)$$

**Proof.** The proof is analogous to the same statement in Segal types Lemma 3.3.4 (proven in [2, Proposition 5.10]). We have that the following two types are contractible and thus equivalent:

$$\sum\begin{pmatrix} h' : \mathrm{dhom}_B^f(x', y') \\ g = h' \end{pmatrix} \simeq \sum\begin{pmatrix} h' : \mathrm{dhom}_B^f(x', y') \\ \mathrm{dhom}_B^{2,\lambda t.\ f(\pi_1 t)}\left(\begin{array}{ccc} & x' & \\ {}^{\mathsf{id}_{x'}}\nearrow & & \searrow^{g} \\ x' & \xrightarrow[h']{} & y' \end{array}\right) \end{pmatrix}$$

The former is contractible by Lemma 2.3.3, the latter is contractible, because $B$ is inner. For this reason the map that is identity in the first component and the natural path induction in the second one is an equivalence. Since the first component map is the identity, it follows that the second map must also be an equivalence fiberwise. □

## 4.4 Dependent Squares

We now give the dependent versions of some of the statements from Section 3.4. The proofs and ideas translate very naturally, in fact we omit the proofs, as they are exactly the same modulo dependency. Luckily we do not require the induction principle for dependent squares, which we thus omit. We first need to define dependent squares:

**Definition 4.4.1:** Dependent Square

Given $B : A \to \mathcal{U}$, $\sigma : \Delta^1 \times \Delta^1 \to A$ and lifts on the boundary:

$$x : B(\sigma(0,0)), y : B(\sigma(1,0)), z : B(\sigma(0,1)), w : B(\sigma(1,1))$$
$$f_1 : \mathrm{dhom}_B^{\lambda t.\ \sigma(t,0)}(x, y), g_1 : \mathrm{dhom}_B^{\lambda t.\ \sigma(0,t)}(x, z)$$
$$f_2 : \mathrm{dhom}_B^{\lambda t.\ \sigma(t,1)}(z, w), g_2 : \mathrm{dhom}_B^{\lambda t.\ \sigma(1,t)}(y, w)$$

we define:

$$\mathsf{dsquare}_B^\sigma\left(\begin{array}{ccc} x & \xrightarrow{f_1} & y \\ {}^{g_1}\downarrow & & \downarrow^{g_2} \\ z & \xrightarrow[f_2]{} & w \end{array}\right) := \left\langle (t : \Delta^1 \times \Delta^1) \to B(\sigma(t)) \;\middle|\; {}^{\partial\Delta^1\times\Delta^1 \,\cup\, \Delta^1\times\partial\Delta^1}_{[f_1,f_2] \,\cup\, [g_1,g_2]} \right\rangle$$

**Lemma 4.4.2**

The following map is an equivalence:

$$\mathsf{dsquare}^{\sigma}_{B}\begin{pmatrix} x' & \xrightarrow{f_1'} & y' \\ g_1'\downarrow & & \downarrow g_2' \\ z' & \xrightarrow[f_2']{} & w' \end{pmatrix} \to \sum \begin{pmatrix} d : \mathrm{dhom}_B^{\lambda t.\ \sigma(t,t)}(x', w') \\ \mathrm{dhom}_B^{2,\sigma(t)}\begin{pmatrix} & f_1' \nearrow y' \searrow g_2' & \\ x' & \xrightarrow[d]{\qquad} & w' \end{pmatrix} \\ \mathrm{dhom}_B^{2,\lambda(t,s).\ \sigma(s,t)}\begin{pmatrix} & g_1' \nearrow z' \searrow f_2' & \\ x' & \xrightarrow[d]{\qquad} & w' \end{pmatrix} \end{pmatrix}$$

$$\sigma' \mapsto \begin{pmatrix} \lambda t.\ \sigma'(t,t) \\ \lambda(t,s).\ \sigma'(t,s) \\ \lambda(t,s).\ \sigma'(s,t) \end{pmatrix}$$

**Lemma 4.4.3**

Let $B : A \to \mathcal{U}$ be an inner family, then there is an equivalence:

$$\mathsf{dsquare}^{\sigma}_{B}\begin{pmatrix} x & \xrightarrow{f_1} & y \\ \mathsf{id}_x\downarrow & & \downarrow g_2 \\ x & \xrightarrow[f_2]{} & z \end{pmatrix} \simeq \mathrm{dhom}_B^{2,\tau}\begin{pmatrix} & f_1 \nearrow y \searrow g_2 & \\ x & \xrightarrow[f_2]{\qquad} & z \end{pmatrix}$$

# 5 Initial Sections

In this section we study *initial sections*. They are equivalent to LARI adjunctions which we show in the next section. We use them as they are significantly simpler than LARI adjunctions and allow us to more easily prove the closure properties of LARI families. To show the equivalence of initial sections and LARI adjunctions, we introduce *dependent initial section*. It is an a priori stronger notion that we show to be equivalent to initial section under the assumption that the family is inner.

**Remark 5.0.1**

The name "section" comes from category theory: a section of a map $f : A \to B$ is a map going back $s : B \to A$ with $f \circ s = \mathrm{id}_B$. To translate this situation into type theory, we just ask for extra data for the equality in the form of a homotopy $p : (b : B) \to f(s(b)) = b$.

We can then translate this into a nicer situation by using the Grothendieck construction: consider the family $\mathrm{fib}_f : B \to \mathcal{U}$ derived from $f$. We can construct a dependent map of this family using $s$:

$$\begin{aligned} s' : (b : B) \to & \overbrace{\mathrm{fib}_f(b)}^{\equiv \Sigma(a:A, f(a)=b)} \\ b \mapsto & (s(b), p(b)) \end{aligned}$$

We call such a map a section of the family $\mathrm{fib}_f$. We now have produced a situation where we have a *strict* section:

$$\begin{aligned} \pi_{\mathrm{fib}_f} : \Sigma\mathrm{fib}_f \to B \qquad\qquad & \overline{s'} : B \to \Sigma\mathrm{fib}_f \\ (b, (a, q)) \mapsto b \qquad\qquad & b \mapsto (b, s'(b)) \end{aligned}$$

This section is called "strict", because $\pi_{\mathrm{fib}_f} \circ \overline{s'} \equiv \mathrm{id}_B$ and thus we can choose the homotopy to be $\mathsf{refl}$. We also note that $\Sigma\mathrm{fib}_f \simeq A$ due to Lemma 2.8.3.

Because of this observation, we only consider strict sections from now on, since they make our work significantly easier due to not needing the extra data of the homotopy.

**Definition 5.0.2:** Section of a family

Given a type family $B : A \to \mathcal{U}$, we call a dependent map $s : (a : A) \to B(a)$ a section of the type family $B$.

The definition of an initial section is quite simple. It is a section that only takes on initial values:

**Definition 5.0.3:** Initial Section

Given a type family $B : A \to \mathcal{U}$. We call a section $s : (a : A) \to B(a)$ *initial* iff the following holds:

$$\text{is-initial-section}(s) :\equiv (a : A) \to \text{is-initial}_{B(a)}(s(a))$$

To define dependent initial sections, we first unfold the definition for regular initial sections:

$$\text{is-initial-section}(s) \equiv (a : A, b : B(a)) \to \text{is-contr}\Big(\hom_{B(a)}(s(a), b)\Big)$$

For a dependent initial section, we want to require the dependent arrows to be contractible. Thus, we introduce the notion of "dependent initial object":

**Definition 5.0.4**

Given a type family $B : A \to \mathcal{U}$ and $b : B(a)$ with $a : A$ we define:

$$\text{is-dependent-initial}_B(b) :\equiv \begin{pmatrix} a' : A \\ f : \hom_A(a, a') \\ b' : B(a') \end{pmatrix} \to \text{is-contr}\Big(\text{dhom}_B^f(b, b')\Big)$$

**Definition 5.0.5:** Dependent Initial Section

We call a section a *dependent initial section* iff the following holds:

$$\text{is-dependent-initial-section}(s) :\equiv (a : A) \to \text{is-dependent-initial}_B(s(a))$$

Since dependent arrows generalize regular arrows, we can prove the following lemma easily:

**Lemma 5.0.6**

A dependent initial section is an initial section.

**Proof.** Given prf : is-dependent-initial-section$(s)$, we construct:

$$(a : A) \to \text{is-initial}_{B(a)}(s(a))$$
$$\lambda a.\ \lambda b.\ \text{prf}(a, a, \mathsf{id}_a, b)$$

By Lemma 3.3.8 we have $\text{dhom}_B^{\mathsf{id}_a}(s(a), b) \equiv \hom_{B(a)}(s(a), b)$ and thus the term is well-typed. □

## 5.1 Closure Properties

The closure properties of (dependent) initial sections are fairly simple to prove. They are closed under: Π-types, pullbacks, and composition. All proofs are easy to write down in `rzk`.

Since dependent initial sections and initial sections are only equivalent if the family is inner, we give the theorems and proofs for both notions.

**Theorem 5.1.1:** (Dependent) initial sections are closed under Π-types

Given a family of sections $s : (i : I, a : A(i)) \to B(i, a)$ over the index type $I : \mathcal{U}$ with bases $A : I \to \mathcal{U}$ and codomains $B : (i : I, a : A(i)) \to \mathcal{U}$. If all $s(i)$ are (dependent) initial sections, then the following map also is a (dependent) initial section:

$$(a : \Pi A) \to ((i : I) \to B(i, a(i)))$$
$$a \mapsto \lambda i.\ s(i, a(i))$$

**Proof.** We fix $a : \Pi A$ and $x : (i : I) \to B(i, a(i))$. We need to show that the following type is contractible:

$$\hom_{(i:I)\to B(i,a(i))}(\lambda i.\ s(i, a(i)), x)$$

We can flip the cube coordinate of hom with the $(i : I)$ argument and can equivalently show that this type is contractible:

$$(i : I) \to \hom_{B(i,a(i))}(s(i, a(i)), x(i))$$

By function extensionality, it suffices to show that each of the hom is contractible. But that is the case given that we have:

$$\text{is-initial}_{B(i,a(i))}(s(i,a(i)))$$

From the assumption that $s(i)$ is an initial section.

**Proof** for the dependent version:

We fix $x, y : (i : I) \to A(i)$ and $f : \hom_{(i:I)\to A(i)}(x, y)$ and $\overline{y} : (i : I) \to B(i, y(i))$. We also introduce the shorthand notation of $Q := \lambda a.\ (i : I) \to B(i, a(i))$ exhibiting dependent functions of $B$ as a family over dependent functions of $A$. We need to show that the following type is contractible:

$$\text{dhom}^{f}_{Q}(\lambda i.\ s(i, x(i)), \overline{y})$$

We can flip the cube coordinate of dhom with the $(i : I)$ argument and can equivalently show that this type is contractible:

$$(i : I) \to \text{dhom}^{\lambda t.\ f(t,i)}_{B(i)}(s(i, x(i)), \overline{y}(i))$$

By function extensionality, it suffices to show that each of the dhom is contractible. But that is the case given that we have:

$$\text{is-dependent-initial}_{B(i)}(s(i, x(i)))$$

From the assumption that $s(i)$ is a dependent initial section. □

**Theorem 5.1.2:** (Dependent) initial sections are closed under pullbacks

Given a (dependent) initial section $s : (a : A) \to B(a)$ and a map $k : A' \to A$ the pullback of $s$ along $k$

$$k \circ s : (a' : A') \to B(k(a'))$$

is a (dependent) initial section.

**Proof.** Trivial: we need to show that $s(k(a'))$ is initial in $B(k(a'))$ for any $a' : A'$, but that's obviously true, since $s$ is an initial section.

Analogous for the dependent version. □

For closure under composition, we need an auxiliary lemma about initial objects in $\Sigma$-types:

**Lemma 5.1.3:** Initial object in $\Sigma$-types

Given a type family $B : A \to \mathcal{U}$, an initial object $a : A$ and a dependent initial object $b : B(a)$ we have that $(a, b)$ is initial in $\Sigma B$.

**Proof.** To show that $\text{is-initial}_{\Sigma B}(a, b)$ holds, we have to show that the following type is contractible for any $(a', b') : \Sigma B$:

$$\hom_{\Sigma B}((a, b), (a', b'))$$

By Proposition 2.5.7 this type is equivalent to:

$$\Sigma\big(f : \hom_A(a, a'), \mathrm{dhom}_B^f(b, b')\big)$$

This type is contractible since both components are contractible: the first one since $a$ is initial and the second one since $b$ is dependent initial. □

**Theorem 5.1.4:** (Dependent) initial sections are closed under composition

Given two composable type families $B : A \to \mathcal{U}$ and $C : \Sigma B \to \mathcal{U}$ as well as sections of each $s : (a : A) \to B(a)$ and $s' : (x : \Sigma B) \to C(x)$. If $s$ is a (dependent) initial sections and $s'$ is a dependent initial section, then their composition is a (dependent) initial section:

$$\begin{aligned}(a : A) &\to \Sigma(b : B(a), C(a, b))\\ a &\mapsto (s(a), s'(a, s(a)))\end{aligned}$$

**Proof.** We need to show that

$$\text{is-initial}_{\Sigma(b:B(a),C(a,b))}(s(a), s'(a, s(a)))$$

holds. For this we use Lemma 5.1.3 and thus need to show that $s(a)$ is initial in $B(a)$ and that $s'(a, s(a))$ is dependent initial. Both of which are true by assumption.

**Proof** for the dependent version:

We fix $x, y : A$, $f : \hom_A(x, y)$ and $(b', c') : \Sigma(b : B(y), C(y, b))$. Now we need to show that the following type is contractible:

$$\mathrm{dhom}^f_{\lambda a.\ \Sigma(b:B(a),C(a,b))}((s(a), s'(a, s(a))), (b', c'))$$

That type is equivalent via Proposition 2.5.7 to:

$$\sum\Big(F : \mathrm{dhom}_B^f(s(a), b'), \mathrm{dhom}_C^{\lambda t.\ (f(t),F(t))}(s'(a, s(a)), c')\Big)$$

By assumption, $s'$ is a dependent initial section and thus the second components are contractible. Again by assumption, $s$ is also a dependent initial section, so the first component also is contractible, thus showing that the entire type is contractible. □

In the next subsection, we show that initial sections are dependent initial sections if the family is inner. With that result, we obtain that initial sections are closed under composition in an inner family:

**Corollary 5.1.5**

Given two composable type families $B : A \to \mathcal{U}$ and $C : \Sigma B \to \mathcal{U}$ where $C$ is inner and sections of each $s : (a : A) \to B(a)$ and $s' : (x : \Sigma B) \to C(x)$. If $s$ and $s'$ are initial section, then their composition is too:

$$\begin{aligned}(a : A) &\to \Sigma(b : B(a), C(a, b))\\ a &\mapsto (s(a), s'(a, s(a)))\end{aligned}$$

## 5.2 Initial Sections are Dependent Initial Sections

For the last theorem of this section, we turn to the very interesting result that an initial section already is a dependent initial section (given that the type family is inner). This result is at first surprising, since a dependent initial section is a much stronger statement. One way to explain this fact is that Π-types have something akin to continuity from topology. Their values are not allowed to vary arbitrarily, but only in a certain consistent manner. If we consider again the unfolded definition of an initial section, we obtain the type:

$$(a : A, b : B(a)) \to \text{is-contr}\Big(\hom_{B(a)}(s(a), b)\Big)$$

It is now these two input variables where we can observe this *continuity:* plugging in a dependent hom $\gamma$ for the second variable (and the hom $f$ that it is over in the first), we obtain:

$$(t : \Delta^1) \to \text{is-contr}\Big(\hom_{B(f(t))}(s(f(t)), \gamma(t))\Big)$$

If we now split this into the center and the homotopy of the contraction (by Proposition 2.5.7), we obtain:

$$\sum\left(\begin{array}{c} C : (t : \Delta^1) \to \hom_{B(f(t))}(s(f), \gamma(t)) \\ \Big(t : \Delta^1, F : \hom_{B(f(t))}(s(f), \gamma(t))\Big) \to C(t) = F \end{array}\right)$$

Inspecting the first component more closely, we observe that it is a dependent curried square:

$$(t : \Delta^1) \to \hom_{B(f(t))}(s(f), \gamma(t)) \simeq \mathsf{dcsquare}_B^{\lambda t.\ f(t)}\left(\begin{array}{ccc} s(x) & \xrightarrow{s(f)} & s(y) \\ \downarrow & & \downarrow \\ s(x) & \xrightarrow[\gamma]{} & \gamma(1) \end{array}\right)$$

The vertical maps are given by $C(t)$ and thus unique (by initiality of $s$). By the homotopy, we have that the left map must be $\mathsf{id}_x$ (since $\hom_{B(x)}(s(x), s(x))$ is contractible). And the right map is just some unique morphism. To finish the proof one just needs to derive that all $\gamma$ are equal to the composition of the upper two maps, which is the case when $B$ is an inner family (otherwise there would be no composition to talk about).

**Theorem 5.2.1**

If $B : A \to \mathcal{U}$ is an inner family, then an initial section $s : (a : A) \to B(a)$ is a dependent initial section.

**Proof.** We fix $x, y : A$ and $f : \hom_A(x, y)$ as well as $\overline{y} : B(y)$. Now, we need to show that $\text{dhom}_B^f(s(x), \overline{y})$ is contractible. To do this, we directly construct a center $c : \text{dhom}_B^f(s(x), \overline{y})$ and a contraction: $\Big(\overline{f} : \text{dhom}_B^f(s(x), \overline{y})\Big) \to \Big(c = \overline{f}\Big)$. We define $c$ to be the unique lift (since $B$ is inner) of:

$$D := \begin{array}{ccc} x & \xrightarrow{f} & y \\ & {\scriptstyle f}\searrow & \Big\downarrow {\scriptstyle \mathsf{id}_y} \\ & & y \end{array} \qquad \overline{D} := \begin{array}{ccc} s(x) & \xrightarrow{s(f)} & s(y) \\ & & \Big\downarrow {\scriptstyle \text{ctr}(\overline{y})} \\ & & \overline{y} \end{array}$$

$$D :\equiv \lambda(t, _)\colon \Delta^2.\ f(t) \qquad \overline{D} :\equiv [s(f), \mathrm{ctr}(\overline{y})]$$

Where $\mathrm{ctr}(\overline{y})$ is the center of the contraction of $\hom_{B(y)}(s(y), \overline{y})$ that we obtain from the fact that $s$ is an initial section.

To define the homotopy of our contraction of $\mathrm{dhom}_B^f(s(x), \overline{y})$, we fix an $\overline{f} : \mathrm{dhom}_B^f(s(x), \overline{y})$ and note the following:

$$\beta :\equiv \lambda t.\ \pi_1\left(\text{is-initial-section-s}\left(f(t), \overline{f}(t)\right)\right)$$
$$\beta : (t : \Delta^1) \to \hom_{B(f(t))}\left(s(f(t)), \overline{f}(t)\right)$$

This is a curried square with contractible vertical slices. Since $\overline{f}(0) \equiv s(x)$, we have $\beta(0) \equiv \mathsf{id}_{s(x)}$ thus we have:

$$\beta : \mathsf{dcsquare}_B^{\lambda(t,_).\ f(t)}\left(\begin{array}{ccc} s(x) & \xrightarrow{s(f)} & s(y) \\ {\scriptstyle \mathsf{id}_{s(x)}}\big\downarrow & & \big\downarrow{\scriptstyle \mathrm{ctr}(\overline{y})} \\ s(x) & \xrightarrow[\overline{f}]{} & \overline{y} \end{array}\right)$$

Applying Lemma 4.4.3 we thus obtain

$$\beta' : \mathrm{dhom}_B^{2,\lambda(t,_).\ f(t)}\left(\begin{array}{ccc} & s(y) & \\ {\scriptstyle s(f)}\nearrow & & \searrow{\scriptstyle \mathrm{ctr}(\overline{y})} \\ s(x) & \xrightarrow[\overline{f}]{} & \overline{y} \end{array}\right)$$

Which by Lemma 4.3.4 yields $c = \overline{f}$. □

## 5.3 Currying of Initial Sections

Currying also works for initial sections, but a particularly useful statement is that we can include the proofs of initiality in the values of the function. So instead of having this data:

$$\Sigma(s : (a : A) \to B(a), \text{is-initial-section}(s))$$

We can use Proposition 2.5.7 in reverse to obtain an equivalent type:

$$(a : A) \to \text{has-initial}(B(a))$$
$$\text{has-initial}(T : \mathcal{U}) :\equiv \Sigma(t : T, \text{is-initial}_T(t))$$

This view allows us to more easily support multiple parameters:

$$(a : A, b : B(a)) \to \text{has-initial}(C(a, b))$$

We need to uncurry it before we can apply Proposition 2.5.7, but ultimately it is equivalent to this:

$$\Sigma(s : ((a, b) : \Sigma B) \to C(a, b), \text{is-initial-section}(s))$$

We prefer the curried version that uses has-initial.

The same situation applies also for dependent initial sections. However, has-dependent-initial is a bit more complicated to define:

$$\Sigma(s : (a : A) \to B(a), \text{is-dependent-initial-section}(s)) \simeq (a : A) \to \text{has-dependent-initial}_a(B)$$

With $\text{has-dependent-initial}_t(R : T \to \mathcal{U}) :\equiv \Sigma(r : R(t), \text{is-dependent-initial}_R(r))$.

# 6 LARI Adjunctions

In this section, we define LARI adjunctions and show their equivalence to initial sections. LARI adjunctions are the backbone of LARI families, which generalize orthogonal families. Ultimately, we want to describe cocartesian families which are a special case of LARI families. From Lemma 4.1.2 we know that a map is orthogonal iff the Leibniz cotensor map is an equivalence. Weakening this requirement to only a LARI adjunction will later be the definition of a LARI family—hence we study LARI adjunctions now.

We follow the work from Buchholtz and Weinberger [5] with one important difference: we study LARI adjunctions through an alternative definition of initial sections which we introduced in the previous section. This is because LARI adjunctions are a relatively elaborate notion—a transposing adjunction plus an additional property. We mainly chose this way for improving the simplicity of our formalization work. A nice side effect of this effort is overall improving the readability of our proofs. Additionally, we proved several interesting results by taking this detour.

To prove the equivalence between LARI adjunctions and initial sections we use dependent initial sections also introduced in the previous section. Figure 4 shows which statements correspond to which theorems.

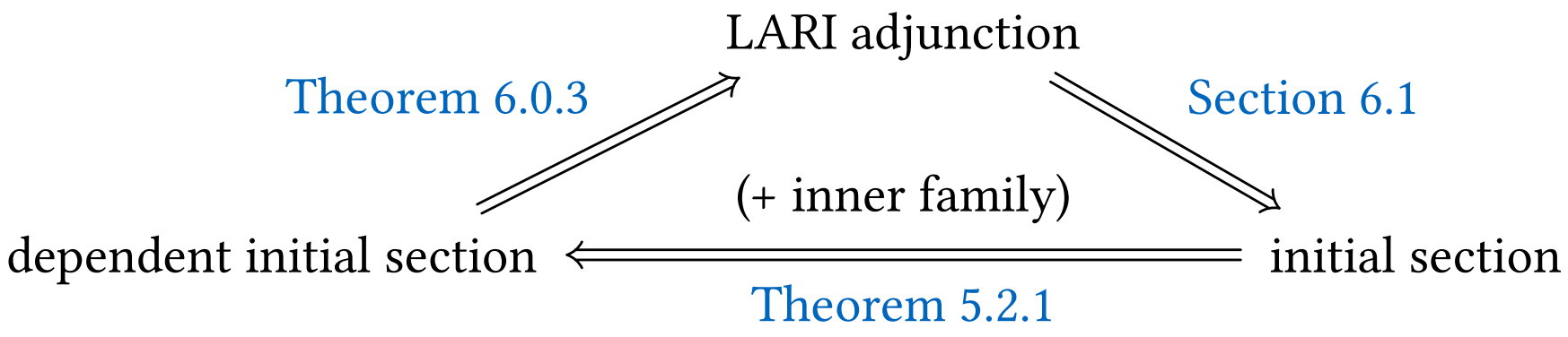

Figure 4: Proof structure for showing equivalence of initial sections and LARI adjunctions

LARI stands for *left adjoint right inverse;* a LARI adjunction is an adjunction where the right adjoint is an inverse of the left adjoint in a way that is compatible with the adjunction – in other words, the unit of the adjunction is an isomorphism.

**Definition 6.0.1:** LARI Adjunctions (cf. [5, Definition B.1.1])

We call a transposing adjunction $f : A \to B, u : B \to A$ with

$$\varphi : (a : A, b : B) \to \hom_B(f(a), b) \simeq \hom_A(a, u(b))$$

a *transposing left adjoint right inverse adjunction* (often shortened to *LARI adjunction*) iff

$$(a : A) \to \text{is-iso-arrow}_A\big(\eta_\varphi(a)\big)$$

holds (where $\eta_\varphi$ is the unit of the adjunction $\varphi$).

**Remark 6.0.2**

- Given a transposing LARI adjunction where the type $A$ is Rezk, we obtain a homotopy $H : \mathrm{id}_A \sim u \circ f$ witnessing the LARI nature of the adjunction.
- For a transposing adjunction it is only a *property* to be a transposing LARI adjunction, since is-iso-arrow is a proposition and Π-type construction preserves propositions.

The first direction that we're going to prove is the easiest. A dependent initial section trivially gives rise to a LARI adjunction.

**Theorem 6.0.3**

A dependent initial section $s$ gives rise to a LARI adjunction:

$$A \underset{\pi_B}{\overset{\overline{s}}{\underset{\bot}{\rightleftarrows}}} \Sigma B$$

where $\overline{s}(a) :\equiv (a, s(a))$.

**Proof.** We first show that $\overline{s}$ and $\pi_B$ form an adjunction: Given $x : A$ and $(y, \overline{y}) : \Sigma B$ we have:

$$\hom_{\Sigma B}(\overline{s}(a), (y, \overline{y})) \xrightarrow{\simeq} \left(f : \hom_A(a, y), \mathrm{dhom}^f_B(s(a), \overline{y})\right) \xrightarrow[\pi_1]{\simeq} \hom_A(a, y)$$

The first map is the equivalence from Lemma 3.3.9. The second map is an equivalence since $s$ is a dependent initial section and thus all these dependent homs are contractible. To see that this adjunction is LARI, it suffices to observe that by construction, $\mathsf{id}_{\overline{s}(a)}$ is mapped to $\mathsf{id}_a$ under this equivalence. □

Now we can combine this with the result from the previous section:

**Corollary 6.0.4**

An initial section $s : (a : A) \to B(a)$ of an inner family $B : A \to \mathcal{U}$ gives rise to a LARI adjunction:

$$A \underset{\pi_B}{\overset{\overline{s}}{\underset{\bot}{\rightleftarrows}}} \Sigma B$$

where $\overline{s}(a) :\equiv (a, s(a))$.

## 6.1 LARI Adjunctions are Initial Sections

Next we're going to show that a LARI adjunction gives rise to an initial section and that if a section has its canonical adjunction that it is an initial section.

**Theorem 6.1.1**

A transposing LARI adjunction from a Rezk type $A$ to a Segal type $B$, so

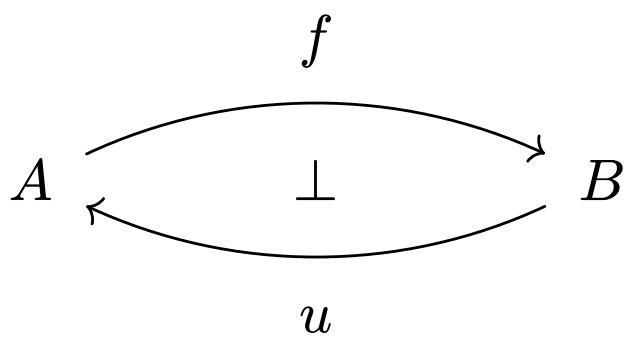

$$\varphi_{a,b} : \hom_B(f(a), b) \xrightarrow{\cong} \hom_A(a, u(b))$$
$$\text{is-LARI} : (a : A) \to \text{is-iso-arrow}\big(\eta_\varphi(a)\big)$$

Gives rise to an initial section:

$$\begin{aligned} \overline{f} : (a : A) &\to \text{fib}_u(a) \\ a &\mapsto \big(f(a), \eta_\varphi(a)^{-1} : u(f(a)) = a\big) \end{aligned}$$

**Proof.** Instead of directly showing that $\overline{f}$ is an initial section, we instead consider the composition from Figure 5. The squiggly arrows on the right side signify that the function takes on a value that is only *propositionally equal* to the one printed. For example, the first squiggly arrow includes the homotopy from $\text{rev}_A(\text{rev}_A(p)) = p$. Since we are interested in only showing that $\overline{f}$ is initial, we're only interested in this function up to equality (and the path that we would have to include would only distract from the statement).

We first note that this composite is initial, since $B$ is a Segal type and by Lemma 3.5.2 we have that $\big(f(a), \mathsf{id}_{f(a)}\big)$ is initial in the coslice. Additionally, we observe that the composition excluding $\overline{f}$ is a full embedding: by Lemma 3.6.3 all equivalences are full embeddings and Theorem 3.6.7 shows that $(\text{id}, \pi_1)$ is a full embedding. Since full embeddings compose by Lemma 3.6.4, we can apply Lemma 3.6.5 and conclude that $\overline{f}$ is an initial section. □

$$\begin{array}{cc}
(a : A) & a \\
\big\downarrow \overline{f} & \big\downarrow \\
\Sigma(b : B, u(b) = a) & \big(f(a), \eta_\varphi(a)^{-1} : u(f(a)) = a\big) \\
\simeq \big\downarrow (\text{id}, \text{rev}_A) & \rightsquigarrow \\
\Sigma(b : B, a = u(b)) & \big(f(a), \eta_\varphi(a) : a = u(f(a))\big) \\
\simeq \big\downarrow (\text{id}, \text{iso-eq}_A) & \rightsquigarrow \\
\Sigma(b : B, \text{Iso}_A(a, u(b))) & \big(f(a), \eta_\varphi(a) : a \cong u(f(a))\big) \\
\big\downarrow (\text{id}, \pi_1) & \big\downarrow \\
\Sigma(b : B, \hom_A(a, u(b))) & \big(f(a), \eta_\varphi(a)\big) \\
\simeq \big\downarrow \big(\text{id}, \varphi^{-1}_{a,f(a)}\big) & \rightsquigarrow \\
\Sigma(b : B, \hom_B(f(a), b)) & \big(f(a), \mathsf{id}_{f(a)}\big)
\end{array}$$

Figure 5: Composition for proof of Theorem 6.1.1.

An interesting and more common situation that we encounter is when the right adjoint is a type family projection $\pi_C$ and $B \equiv \Sigma C$ for some $C : A \to \mathcal{U}$. In this case, we can strictify $f$ to be a section and directly observe that the theorem above results in this section being initial. Even more useful, this section ends up being a dependent initial section.

**Proposition 6.1.2**

Given an inner type family $C : A \to \mathcal{U}$ over a Rezk type $A$, where the projection map $\pi_C$ is part of a LARI adjunction:

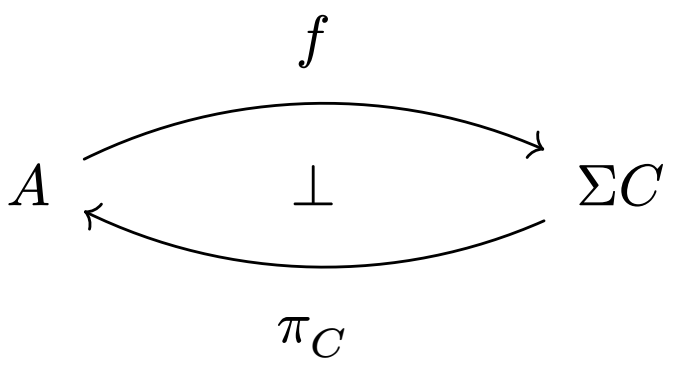

$$\varphi_{a,c} : \hom_{\Sigma C}(f(a), c) \xrightarrow{\simeq} \hom_A(a, \pi_C(c))$$

$$\text{is-LARI} : (a : A) \to \text{is-iso-arrow}\big(\eta_\varphi(a)\big)$$

Then there is a dependent initial section $g : \Pi C$.

**Proof.** We now transform $f$ into a strict section while keeping the adjunction. Since $\eta_\varphi(a)$ is an isomorphism and $A$ is Rezk, we get a path $p_a : a = \pi_C(f(a))$. Combining this with the decomposition $f(a) \equiv (\pi_C(f(a)), \pi_2(f(a)))$, we get:

$$p'_a : \left(a, \mathsf{tr}_C^{p_a^{-1}}(\pi_2(f(a)))\right) = (\pi_C(f(a)), \pi_2(f(a))) \equiv f(a)$$

We thus set $f'(a) :\equiv \left(a, \mathsf{tr}_C^{p_a^{-1}}(\pi_2(f(a)))\right)$ and by applying transport to $p'_a$, we get:

$$\varphi'_{a,c} : \hom_{\Sigma C}(f'(a), c) \xrightarrow{\simeq} \hom_A(a, \pi_C(c))$$

$$\varphi'_{a,c} :\equiv \varphi_{a,c} \circ \mathsf{tr}^{p'_a}_{\hom_{\Sigma C}(-,c)}$$

This adjunction also is LARI and in particular we have that

$$\eta_{\varphi'}(a) \equiv \varphi'_{a,f'(a)}\left(\mathsf{id}_{f'(a)}\right) \equiv \varphi_{a,f'(a)}\left(\mathsf{tr}^{p'_a}_{\hom_{\Sigma C}(-,c)}\left(\mathsf{id}_{f'(a)}\right)\right)$$

Since $\mathsf{tr}^{p'_a}_{\hom_{\Sigma C}(-,c)}$ is precomposition with $(p'_a)^{-1}$, we further have that:

$$\varphi_{a,f'(a)}\left(\mathsf{tr}^{p'_a}_{\hom_{\Sigma C}(-,c)}\left(\mathsf{id}_{f'(a)}\right)\right) \equiv \varphi_{a,f'(a)}\left((p'_a)^{-1}\right)$$

Since adjunctions are given by first applying the right adjoint to the hom-types and then precomposing with the unit, we have:

$$\varphi_{a,f'(a)}\left((p'_a)^{-1}\right) \equiv \eta_\varphi(a) \circ (\pi_C)_\#\left((p'_a)^{-1}\right) \equiv p_a \circ p_a^{-1} \equiv \mathsf{refl}_a$$

So the adjunction that we now have is a strict adjunction. We now apply Theorem 6.1.1 and thus obtain a $\overline{f'} : (a : A) \to \text{fib}_{\pi_C}(a)$ which by Lemma 2.8.2 and the calculation we just made results in the map $g :\equiv \overline{f'} : \Pi C$. □

**Remark 6.1.3**

In fact we could prove an even stronger statement using the following lemma from category theory:

**Lemma 6.1.4:** (cf. [27, Lemma 3.3.1])

Let $F : \mathcal{C} \leftrightarrows \mathcal{D} : G$ be an adjunction with $GF$ being an equivalence, then $F$ is fully faithful. Or equivalently, the unit $\eta : 1 \to GF$ is an isomorphism.

Since the equivalence of an adjunction decomposes into the following composition:

$$\hom_{\mathcal{D}}(Fx, y) \xrightarrow{G} \hom_{\mathcal{C}}(GFx, Gy) \xrightarrow{(\eta_x)_*} \hom_{\mathcal{C}}(x, Gy)$$

And we know that $\eta$ is an isomorphism, we also have that precomposing with $\eta_x$ is an isomorphism, from this it follows that $\hom_{\mathcal{D}}(Fx, y) \xrightarrow{G} \hom_{\mathcal{C}}(GFx, Gy)$ must also be an isomorphism.

In the concrete case where $\mathcal{D}$ is $\Sigma B$ for a family $B : A \to \mathcal{U}$ and $\mathcal{C}$ is $A$ with $G \equiv \pi_B$ and $F \equiv (\mathrm{id}, f)$ this means that the following map is an equivalence:

$$\hom_{\Sigma B}((a, f(a)), (a', b')) \xrightarrow{(\pi_B)_\#} \hom_A(a, a')$$

From this equivalence, we can deduce that $f$ is a dependent initial section, since the backwards map and homotopy provides us with a contraction of $\mathrm{dhom}^g_B(f(a), b')$.

# 7 LARI Families and Cocartesian Families

In this section we cover LARI cells, LARI families and cocartesian fibrations. We show that the notion of "having enough LARI lifts" from [7] is equivalent to our definition of LARI families as initial sections. We then show that having enough $i_0$-LARI lifts is equivalent to the notion of cocartesian family.

## 7.1 From Adjunction to Section

To enjoy the fruits of our labor from the two previous sections, we have to transform the notion of LARI-family given by U. Buchholtz and J. Weinberger [5, Definition 3.2.2] that uses LARI-adjunction into a definition that uses our notion of initial section. We begin with their definition, let $P : B \to \mathcal{U}$:

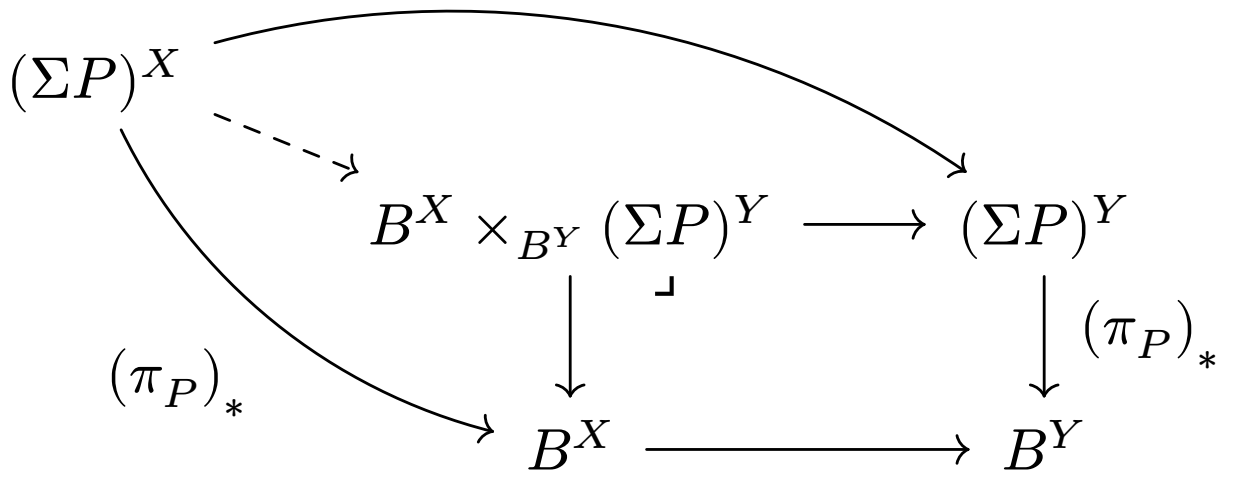

We want to obtain that the dashed map actually is the total-type projection of a family such that we can apply our initial section criterion to the map going back up. Such a map back along with the proof for it being initial is then the definition of a LARI family.

Applying Proposition 2.5.7, we obtain the following equivalences:

$$(\Sigma P)^Y \simeq \Sigma(g : Y \to B, \underbrace{(y : Y) \to P(g(y))}_{=: Q(g)}) \equiv \Sigma Q$$

$$(\Sigma P)^X \simeq \Sigma(g : X \to B, \underbrace{(x : X) \to P(g(x))}_{=: R(g)}) \equiv \Sigma R$$

Explicitly constructing the pullback and then inspecting it yields:

$$\begin{aligned}
& B^X \times_{B^Y} (\Sigma P)^Y \\
\simeq\ & \Sigma\Big(b : X \to B, p : Y \to \Sigma P, b \underset{Y \to B}{=} (\pi_P)_*(p)\Big) \\
\simeq\ & \Sigma\Big(b : X \to B, p_0 : Y \to B, p_1 : (y : Y) \to P(p_0(y)), b \underset{Y \to B}{=} p_0\Big) \\
\simeq\ & \Sigma(b : X \to B, \underbrace{(y : Y) \to P(b(y))}_{=: T(b)}) \equiv \Sigma T
\end{aligned}$$

Plugging all of these into the diagram gives us:

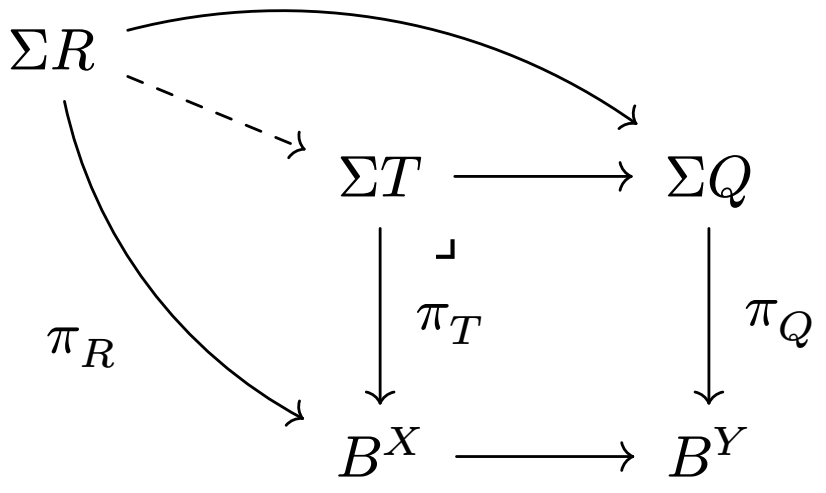

Now we can observe that we can express $R$ as a family over $\Sigma T$ and thus also ensure that the dashed map is a projection:

$$\begin{aligned}\Sigma R &\equiv \Sigma(g : X \to B, (x : X) \to P(g(x))) \\ &\simeq \sum\left(g : X \to B, p_0 : (y : Y) \to P(g(y)), \left\langle (x : X) \to P(g(x)) \;\middle|\; {}^{Y}_{p_0} \right\rangle\right) \\ &\simeq \sum\left((g, p_0) : \Sigma T, \left\langle (x : X) \to P(g(x)) \;\middle|\; {}^{Y}_{p_0} \right\rangle\right) =: \Sigma W\end{aligned}$$

Ultimately we obtain:

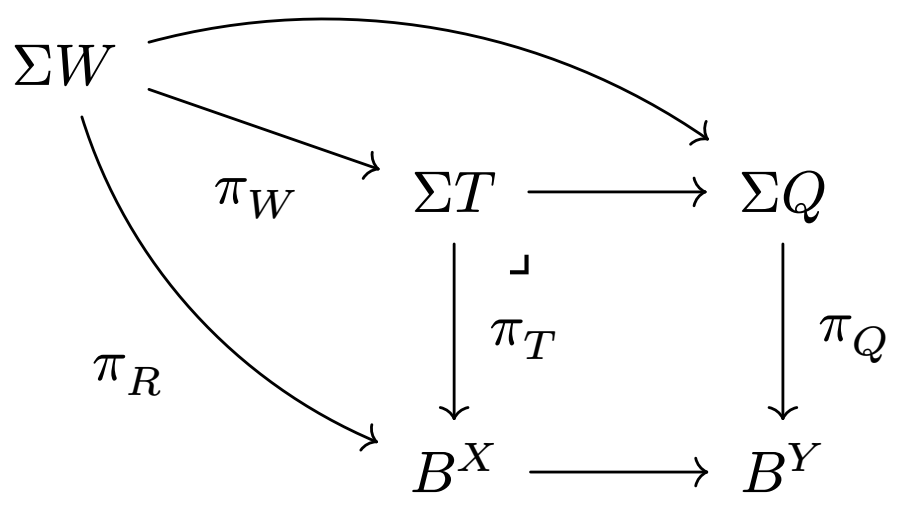

Where we can apply $\pi_R$ because $\Sigma W \simeq \Sigma R$. Finally, we can give our definition of LARI family using initial sections, by requiring $\pi_W$ to be a dependent initial section.

**Definition 7.1.1:** LARI Family

Given a shape inclusion $j : Y \subset X$, we say that $P : B \to \mathcal{U}$ is a $j$-LARI-family if the following holds:

$$\text{is-LARI-family}(P) := \begin{pmatrix} g : X \to B \\ f_0 : (y : Y) \to P(g(y)) \end{pmatrix} \to \text{has-dependent-initial}_{(g, f_0)}(Q)$$

$$\text{with } Q : \sum\begin{pmatrix} g : X \to B \\ (y : Y) \to P(g(y)) \end{pmatrix} \to \mathcal{U} \qquad (g, f_0) \mapsto \left\langle (x : X) \to P(g(x)) \;\middle|\; {}^{Y}_{f_0} \right\rangle$$

**Lemma 7.1.2**

Given a shape inclusion $j : Y \subset X$ and a family $P : B \to \mathcal{U}$. Then it is a $j$-LARI family if and only if the Leibniz cotensor map $\widehat{\pitchfork}_Y \pi_P$ is a $j$-LARI adjunction.

**Proof.** For the if part, we apply Proposition 6.1.2. For the only if part, we apply Theorem 6.0.3. □

## 7.2 Closure Properties

Since we defined LARI families as dependent initial sections, we can simply translate their closure from Section 5.1 to this concrete definition.

**Theorem 7.2.1:** LARI families are closed under Π-types (cf. [5, Proposition 3.2.5])

Given a shape inclusion $j : Y \subset X$ and given a family of families: $P : (i : I, b : B(i)) \to \mathcal{U}$ over bases $B : I \to \mathcal{U}$ where $P(i)$ is a $j$-LARI family. Then the family $P' : \Pi B \to \mathcal{U}$ with $P'(b) := (i : I) \to P(i, b(i)))$ also is a $j$-LARI family.

**Proof.** We apply Theorem 5.1.1 to the dependent initial sections that define LARI families. We have

$$\begin{aligned}
&\text{is-LARI-family}(P') \equiv \\
&\quad (g : X \to \Pi B, f_0 : (y : Y) \to P'(g(y))) \\
&\quad \to \text{has-dependent-initial}_{(g,f_0)}\left(\lambda(g, f_0).\ \left\langle (x : X) \to P'(g(x)) \;\middle|\; {}^{Y}_{f_0} \right\rangle\right) \\
&\equiv (g : X \to \Pi B, f_0 : (y : Y, i : I) \to P(i, g(y, i))) \\
&\quad \to \text{has-dependent-initial}_{(g,f_0)}\left(\lambda(g, f_0).\ \left\langle (x : X) \to (i : I) \to P(i, g(x, i)) \;\middle|\; {}^{Y}_{f_0} \right\rangle\right) \\
&\simeq (g : X \to \Pi B, f_0 : (y : Y, i : I) \to P(i, g(y, i)), i : I) \\
&\quad \to \text{has-dependent-initial}_{(g(-,i), f_0(-,i))}\left(\lambda(g', f_0').\ \left\langle (x : X) \to P(i, g'(x)) \;\middle|\; {}^{Y}_{f_0'} \right\rangle\right)
\end{aligned}$$

The last type is inhabited, because for every $i$ we have is-LARI-family$(P(i))$ which unfolded is:

$$\begin{aligned}
&(g : X \to B(i), f_0 : (y : Y) \to P(i, g(y))) \\
&\to \text{has-dependent-initial}_{(g,f_0)}\left(\lambda(g, f_0).\ \left\langle (x : X) \to P(i, g(x)) \;\middle|\; {}^{Y}_{f_0} \right\rangle\right)
\end{aligned}$$

□

**Theorem 7.2.2:** LARI families are closed under pullbacks (cf. [5, Theorem 3.2.8])

Given a shape inclusion $j : Y \subset X$ and a $j$-LARI family $P : B \to \mathcal{U}$. For any map $k : A \to B$, the pullback $k^*P : A \to \mathcal{U}$ also is a $j$-LARI family.

**Proof.** Unfolding the definition of is-LARI-family$_j(k^*P)$, we get:

$$\begin{pmatrix} g : X \to A \\ f_0 : (y : Y) \to P(k(g(y))) \end{pmatrix} \to \text{has-dependent-initial}_{(g,f_0)}(Q')$$
$$Q' :\equiv \lambda(g, f_0).\ \left\langle (x : X) \to P(k(g(x))) \;\middle|\; {}^{Y}_{f_0} \right\rangle$$

We compare this with is-LARI-family$_j(P)$ and note that the only difference is the appearance of $k$ and the use of $A$ instead of $B$. Now we can simply apply Theorem 5.1.2 to is-LARI-family$_j(P)$ which we pull back along the following map:

$$\begin{pmatrix} g : X \to A \\ f_0 : (y : Y) \to P(k(g(y))) \end{pmatrix} \to \begin{pmatrix} g : X \to B \\ f_0 : (y : Y) \to P(g(y)) \end{pmatrix}$$
$$(g, f_0) \mapsto (k \circ g, f_0)$$

□

**Theorem 7.2.3:** LARI families are closed under composition (cf. [5, Proposition 3.2.7])

Given a shape inclusion $j : Y \subset X$ and two composable $j$-LARI families $P : B \to \mathcal{U}$ and $R : \Sigma P \to \mathcal{U}$, then their composition $R \odot P : B \to \mathcal{U}$ is a $j$-LARI family.

**Proof.** Note that $(R \odot P)(b) \equiv \Sigma(p : P(b), R(b, p))$. We directly construct a value of type is-LARI-family$_j(R \odot P)$ using the assumptions. Let $\mathfrak{p}$ : is-LARI-family$_j(P)$ and $\mathfrak{r}$ : is-LARI-family$_j(R)$. We begin by fixing $g : X \to A$ and $f_0 : (y : Y) \to (R \odot P)(g(y))$. We now need to construct a value of type $\left\langle (x : X) \to (R \odot P)(g(x)) \;\middle|\; {}^{Y}_{f_0} \right\rangle$ and show that this is dependently initial in the family over $(g, f_0)$.

Consider $\mathfrak{p}_1(g, \pi_1 \circ f_0) : \left\langle (x : X) \to P(g(x)) \;\middle|\; {}^{Y}_{\pi_1 \circ f_0} \right\rangle$ which allows us to define $g' : X \to \Sigma P$ as $\lambda x.\ (g(x), \mathfrak{p}(g, \pi_1 \circ f_0, x))$. Now we obtain $\mathfrak{r}_1(g', \pi_2 \circ f_0) : \left\langle (x : X) \to R(g'(x)) \;\middle|\; {}^{Y}_{\pi_2 \circ f_0} \right\rangle$, which combines with $\mathfrak{p}$ to give us:

$$f := \lambda x.\ (\mathfrak{p}_1(g, \pi_1 \circ f_0, x), \mathfrak{r}_1(g', \pi_2 \circ f_0)) : \left\langle (x : X) \to (R \odot P)(g(x)) \;\middle|\; {}^{Y}_{f_0} \right\rangle$$

The last step is to show that $f$ is dependently initial, but that is the case, since $\mathfrak{p}_2$ and $\mathfrak{r}_2$ show that the ingredients are dependently initial. In particular we invoke Theorem 5.1.4. □

## 7.3 LARI Cells

To bridge the gap between our definition of LARI families and cocartesian fibrations, we use the notion of "having enough $j$-LARI lifts" studied by J. Weinberger [7]. This is because the usual definition of cocartesian fibrations is not as a $i_0$-LARI family (where $i_0 : \Delta^0 \to \Delta^1$ is the inclusion of the initial vertex). The definition that we use for $j$-LARI cells is not be the definition from Weinberger, instead we use one that he showed to be equivalent to his definition.

The definition itself is a mouthful and it is best understood using the graphic in Figure 6. There we consider a type family $P : B \to \mathcal{U}$ and the shape inclusion $\Delta^1 \subset \Lambda^2_1$ (we include it as the first morphism). We also have a $v : \Lambda^2_1 \to B$ in the base and a $g : (t : \Lambda^2_1) \to P(v(t))$. Note that we could specify this using a single function $\Lambda^2_1 \to \Sigma P$, but we need the data specified by $g$ separately. In the figure all of this data is marked in bold.

Now with all this data fixed, the following has to hold: for any other map $(w, m) : \Lambda^2_1 \to \Sigma P$ and morphisms in the base $\alpha : \hom_{\Lambda^2_1 \to B}(v, w)$ along with a morphism in the fiber that is restricted to $\Delta^1$ called $\beta$ (its type is rather complicated, so we leave this detail to the formal definition). We have that extending $\beta$ to $\Lambda^2_1$ in the fiber is unique (so the type of these extensions is contractible).

**Definition 7.3.1:** $j$-LARI cell (cf. [7, Definition 4.2, equation (3)])

Let $j : Y \subset X$ be a shape inclusion and $P : B \to \mathcal{U}$ a type family. We call a map $(v, g) : X \to \Sigma P$ a $j$-LARI cell if the following holds:

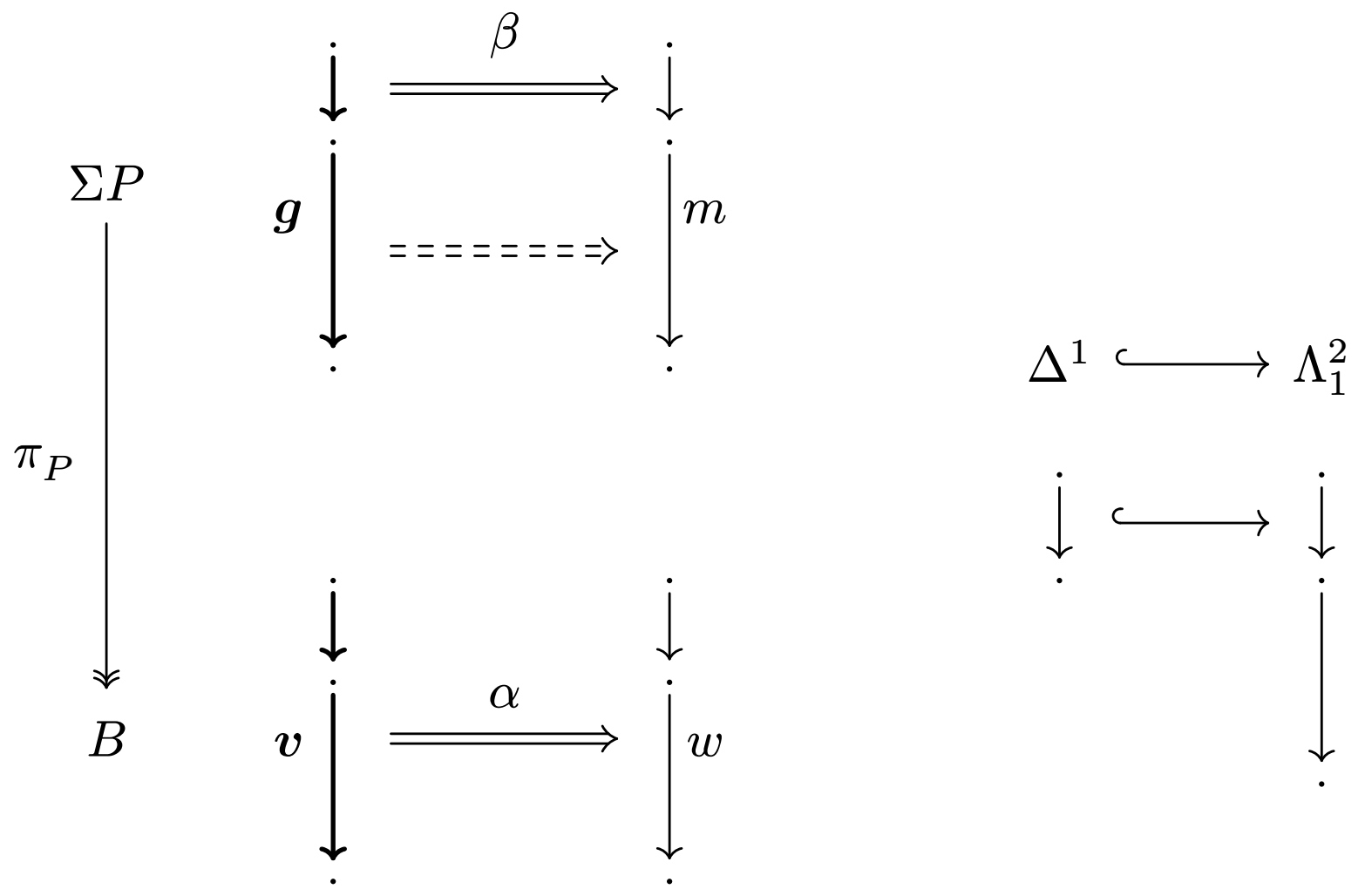

Figure 6: $j$-LARI cell condition with $j : \Delta^1 \subset \Lambda^2_1$. For the bold maps, we require the dashed map to be unique for any $m$, $w$, $\alpha$ and $\beta$.

$$\begin{aligned}\text{is-LARI-cell}_j(g) :\equiv \Big(& w : X \to B \\ & , m : (x : X) \to P(w(x)) \\ & , \alpha : \left\langle \Delta^1 \times X \to B \;\middle|\; {}^{\partial\Delta^1 \times X}_{[v,w]} \right\rangle \\ & , \beta : \left\langle ((t,y) : \Delta^1 \times Y) \to P(\alpha(t,y)) \;\middle|\; {}^{\partial\Delta^1 \times Y}_{[g,m]} \right\rangle\Big) \\ & \to \text{is-contr}\Big(\left\langle ((t,x) : \Delta^1 \times X) \to P(\alpha(t,x)) \;\middle|\; {}^{\partial\Delta^1 \times X \,\cup\, \Delta^1 \times Y}_{[g,m] \,\vee\, \beta} \right\rangle\Big)\end{aligned}$$

**Definition 7.3.2:** (cf. [7, Definition 4.3])

Let $j : Y \subset X$ be a shape inclusion and $P : B \to \mathcal{U}$ a type family. We say that $P$ *has enough $j$-LARI cells* if the following type is inhabited:

$$\text{has-enough-LARI-lifts}_j(P) :\equiv \begin{pmatrix} v : X \to B \\ f : (y : Y) \to P(v(y)) \end{pmatrix} \to \sum \begin{pmatrix} g : \left\langle (x : X) \to P(v(x)) \;\middle|\; {}^{Y}_{f} \right\rangle \\ \text{is-LARI-cell}_j(g) \end{pmatrix}$$

This definition is now equivalent with our earlier definition of is-LARI-family. The data is the same, it is just packaged differently. The proof takes quite some space to write down, but there are no interesting ideas involved.

**Theorem 7.3.3**

Let $j : Y \subset X$ be a shape inclusion and $P : B \to \mathcal{U}$ a type family. $P$ has enough $j$-LARI lifts if and only if $P$ is a $j$-LARI family.

**Proof.** We begin by unfolding the definitions and observing that they both start with a Π-type and the same two arguments. Thus, we fix both arguments $v : X \to B$ and $f : (y : Y) \to P(v(y))$ and then show that their outputs are equivalent. In particular, we need to show that the following two types are equivalent:

$$\text{has-dependent-initial}_{(v,f)}\Big(\lambda(v,f).\ \left\langle (x : X) \to P(v(x)) \;\middle|\; {}^{Y}_{f} \right\rangle\Big)$$

$$\sum\Big(g : \left\langle (x : X) \to P(v(x)) \;\middle|\; {}^{Y}_{f} \right\rangle, \text{is-LARI-cell}_j(g)\Big)$$

Unpacking has-dependent-initial, reveals that it also is a Σ-type with the same base and thus we need to show that for any $g : \left\langle (x : X) \to P(v(x)) \;\middle|\; {}^{Y}_{f} \right\rangle$ we have an equivalence between:

$$\text{is-dependent-initial}_{\lambda(v,f).\ \left\langle (x:X) \to P(v(x)) \middle| {}^{Y}_{f} \right\rangle}(g) \qquad\qquad \text{is-LARI-cell}_j(g)$$

Unfolding both definitions, we need to show that these two types are equivalent:

$$\begin{aligned}\Big(& (v', f') : \Sigma Q \\ & , F : \hom_{\Sigma Q}((v,f),(v',f')) \\ & , g' : \left\langle (x : X) \to P(v'(x)) \;\middle|\; {}^{Y}_{f'} \right\rangle\Big) \\ & \to \text{is-contr}\big(\text{dhom}^F_R(g,g')\big)\end{aligned} \quad \overset{!}{\simeq} \quad \begin{aligned}\Big(& w : X \to B \\ & , m : (x : X) \to P(w(x)) \\ & , \alpha : \left\langle \Delta^1 \times X \to B \;\middle|\; {}^{\partial\Delta^1 \times X}_{[v,w]} \right\rangle \\ & , \beta : \left\langle ((t,x) : \Delta^1 \times Y) \to P(\alpha(t,x)) \;\middle|\; {}^{\partial\Delta^1 \times Y}_{[g,m]} \right\rangle\Big) \\ & \to \text{is-contr}(\left\langle ((t,x) : \Delta^1 \times X) \to P(\alpha(t,x)) \;\middle|\; {}^{\partial\Delta^1 \times X \cup \Delta^1 \times Y}_{[g,m] \vee \beta} \right\rangle\end{aligned}$$

With the following shorthand notation:

$$
\begin{aligned}
Q &:= \lambda v' \colon X \to B.\ (y : Y) \to P(v'(y)) \\
R &:= \lambda (v, f) \colon \Sigma Q.\ \left\langle (x : X) \to P(v(x)) \,\middle|\, {}^{Y}_{f} \right\rangle
\end{aligned}
$$

Analyzing these two types, we can unify the following values:

- $v'$ and $w$,
- $g'$ and $m$ (and thus also $f'$ and $m$ on $Y$),
- $F$ and $(\alpha, \beta)$,

Finally consider the following chain of equalities and equivalences:

$$
\begin{aligned}
&\mathrm{dhom}^{F}_{R}(g, g') \\
\equiv\ & \left\langle (t : \Delta^1) \to R(F(t)) \,\middle|\, {}^{\partial\Delta^1}_{[g, g']} \right\rangle \\
\equiv\ & \left\langle (t : \Delta^1) \to \left\langle (x : X) \to P\big((\pi_1(F(t)))(x)\big) \,\middle|\, {}^{Y}_{\pi_2(F(t))} \right\rangle \,\middle|\, {}^{\partial\Delta^1}_{[g, g']} \right\rangle \\
\simeq\ & \left\langle (t : \Delta^1) \to \left\langle (x : X) \to P(\alpha(t, x)) \,\middle|\, {}^{Y}_{\beta(t, -)} \right\rangle \,\middle|\, {}^{\partial\Delta^1}_{[g, g']} \right\rangle \\
\simeq\ & \left\langle ((t, x) : \Delta^1 \times X) \to P(\alpha(t, x)) \,\middle|\, {}^{\partial\Delta^1 \times X \cup \Delta^1 \times Y}_{[g, m] \vee \beta} \right\rangle
\end{aligned}
$$

□

With this equivalence proven, we can finally define cocartesian families.

## 7.4 Cocartesian Families

As our final step to formalize cocartesian families, we have to show that having enough LARI lifts is equivalent to the usual condition of being cocartesian. The crux of this statement lies in our work about dependent squares from Section 4.

**Definition 7.4.1:** Cocartesian family

A type family $P : B \to \mathcal{U}$ is called *cocartesian* if it is an inner family and all arrows in the base have a *cocartesian lift* for any starting point. A cocartesian lift of an arrow $f : \mathrm{hom}_B(b, b')$ with a starting point $e : P(b)$ is an arrow $\overline{f} : \mathrm{dhom}_P(e, e')$ (where $e' : P(b')$ can be chosen freely by the morphism) with the following universal property:

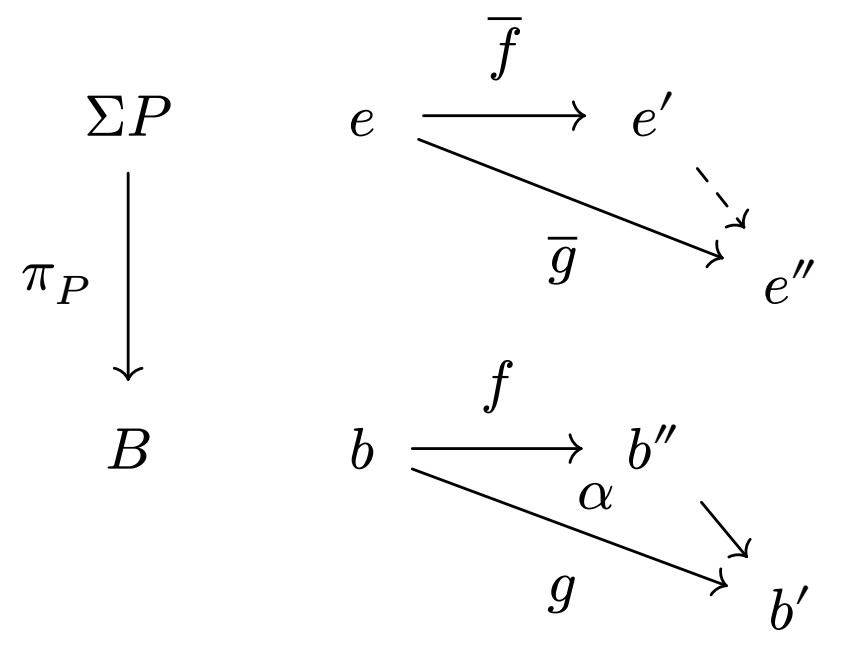

So for any triangle in the base and any lift of the $g$ edge to $\overline{g}$, the type of fillers in the fibers above should be contractible.

To define this notion formally, we begin with defining cocartesian lifts. Let $b, b' : B$ and $e : P(b), e' : P(b')$ as well as $f : \hom_B(b, b')$ and $\overline{f} : \mathrm{dhom}^f_P(e, e')$. We call $\overline{f}$ a cocartesian arrow if the following holds:

$$\text{is-cocartesian-arrow}\left(\overline{f}\right) := \left(\begin{array}{c} b'' : B \\ \tau : \hom^2_B\left(\begin{array}{ccc} & f \nearrow b' \searrow & \\ b & \longrightarrow & b'' \end{array}\right) \\ e'' : P(b'') \\ \overline{g} : \mathrm{dhom}_P^{\lambda t.\ \tau(t,t)}(e, e'') \end{array}\right) \to \text{is-contr}\left(\mathrm{dhom}^{2,\tau}_P\left(\begin{array}{ccc} & \overline{f} \nearrow e' \searrow & \\ e & \xrightarrow[\overline{g}]{} & e'' \end{array}\right)\right)$$

A cocartesian family then always has a cocartesian lift for any morphism in the base and a starting point in the fiber:

$$\text{is-cocartesian-fam}(P) := \text{is-inner-family}(P) \times \left(\left(\begin{array}{c} f : \Delta^1 \to B \\ e : P(f(0)) \end{array}\right) \to \sum\left(\begin{array}{c} \overline{f} : \langle (t : \Delta^1) \to P(f(t)) \mid {}^0_e \rangle \\ \text{is-cocartesian-arrow}\left(\overline{f}\right) \end{array}\right)\right)$$

In order to show that cocartesian families and $i_0$-LARI families are the same, we show that being a cocartesian arrow is the same as being a dependent initial object in a certain type family. While challenging in a technical and formal sense, this proof is just repackaging the data from one form into another.

**Lemma 7.4.2**

Given an inner type family $P : B \to \mathcal{U}$, we write $G := \Sigma(f : \Delta^1 \to B, P(f(0)))$ for the type of arrows in the base and a starting point lift. Also let $Q : G \to \mathcal{U}$ be defined as $\lambda(f, e).\ \langle (t : \Delta^1) \to P(f(t)) \mid {}^0_e \rangle$.

Now given $(f, e) : G$ and $F : Q(f, e)$, we show that the following two propositions are equivalent:

$$\text{is-cocartesian-arrow}(F) \simeq \text{is-dependent-initial}_Q(F)$$

**Proof.** We apply Lemma 2.5.11 to the uncurried versions of the unfolded definitions of the two propositions. We define $A, A' : U$ and $R : A \to \mathcal{U}, R' : A' \to \mathcal{U}$ such that $\text{is-cocartesian-arrow}(F) \equiv \Pi R$ and $\text{is-dependent-initial}_Q(F) \equiv \Pi R'$.

On the left side we have the following base type:

$$A := \sum\left(\begin{array}{c} b'' : B \\ \tau : \hom^2_B\left(\begin{array}{ccc} & f \nearrow b' \searrow & \\ b & \longrightarrow & b'' \end{array}\right) \\ e'' : P(b'') \\ \mathrm{dhom}_P^{\lambda t.\ \tau(t,t)}(e, e'') \end{array}\right)$$

Over which we form the following family:

$$R := \lambda(b'', \tau, e'', F').\ \text{is-contr}\left(\text{dhom}_P^{2;\tau}\left(\begin{array}{c} F \nearrow F(1) \searrow \\ e \xrightarrow[F']{} e'' \end{array}\right)\right)$$

On the right side, we have this base type:

$$A' := \sum \begin{pmatrix} (f', e') : G \\ m : \text{hom}_G((f, e), (f', e')) \\ Q'(f', e') \end{pmatrix}$$

Over which we form the following type:

$$R' := \lambda((f', e'), m, F').\ \text{is-contr}\left(\text{dhom}_Q^m(F, F')\right)$$

To apply Lemma 2.5.11, we need to define a function $\alpha_1 : A \to A'$ with $R'(\alpha_1(a)) \to R(a)$ and a function $\alpha_2 : A' \to A$ with $R(\alpha_2(a')) \to R'(a')$. To help intuition, we give the following visual representation of the two type families:

$$(a : A) \quad \to \quad R(a) \qquad\qquad (a' : A') \quad \to \quad R'(a')$$

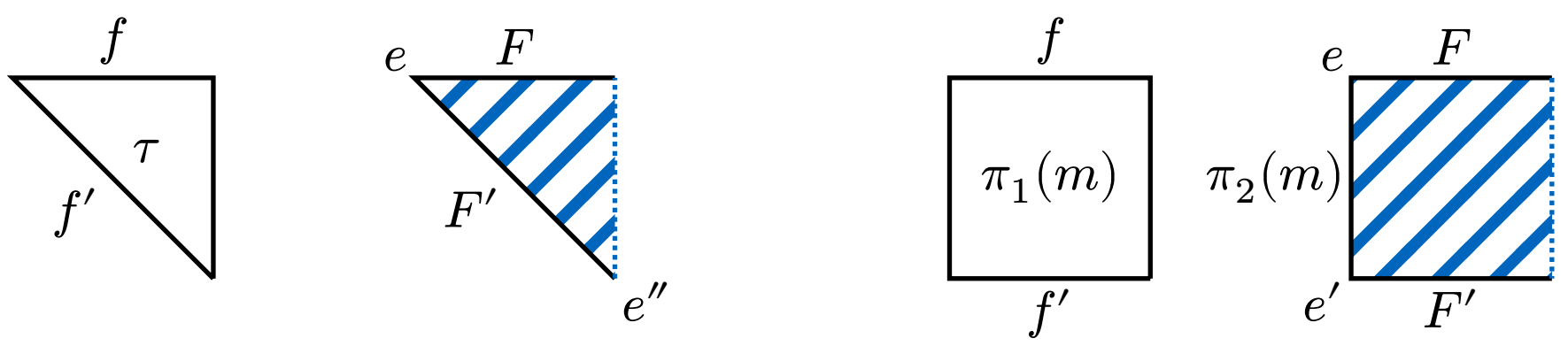

The blue parts correspond to the type that is required to be contractible while the black part corresponds to the base. We now give the first of the two base transformation maps.

$$\alpha_1 : A \to A' \quad \begin{pmatrix} b'' \\ \tau \\ e'' \\ F \end{pmatrix} \mapsto \begin{pmatrix} (\lambda t.\ \tau(t, t), e'') \\ \lambda s.\ (\lambda t.\ \sigma(t, s), e) \\ F' \end{pmatrix}$$

Where $\sigma$ is the square that has $\tau$ as the top right triangle and the degenerate triangle of $f'$ in the lower left. Formally we can write this as $\sigma := \lambda(t, s).\ \text{rec}_\vee^{\Delta^2, \nabla^2}(\tau(t, s), f'(t))$ with $\nabla^2$ being the lower right triangle of the square.

For the map back, we define:

$$\alpha_2 : A' \to A \quad \begin{pmatrix} (f', e') \\ m \\ F' \end{pmatrix} \mapsto \begin{pmatrix} f'(1) \\ \lambda(t, s).\ \pi_1(m(s))(t) \\ F(1) \\ \pi_2(m) \circ F' \end{pmatrix}$$

The composition of $\pi_2(m)$ and $F'$ is defined and unique, since the type family $P$ is inner.

To show show that $R'(\alpha_1(a)) \to R(a)$ holds for all $a : A$, we use Lemma 4.4.3, which shows that both types that ought to be contractible are equivalent, since we instantiated the square in $R'$ to have $\text{id}_e$ as the left side.

To show the same for $R(\alpha_2(a')) \to R'(a')$ for all $a' : A'$, we first destructure $a' \equiv ((f', e'), m, F')$ and then we observe that applying Lemma 4.4.2 we obtain the following equivalence for the type that ought to be contractible:

$$\mathrm{dhom}_Q^m(F, F') \simeq \sum \begin{pmatrix} d : \mathrm{dhom}_P^{\lambda t.\ \pi_2(m)(t,t)}(e, F'(1)) \\ \mathrm{dhom}_P^{2,\pi_2(m)}\begin{pmatrix} & F(1) & \\ {}^{F}\nearrow & & \searrow \\ e & \xrightarrow[d]{} & F'(1) \end{pmatrix} \\ \mathrm{dhom}_P^{2,\lambda(t,s).\ \pi_2(m)(s,t)}\begin{pmatrix} & e' & \\ {}^{\pi_2(m)}\nearrow & & \searrow^{F'} \\ e & \xrightarrow[d]{} & F'(1) \end{pmatrix} \end{pmatrix}$$

Since $P$ is inner, we have that the third type together with $d$ is contractible and thus we get:

$$\mathrm{dhom}_Q^m(F, F') \simeq \mathrm{dhom}_P^{2,\pi_2(m)}\begin{pmatrix} & F(1) & \\ {}^{F}\nearrow & & \searrow \\ e & \xrightarrow[\pi_2(m)]{} & F'(1) \end{pmatrix}$$

Which is contractible since $R(\alpha_2(a'))$ holds which by construction shows the contractibility of this type. □

**Theorem 7.4.3**

Given an inner type family $P : B \to \mathcal{U}$ the following two propositions are equivalent:

$$\text{is-cocartesian-fam}(P) \simeq \text{is-LARI-family}_{\{0\} \subset \Delta^1}(P)$$

**Proof.** We first unfold both definitions, directly unfolding has-dependent-initial in the process:

$$\text{is-cocartesian-fam}(P) \equiv \begin{pmatrix} f : \Delta^1 \to B \\ e : P(f(0)) \end{pmatrix} \to \sum \begin{pmatrix} \overline{f} : \langle (t : \Delta^1) \to P(f(t)) \mid {}^{0}_{e} \rangle \\ \text{is-cocartesian-arrow}(\overline{f}) \end{pmatrix}$$

$$\text{is-LARI-family}_{\{0\} \subset \Delta^1}(P) \equiv \begin{pmatrix} g : \Delta^1 \to B \\ f_0 : (y : \{0\}) \to P(g(y)) \end{pmatrix} \to \sum \begin{pmatrix} F : \langle (t : \Delta^1) \to P(f(t)) \mid {}^{0}_{e} \rangle \\ \text{is-dependent-initial}_Q(F) \end{pmatrix}$$

Where $Q$ is the family from the definition of Lemma 7.4.2. We can unify $g$ with $f$ and $f_0$ with $e$, since a map out of $\{0\}$ is equivalent to the codomain. We thus fix a $f : \Delta^1 \to B$ and an $e : P(g(0))$ and are left with showing that the following two types are equivalent:

$$\sum \begin{pmatrix} F : \langle (t : \Delta^1) \to P(f(t)) \mid {}^{0}_{e} \rangle \\ \text{is-cocartesian-arrow}(F) \end{pmatrix} \overset{!}{\simeq} \sum \begin{pmatrix} F : \langle (t : \Delta^1) \to P(f(t)) \mid {}^{0}_{e} \rangle \\ \text{is-dependent-initial}_Q(F) \end{pmatrix}$$

Here we again can unify the first component of the $\Sigma$-type and then apply Lemma 7.4.2 to finish the proof. □

From our abstract work in the previous sections we have that these cocartesian families satisfy the following closure properties:

**Theorem 7.4.4:** Closure properties of cocartesian families

- Cocartesian families are closed under $\Pi$-types:

Given a family of families: $P : (i : I, b : B(i)) \to \mathcal{U}$ over bases $B : I \to \mathcal{U}$ where each $P(i)$ is a cocartesian family. Then the family $P' : \Pi B \to \mathcal{U}$ with $P'(b) :\equiv (i : I) \to P(i, b(i)))$ also is a cocartesian family.

- Cocartesian families are closed under pullbacks:

  Given a cocartesian family $P : B \to \mathcal{U}$ and any map $k : A \to B$. The pullback $k^* P : A \to \mathcal{U}$ also is a cocartesian family.

- Cocartesian families are closed under composition:

  Given two composable cocartesian families $P : B \to \mathcal{U}$ and $R : \Sigma P \to \mathcal{U}$, then their composition $R \odot P : B \to \mathcal{U}$ is a cocartesian family.

# 8 Formalization with rzk

We already covered many aspects of `rzk` in Section 2 and Section 3 while introducing our theory. In this section we give a short overview of how one formalizes statements in `rzk` and then describe our experience of formalizing this thesis. Additionally, we list all statements that we have proven in `rzk` and report on the soundness of `rzk` itself, since we discovered a bug during our formalization work. This bug was promptly fixed by the maintainer after our report and the fix is included in `rzk` version `0.7.7`.

`rzk` always requires one to declare which version of the language one is using with the construct `#lang <version>`. At the moment the only available version is `rzk-1` (this is different from the version of the command-line tool and any extensions). Creating a new one in the future will allow projects to gradually change their existing files to take advantage of any compatibility-breaking changes.

## 8.1 How to Formalize with rzk

Statements in `rzk` are declared using the `#def` keyword after which one writes a name and then an arbitrary number of context elements `(X : Y)` where `X` is a variable that is bound to the type `Y` from that point onwards. After the context, one writes a `:` and then the type of the definition. Following that is a `:=` and the body of the definition.

As a simple example, we give the statement for dependent function currying:

```
#def dependent-curry
  ( A : U)
  ( B : A → U)
  ( C : (a : A) → B a → U)
  ( f : ( p : Σ (a : A) , (B a)) → C (first p) (second p))
  : ( a : A) → (b : B a) → C a b
  := \ s a b → s (a , b)
```

The file extension is `.rzk.md` which gets treated as a Markdown file by most programs. This allows prose text to be placed between code blocks that explain the intuition or statement in natural language. Additionally, lots of formatting options exist in Markdown such as including images and LaTeX. Proof assistants usually have the downside that seemingly simple and concise mathematical statements can take up vast amounts of space. For these cases it is very helpful to explain what is going on in prose before giving the mostly unreadable code.

### 8.1.1 Formatting Style

E. Riehl *et al.* [3] use a special formatting style that aims to make formalization easier to read and write. The most important aspect of this style guide is that expressions ought to be written down using their tree structure and using *indentation* to help guide the reader to discern which term belongs to which subexpression. An example the style guide gives for good style is the following formalization:

```
#def is-segal-is-local-horn-inclusion
  ( A : U)
  ( is-local-horn-inclusion-A : is-local-horn-inclusion A)
  : is-segal A
  :=
  \ x y z f g →
  contractible-fibers-is-equiv-projection
  ( Λ → A)
  ( \ k →
    Σ ( h : hom A (k (0₂ , 0₂)) (k (1₂ , 1₂)))
```

```
  , ( hom2 A
      ( k (0₂ , 0₂)) (k (1₂ , 0₂)) (k (1₂ , 1₂))
      ( \ t → k (t , 0₂))
      ( \ t → k (1₂ , t))
      ( h)))
( second
  ( equiv-comp
    ( Σ ( k : Λ → A)
      , Σ ( h : hom A (k (0₂ , 0₂)) (k (1₂ , 1₂)))
        , ( hom2 A
            ( k (0₂ , 0₂)) (k (1₂ , 0₂)) (k (1₂ , 1₂))
            ( \ t → k (t , 0₂))
            ( \ t → k (1₂ , t))
            ( h)))
    ( Δ² → A)
    ( Λ  → A)
    ( inv-equiv
      ( Δ² → A)
      ( Σ ( k : Λ → A)
        , Σ ( h : hom A (k (0₂ , 0₂)) (k (1₂ , 1₂)))
          , ( hom2 A
              ( k (0₂ , 0₂)) (k (1₂ , 0₂)) (k (1₂ , 1₂))
              ( \ t → k (t , 0₂))
              ( \ t → k (1₂ , t))
              ( h)))
      ( equiv-horn-restriction A))
    ( horn-restriction A , is-local-horn-inclusion-A)))
( horn A x y z f g)
```

The tree structure of the nodes are clearly represented in the indentation, subexpressions begin at the same indentation level and their own subexpressions start one more indentation further right. Parenthesis and commas match up due to the choice of using two tabs to indent content.

## 8.2 Experience of Using `rzk`

`rzk` is an amazing piece of software that we were very happy to use for this thesis. We wrote over 6 thousand lines of `rzk` and Markdown code. The experimental and early-stage status of `rzk` resulted in several issues that we had to overcome.

The first issue one encounters is that of cryptic error messages. The type checker prints the entire context and then a stacktrace of where the error occurred. The final part of such a stacktrace is often not where the actual problem lies and the very start also rarely is the culprit. It takes a lot of time and effort to get used to finding the exact location of the error. At the end we got accustomed with this and increased our speed of fixing issues by a lot.

The second issue is that `rzk` type-checking takes quite some time even for a few thousand lines of code. Running the check on a laptop can easily take 30 seconds or more.

The final big issue for us is the lack of certain features such as local variables, local definitions, implicit arguments and type inference. Additionally, the tooling could be improved, as the LSP did not work in our editor.

Formalizing mathematics in `rzk` works best when one already has a good understanding of how the proof works intuitively. However, we were also able to formalize statements without much knowledge of the theory itself. After working with `rzk` for a while, we felt that our intuition and manual type-checking had improved by a lot leading to less mistakes during formalization and faster discovery of proofs in the colloquial sense.

Compared with Isabelle/HOL [18] – the proof assistant that we had the most experience before working with `rzk` – `rzk` definitely is lacking in terms of ergonomics. Isabelle offers a number of different automated provers, which are able to prove quite advanced statements without any manual intervention. This stark difference to `rzk`'s purely type-theoretic proofs was a bit jarring in the beginning, but we adjusted rather quickly. Overall, we believe that tactics-based proofs simplify the implementation work immensely, but hinder the readability and will in many cases result in more difficult work down the road unless applied *correctly*. In Isabelle more fine-grained proofs using just introduction/elimination rules and logic is possible and offers a similar experience to the `rzk`.

## 8.3 Statements proven in rzk

The majority of statements that we gave in this thesis have been formalized in `rzk`. The exceptions are most of the examples, since `rzk` does not offer ways of defining inductive types without manual axioms, which is rather cumbersome. Concrete instantiations of the theory are also much more complicated in nature to formalize, so we refrained from attempting this in our limited time.

We built upon the work of E. Riehl *et al.* [3] which already contained many formalized statements from Section 2 and Section 3. Two important parts that were missing, were squares (Section 3.4) and full embeddings (Section 3.6). Many of the statements from Section 4 were already available as fibrations, but because we needed the family versions, we formalized those as well, except for their closure properties given in Theorem 4.2.1. Our main work from Section 5, Section 6 and Section 7 is fully formalized except for Lemma 6.1.4, which we did not use and only stated as a curiosity.

Our formalization is publicly available at https://github.com/BennoLossin/sHoTT. We note that this repository might still contain the working version of our formalization that contains unused and unfinished statements along with missing useful intermediate texts. We direct the reader to view any open or closed pull requests by us under the following links:

- Currently open pull requests at `sHoTT` (cleaned up code that is submitted to be included upstream),
- Closed pull requests at `sHoTT` (code that has been accepted upstream).

We hope that everything goes well and our formalization work is accepted upstream in the near future. Then they can be found in a more readable form under https://rzk-lang.github.io/sHoTT.

## 8.4 A Bug in rzk

We discovered a bug in `rzk` version `0.7.6` that allows us to prove the following:

```
#def issue
  ( A : U)
  ( x : A)
  : Equiv
    ( (t : Δ¹) → A [t ≡ 0₂ ↦ x])
    ( Δ¹ → A)
  :=
  ( identity ((t : Δ¹) → A [t ≡ 0₂ ↦ x])
  , is-equiv-identity ((t : Δ¹) → A [t ≡ 0₂ ↦ x]))
```

From this lemma we can derive that any inhabited type is contractible, which clearly is not the case.

The problem with this code is, that `rzk` type-checks it without any issues, but it should complain about the codomain of the identity function not matching. Outside of the body it is `Δ¹ → A` and inside it is

`(t : Δ¹) → A [t ≡ 0₂ ↦ x]`. The important difference is that the latter has a restriction part in the extension type.

After we reported this issue to the maintainer, he fixed it and released a new version. None of our or the existing formalization had used this buggy behavior, so all of those formalized statements were correct even after the update. Since forgetting the restriction part of extension types can easily happen accidentally, this is a serious issue. Our discovery is however not the only soundness issue in `rzk`.

### 8.4.1 Other Known Soundness Issues

Other than the issue that we found, `rzk` in its `rzk-1` version of the language contains multiple soundness holes. Fortunately, these are irrelevant for day to day formalization work. For this reason they have not yet been fixed, as opposed to our bug, which can be easily triggered by accidentally writing incomplete code.

#### `U : U`

The first soundness issue is that the judgment `U : U` holds. From it, one can exhibit a version of Russell's paradox. This can be fixed in a future version using the cumulative hierarchy strategy and does not impact normal mathematics as it is very difficult to *accidentally* create Russell's paradox.

#### Cubes depending on Types

In `rzk-1`, cubes are allowed to depend on types:

```
#define weird
    (A : U)
    (I : A → CUBE)
    (x y : A)
  : CUBE
  := I x × I y
```

[4] claims that this likely results in another inconsistency. Similarly to the previous hole, this one will not appear in actual proofs, since one has to explicitly use type-dependent cubes.

# 9 Conclusion and Future Work

In this thesis, we studied $\infty$-categories through a synthetic and type theoretic approach while formalizing the majority of statements using the experimental proof assistant `rzk` [4]. We introduced classical homotopy type theory as defined by The Univalent Foundations Program [1] and R. Harper *et al.* [19]. We followed with defining the simplicial type theory by E. Riehl and M. Shulman [2]. We then covered orthogonal families and other preliminary concepts for the rest of our work. As the main part of our work, we contributed an equivalence between initial sections and LARI adjunctions as defined by U. Buchholtz and J. Weinberger [5]. This equivalence enables us to define LARI-families in terms of initial sections, which simplifies the proofs for their closure properties immensely, especially for the purposes of formalization. As the final step to defining cocartesian families, we used the notion of "having enough LARI lifts" from J. Weinberger [7]. Our formalization work is based on top of the formalization by E. Riehl *et al.* [3] and we are upstreaming our work there.

### Future Work

There are several directions in which this work can be continued, both in small increments and in giant leaps. One incremental improvement would be the formalization of more concepts from category theory using `rzk`. For example, locally cocartesian fibrations and their connections to cocartesian fibrations, symmetric monoidial categories, or the analogous cancellation condition for LARI families with orthogonal fibrations. This amounts to showing the general statement of "Suppose that $g \circ f$ is cocartesian and $g$ is orthogonal, then $f$ is also cocartesian".

Many great improvements could also be performed on the `rzk` proof assistant. As we discussed in Section 8.2, there are quite a lot of flaws that, when addressed, could improve the usability of `rzk` immensely. Additionally, a version of directed univalence could be added to `rzk`.

A bigger undertaking would be the formalization of the alternative approach to adding directionality by D.-C. Cisinski, B. Cnossen, K. Nguyen, and T. Walde [20]. This would involve creating a novel experimental proof assistant similar to `rzk` that would not assume that $\Pi$-type formation is always allowed.

Since this field of study is still very young, there is always the option for a novel theory to disrupt the field. Although unlikely to achieve the "holy grail" of unifying type theory and higher categories like HoTT, there are still many things left to explore.

# A Glossary

| Syntax | Full Name | Introduced in |
|---|---|---|
| $x \equiv y$ | Judgmental equality | Section 2.1.4 |
| $\Sigma(a : A, B(a))$ | Sigma Type | Definition 2.1.8 |
| $(a : A) \to B(a)$ | Pi Type | Definition 2.1.11 |
| $g \circ f$ | Function or arrow composition | Definition 2.1.13, Definition 3.3.3 |
| $\Sigma B$ | Total type of the type family $B$ | Definition 2.1.17 |
| $\Pi B$ | Dependent functions or sections of the type family $B$ | Definition 2.1.17 |
| $\pi_B$ | Total type projection to the base of a type family $B$ | Definition 2.1.17 |
| $\mathrm{id}_A$ | Identity function of $A$ | Definition 2.1.18 |
| $A \simeq B$ | An equivalence between $A$ and $B$ | Definition 2.5.1 |
| $Q \odot P$ | Composition of two type families | Definition 2.8.6 |
| $\mathbf{1}, \mathbb{2}$ | Primitive cubes | Definition 3.1.1 |
| $\Delta^n, \partial\Delta^n$ | Simplicies and their boundary | Section 3.1.3 |
| $\partial\Delta^1$ | Simplicies and their boundary | Section 3.1.3 |
| $\psi \cup \zeta$ | Shape union | Definition 3.1.4 |
| $\psi \times \zeta$ | Type or Shape product | Section 2.1.8, Definition 3.1.5 |
| $\varphi \otimes \chi$ | Pushout product or subshape product | Definition 3.1.6 |
| $x \vee y$ | type-layer recursor for tope-disjunction | Definition 3.2.1 |
| $[x, y]$ | type-layer recursor for $\partial\Delta^1$ or $\Lambda^2_1$ | Definition 3.2.1 |
| $\langle (x : \psi) \to A \mid {}^{\varphi}_{a} \rangle$ | Extension type | Definition 3.2.2 |
| $\mathsf{id}_x$ | Identity Arrow | Definition 3.3.1 |
| $\eta_\varphi, \varepsilon_\varphi$ | Unit and counit of the adjunction $\varphi$ | Definition 3.7.2 |

# References

[1] The Univalent Foundations Program, *Homotopy Type Theory: Univalent Foundations of Mathematics*. Institute for Advanced Study, 2013. [Online]. Available: https://homotopytypetheory.org/book

[2] E. Riehl and M. Shulman, "A type theory for synthetic $\infty$-categories". [Online]. Available: https://arxiv.org/abs/1705.07442

[3] E. Riehl *et al.*, "Simplicial HoTT and synthetic $\infty$-categories." [Online]. Available: https://rzk-lang.github.io/sHoTT/

[4] N. Kudasov, A. Abounegm, and D. Danko, "Rzk proof assistant." [Online]. Available: https://rzk-lang.github.io/rzk/en/latest/

[5] U. Buchholtz and J. Weinberger, "Synthetic fibered $(\infty, 1)$-category theory". [Online]. Available: https://arxiv.org/abs/2105.01724

[6] C. Bardomiano Martínez, "Limits and colimits in synthetic $\infty$-categories", *Mathematical Structures in Computer Science*, vol. 35, 2025, doi: 10.1017/s0960129525100248.

[7] J. Weinberger, "Generalized Chevalley criteria in simplicial homotopy type theory." [Online]. Available: https://arxiv.org/abs/2403.08190

[8] A. N. Whitehead and B. Russel, *Principia Mathematica*, vol. 3. Cambridge University Press, 1910.

[9] A. Church, "A formulation of the simple theory of types," *Journal of Symbolic Logic*, vol. 5, no. 2, pp. 56–68, 1940, doi: 10.2307/2266170.

[10] P. Martin-Löf *et al.*, *Intuitionistic Type Theory*. Bibliopolis, 1984.

[11] V. Voevodsky, "Univalent Foundations Project," *a modified version of an NSF grant application*, pp. 1–12, Oct. 2010, [Online]. Available: http://www.math.ias.edu/vladimir/files/univalent_foundations_project.pdf

[12] L. de Moura, S. Kong, F. Doorn, and J. Raumer, "The Lean Theorem Prover (System Description)," 2015, pp. 378–388. doi: 10.1007/978-3-319-21401-6_26.

[13] L. d. Moura and S. Ullrich, "The Lean 4 Theorem Prover and Programming Language," in *Automated Deduction – CADE 28: 28th International Conference on Automated Deduction, Virtual Event, July 12–15, 2021, Proceedings*, Berlin, Heidelberg: Springer-Verlag, 2021, pp. 625–635. doi: 10.1007/978-3-030-79876-5_37.

[14] Lean Contributors, "Lean." [Online]. Available: http://lean-lang.org/

[15] T. Coquand and G. Huet, "The calculus of constructions," *Information and Computation*, vol. 76, no. 2, pp. 95–120, 1988, doi: https://doi.org/10.1016/0890-5401(88)90005-3.

[16] The Coq Development Team, "The Coq Proof Assistant." [Online]. Available: https://doi.org/10.5281/zenodo.14542673

[17] T. Nipkow and M. Wenzel, *Isabelle/HOL*. Springer Berlin, Heidelberg, 2002. doi: https://doi.org/10.1007/3-540-45949-9.

[18] Isabelle Contributors, "Isabelle." [Online]. Available: https://isabelle.in.tum.de/

[19] R. Harper *et al.*, "Lecture 15-819 Homotopy Type Theory." [Online]. Available: https://www.cs.cmu.edu/~rwh/courses/hott/

[20] D.-C. Cisinski, B. Cnossen, K. Nguyen, and T. Walde, “Synthetic Category Theory,” Nov. 2025. [Online]. Available: https://drive.google.com/file/d/1lKaq7watGGl3xvjqw9qHjm6SDPFJ2-0o/view. *Note: Accessed 27.11.2025*

[21] M. Land, *Introduction to Infinity-Categories*, 1st ed. in Compact Textbooks in Mathematics. Birkhäuser Cham, 2021. [Online]. Available: https://doi.org/10.1007/978-3-030-61524-6

[22] nLab authors, “dependent sum.” [Online]. Available: https://ncatlab.org/nlab/show/dependent+sum. *Note: Revision 19*

[23] nLab authors, “dependent product.” [Online]. Available: https://ncatlab.org/nlab/show/dependent+product. *Note: Revision 54*

[24] D. Gratzer, J. Weinberger, and U. Buchholtz, “Directed univalence in simplicial homotopy type theory.” [Online]. Available: https://arxiv.org/abs/2407.09146

[25] J. Sterling and C. Angiuli, “Normalization for Cubical Type Theory,” in *2021 36th Annual ACM/IEEE Symposium on Logic in Computer Science (LICS)*, IEEE, June 2021, pp. 1–15. doi: 10.1109/lics52264.2021.9470719.

[26] D. Gratzer, J. Sterling, C. Angiuli, T. Coquand, and L. Birkedal, “Controlling unfolding in type theory.” [Online]. Available: https://arxiv.org/abs/2210.05420

[27] S. Carmeli, T. M. Schlank, and L. Yanovski, “Ambidexterity and Height.” [Online]. Available: https://arxiv.org/abs/2007.13089